\documentclass[a4paper,english,10pt,twoside,article]{memoir}
\usepackage[latin1]{inputenc}
\usepackage[T1]{fontenc}
\usepackage{babel}
\usepackage{url}
\usepackage{amsmath,amssymb}
\usepackage{mathtools}
\usepackage[thmmarks]{ntheorem}
\usepackage{empheq}
\usepackage{graphicx}
\makeevenhead{plain}{\thepage}{}{}
\makeoddhead{plain}{}{}{\thepage}
\makeevenfoot{plain}{}{}{}
\makeoddfoot{plain}{}{}{}
\usepackage{lmodern} 
\usepackage[all]{xy} 
\SelectTips{cm}{} 
\usepackage{microtype} 

\makeatletter 
\newcommand*\wt[2][0.1ex]{%
        \begingroup
        \mathchoice{\wt@helper{#1}{#2}{\displaystyle}{\textfont}}
                   {\wt@helper{#1}{#2}{\textstyle}{\textfont}}
                   {\wt@helper{0.53ex}{#2}{\scriptstyle}{\scriptfont}}
                   {\wt@helper{#1}{#2}{\scriptscriptstyle}{\scriptscriptfont}}%
        \endgroup
        #2%
}
\newcommand*\wtb[2][0.1ex]{%
        \begingroup
        \mathchoice{\wt@helper{#1}{#2}{\displaystyle}{\textfont}}
                   {\wt@helper{#1}{#2}{\textstyle}{\textfont}}
                   {\wt@helper{0.34ex}{#2}{\scriptstyle}{\scriptfont}}
                   {\wt@helper{#1}{#2}{\scriptscriptstyle}{\scriptscriptfont}}%
        \endgroup
        #2%
}
\newcommand*\wt@helper[4]{%
        \def\currentfont{\the#41}%
        \def\currentskewchar{\char\the\skewchar\currentfont}%
        \setbox\tw@\hbox{\currentfont#2\currentskewchar}%
        \dimen@ii\wd\tw@
        \setbox\tw@\hbox{\currentfont#2{}\currentskewchar}%
        \advance\dimen@ii-\wd\tw@
        \rlap{\raisebox{-#1}{$\m@th#3\kern\dimen@ii\widetilde{\phantom{#2}}$}}%
}
\makeatother 

\DeclareMathOperator*{\colim}{{colim}}

\DeclareMathOperator{\ev}{Ev}
\DeclareMathOperator{\dist}{dist}

\DeclareMathOperator{\Sp}{Sp}

\DeclareMathOperator{\id}{Id}

\DeclareMathOperator{\co}{co}

\newcommand{\sm}{\textrm{-}}
\newcommand{\tpd}[1]{\tfrac{\partial}{\partial #1}}
\newcommand{\pd}[1]{\frac{\partial}{\partial #1}}

\newcommand{\epsin}{{\delta_0}}
\newcommand{\Or}{\mathrm{O}}

\newcommand{\U}{\mathrm{U}}

\renewcommand{\epsilon}{\varepsilon}

\newcommand{\ua}{{a_{\textrm{v}}}}
\newcommand{\ub}{{b_{\textrm{v}}}}

\newcommand{\uD}{{\underline{D}}}
\newcommand{\uC}{{\underline{C}}}

\newcommand{\tS}{{\wt{S}}}
\newcommand{\ti}{{\wtb{i}}}
\newcommand{\ts}{{\wt[0.2ex]{s}}}
\newcommand{\tp}{{\wt[0.2ex]{p}}}
\newcommand{\tX}{{\wt{X}}}
\newcommand{\tY}{{\wt{Y}}}
\newcommand{\tA}{{\widetilde{\alpha}}}

\newcommand{\tFL}{{
    \rlap{\raisebox{-0.3ex}{$\widetilde{\phantom{\FL}}$}}
    \FL}}
\newcommand{\ta}{{\wt[0.2ex]{a}}}
\newcommand{\tb}{{\wtb[0.3ex]{b}}}
\newcommand{\tC}{{\widetilde{C}}}
\newcommand{\uta}{{\ta_{\textrm{v}}}}
\newcommand{\utb}{{\tb_{\textrm{v}}}}

\newcommand{\oH}{{\overline{H}}}
\newcommand{\oS}{{\overline{S}}}

\newcommand{\oa}{{\overline{a}}}
\newcommand{\ob}{{\overline{b}}}

\newcommand{\te}{{\widetilde{\epsilon}}}

\newcommand{\HS}{\mathcal{H}}
\newcommand{\Com}{\mathcal{C}}

\newcommand{\As}{\mathcal{A}}
\newcommand{\Bs}{\mathcal{B}}
\newcommand{\GF}{\mathcal{GF}}
\newcommand{\FL}{\mathcal{FL}}

\newcommand{\R}{\mathbb{R}}
\newcommand{\N}{\mathbb{N}}
\newcommand{\F}{\mathbb{F}}
\newcommand{\Q}{\mathbb{Q}}
\newcommand{\C}{\mathbb{C}}

\newcommand{\Z}{\mathbb{Z}}

\renewcommand{\epsilon}{\varepsilon}
\newcommand{\arq}{{\overrightarrow{q}}}
\newcommand{\arp}{{\overrightarrow{p}}}

\newcommand{\arz}{{\overrightarrow{z}}}

\newcommand{\arzet}{{\overrightarrow{z_1}}}
\newcommand{\arzto}{{\overrightarrow{z_2}}}

\newcommand{\Jr}{[0,s_r]}
\newcommand{\Jinf}{[0,\infty[}

\DeclarePairedDelimiter\absv{\lvert}{\rvert}
\DeclarePairedDelimiter\norm{\lVert}{\rVert}

\DeclarePairedDelimiter\pare{\lparen}{\rparen}

\theoremseparator{.}

\newtheorem*{MainTheorem}{Main Theorem}

\newtheorem{Definition}{Definition}[chapter]
\newtheorem{Thm}[Definition]{Theorem}
\newtheorem{Proposition}[Definition]{Proposition}
\newtheorem{Lemma}[Definition]{Lemma}
\newtheorem{Corollary}[Definition]{Corollary}
\theoremsymbol{$\ast$}
\theorembodyfont{}
\newtheorem{Remark}[Definition]{Remark}

\theoremsymbol{$\square$}
\newtheorem*{bevis}{Proof}
\newtheorem*{bevisE1}{Proof of Proposition~\ref{prop:main}}

\newtheorem*{bevis:lem:11}{Proof of Lemma \ref{lem:11}}
\newtheorem*{mainbevis}{Proof of Main Theorem}
\theoremsymbol{``$\square$''}

\newcommand{\bR}{{\mathbb R}}
\newcommand{\bZ}{{\mathbb Z}}

\newcommand{\TN}{T^{*}N}
\newcommand{\Tn}{T^{*}_{q} N}
\newcommand{\TtN}{T^{*}N'}
\newcommand{\Ttn}{T^{*}_{\tq} \tN}
\newcommand{\tL}{L'}
\newcommand{\tN}{N'}
\newcommand{\tq}{q'}
\newcommand{\tx}{x'}
\def\co{\colon\thinspace}

\begin{document}

\pagestyle{plain}

\hyphenation{pa-ra-me-tri-zed}

$\strut$ \\

$\strut$ \\

{\center

 {\Huge \noindent Parametrized Ring-Spectra and the \\
 \noindent Nearby Lagrangian Conjecture}

$\strut$ \\

 Thomas Kragh\footnote{Supported by Carlsberg-fondet, MIT, Topology in
   Norway, and The University of Oslo.}\\
}
$\strut$ \\

$\strut$ \\

\noindent \textbf{Abstract:} We prove that any closed connected exact
Lagrangian manifold $L$ in a connected cotangent bundle $T^*N$ is up
to a finite covering space lift a homology equivalence. We prove this
by constructing a fibrant parametrized family of ring spectra $\FL$
parametrized by the manifold $N$. The homology of $\FL$ will be
(twisted) symplectic cohomology of $T^*L$. The fibrancy property will
imply that there is a Serre spectral sequence converging to the
homology of $\FL$ and the product combined with intersection product
on $N$ induces a product on this spectral sequence. This product
structure and its relation to the intersection product on $L$ is then
used to obtain the result. Combining this result with work of Abouzaid
we arrive at the conclusion that $L\to N$ is always a homotopy
equivalence.

\setlength\cftparskip{-15pt}
\tableofcontents

\chapter{Introduction}

Let $N$ be any closed smooth connected manifold, and let $j\colon L
\subset T^*N$ be a connected closed exact Lagrangian sub-manifold. Let
$\pi \colon T^*N \to N$ be the obvious projection and define $p =
\pi \circ j \colon L \to N$. The Nearby Lagrangian conjecture states
that $L$ is Hamiltonian isotopic to the zero section. In \cite{FSS},
which introduced new powerful methods, there is a good summary of what
was known at that time about this conjecture.

An easy to prove property of exact Lagrangians is that if
$f\colon N' \to N$ is a smooth covering space then the associated pull
back 
\begin{align*}
  \xymatrix{
     L' \ar[r]^{j'} \ar[d] & T^*N' \ar[d]^{T^*f} \\
     L \ar[r]^j & T^*N,
  }
\end{align*}
defines an exact Lagrangian $j'\colon L'\subset T^*N'$ - if $L'$ is
not connected, we may make a choice of connected component. Combining
this with the result in \cite{MR1090163} that $p$ is always surjective
(which also works when $N$ is not closed) one easily gets that
$\pi_1(N)/p_*(\pi_1(L))$ is finite. Note that this surjectivity result
works when $L$ is assumed closed but not $N$, which means that this
possibility is in fact excluded. So the assumption that $N$ is closed 
is in fact superfluous, but we include it so as not to cause
confusion.

\begin{MainTheorem}
  If in addition $N$ is oriented and the induced map $p_*$ on
  fundamental groups is surjective then the induced map $p_*$ on
  homology is an isomorphism.
\end{MainTheorem}

To see just how restrictive this is in general we prove a conjecture
by Arnold on the degree of $p$.

\begin{Corollary}
  \label{cor:3}
  The degree of $p$ is always non-zero.
\end{Corollary}

\begin{bevis}
  Let $N$ be any (possibly non-orientable) smooth connected manifold,
  and let $L\subset T^*N$ be as above. By lifting to the \emph{finite}
  covering of $N$ associated to the image $p_*(\pi_1(L))$ intersected
  with the kernel of the first Stiefel-Whitney class of $N$ we get a
  lifted exact Lagrangian $j'\colon L'\subset T^*N'$. The actual lift
  could have 2 connected components, but as above we let $L'$ be one of
  these components. By the definition of degree (see \cite{MR0192475})
  we have that $p'=\pi' \circ j' \colon L' \to N'$ has non-zero degree
  if and only if $p$ has non-zero degree.

  By construction this lift satisfies the assumptions of the Main
  Theorem and so $p'$ has degree one.
\end{bevis}

We in fact see that the Main Theorem is much stronger because it tells
us that up to a finite lift (and restricting to a component) the map
$p$ is a homology equivalence - and on all further finite lifts it is
also a homology equivalence (this is not true for any homology
equivalence so this strengthens the statement).

In the case of vanishing Maslov class and $N$ simply connected
Fukaya, Seidel and Smith proved in \cite{FSS} that $p$ induces a
homology isomorphism. Later in \cite{Abou1} Abouzaid removed the
assumption on $N$ and strengthened this to a homotopy equivalence. 
However, the assumption on the Maslov class remained. Only a few
and not very general results was known about the non-vanishing case. However, in
Appendix~\ref{cha:nlahe} Abouzaid uses Corollary~\ref{cor:1gfr} below
to prove that in fact this assumption is always true. So as
Theorem~\ref{thm:main} states we now know that $L\to N$ is a homotopy
equivalence.

\begin{Corollary}
  \label{cor:1gfr}
  Lifting $L \to T^*N$ to the universal cover $N' \to N$ the Maslov
  class of a component $L' \to T^*N'$ vanishes.
\end{Corollary}

Note that the universal cover of $N$ may not be closed and thus $L'$ 
may not be closed either.

\begin{bevis}
  Since we have established that some partial lift $L'' \to T^*N''$ is a
  homology equivalence we have that in the diagram
  \begin{align*}
    \xymatrix{
      \pi_1(L'') \ar[r]^{p_*} \ar[d] & \pi_1(N'') \ar[d] \\
      H_1(L'') \ar[r]^{p_*}_{\cong} \ar[d]^{m_L} & H_1(N'') \ar@{.>}[dl] \\
      \Z.
    }
  \end{align*}
  the dashed arrow exists making it commutative. Here the top vertical
  arrows are the Abelianization maps, and $m_L$ is given by evaluating
  the Maslov class. This proves that $m_L$ vanishes on the kernel of
  $p_* \colon \pi_1(L'') \to \pi_1(N'')$. So since the Maslov index of
  any loop in the cover $L'$ is the same as the Maslov index of its
  projection to $L$ the corollary follows.
\end{bevis}

The proof of the Main Theorem uses parts of Viterbos Transfer map
construction in \cite{MR1617648}, which we constructed as a map of
spectra in \cite{hejeh}. The idea is that the target spectrum, which
represents (twisted)\footnote{(twisted) is explained in
  Remark~\ref{rem:coef} and Appendix~\ref{cha:coor}} symplectic
cohomology of $T^*L$, can be created as a fibered spectrum $\FL$ over
the base $N$ and in fact turns out to be fibrant\footnote{Our
  notion of Fibrant is defined in Section~\ref{cha:fib-fl}, but it is
  close to Serre fibrant}. This means that the homology of
the fibers form a local system. It also implies the existence of a
Serre spectral sequence with page 2 given by the homology of $N$ with
coefficients in the homology of the fibers, strongly converging to the
homology of $\FL$ - i.e. the (twisted) symplectic cohomology of
$T^*L$. There are also natural product structures on the fibers
(continuous over the base) and when combined with the intersection
product on $N$ we get a product on the spectral sequence. The idea
that such a spectral sequence might exist (and be important) was
inspired by a description of a Spectral sequence in \cite{FSS} by
Fukaya, Seidel and Smith.

To construct this fibered spectrum $\FL$ we start in
Section~\ref{finred} by introducing finite dimensional approximations
of the action integral $A$, which is defined on paths $\gamma \colon I
\to T^*N$ by
\begin{align} \label{action}
  A(\gamma) = \int_\gamma \lambda - Hdt.
\end{align}
Here $\lambda$ is the canonical Liouville 1-form given by
\begin{align*}
  \lambda_{q,p}(v) = p(\pi_*(v)), \qquad q\in N, p \in T^*_qN,
\end{align*}
and $H\colon T^*N\to \R$ is a Hamiltonian. We will restrict the class
of Hamiltonians that we will consider. These finite dimensional
approximations are versions of Chaperons broken geodesic approach in
\cite{MR765426}, which is similar to the generating functions used by
Viterbo in \cite{MR1617648}. So the finite dimensional
approximations will be functions defined on finite dimensional
manifolds depending on an $r \in \N$ with $r>>1$ depending on the
Hamiltonian $H$ and other structure. These functions will (when
properly handled) have well-defined Morse homologies, which  
is (twisted) Floer homology of $A$. To define symplectic cohomology
one really only need to consider closed loops in $T^*N$. However, 
the fiber over $q\in N$ of our spectrum is constructed by considering
paths starting and ending on the Lagrangian $T^*_qN$ - i.e. a
different Lagrangian boundary condition.

There are some difficulties concerning orientations, which we
summarize in the following remark. For concreteness (and because this
is what we compute with our spectrum) we describe this in terms
of symplectic cohomology. However, the same difference in
orientations is there on the Morse homology of the action and the
finite dimensional approximations. Indeed, symplectic cohomology is as
special case of the Floer (co)homology associated to a Hamiltonian
quadratic at infinity (see \cite{MR2190223}).

\begin{Remark}
  \label{rem:coef}
  There are some orientation discrepancies, which can be slightly
  confusing. We assume here that $N$ and $L$ are oriented. Until
  recently it was thought that this implied that the symplectic
  cohomology $SH^*(TN)$ was isomorphic to $H_*(\Lambda N)$. Indeed,
  this was established in \cite{MR1403954}, \cite{MR2190223}, and
  \cite{Salweb}. There is a slight error in this. The fact is that if
  the second Stiefel-Whitney class of $N$ does not vanish on $\pi_2(N)$
  this is false. Indeed, in \cite{hejeh} it was proven that Viterbo
  functoriality in the case of $L\subset T^*N$ can be realized as a
  map of Thom-spectra
  \begin{align}\label{eq:8sdafuc}
    (\Lambda N)^{TN} \to (\Lambda L)^{TL\oplus \eta},
  \end{align}
  where $(-)^\beta$ denotes Thom-spectrum construction, and $\eta$ is
  a virtual bundle, which is \emph{not oriented} unless the relative
  Stiefel Whitney class vanishes on $\pi_2(L)$ 
  (this last fact is proved in Corollary~\ref{cor:2} and not in
  \cite{hejeh}). For the reader mostly concerned with (co)homology this
  simply means that the homology of the target is not the homology of
  $\Lambda L$, but a twisted version of it (basically because the Thom
  isomorphism only works for oriented bundles). This is contradictory
  to previous results since Viterbo functoriality
  tells us that symplectic cohomology is natural with respect to
  restrictions (as Paul Seidel and Mohammed Abouzaid pointed out to
  us). Indeed, this bundle $\eta$ was thought to \emph{be} oriented
  due to the previous results. The only solution to this problem is
  that the transgression of the second Stiefel Whitney class enters in
  when defining signs on the differentials of the Floer complex. As a
  test Paul Seidel produced a calculation that $SH^*(T^*\C P^3)$ in
  fact vanishes with rational coefficients, and this can be explained
  by it being 2-torsion with $\Z$ coefficients, as one would not
  previously have expected, but with the correct signs this is the
  correct answer. These new signs in the differentials has further
  been verified by Abouzaid in Appendix A of \cite{abou2}.

  For the purpose of this paper we mention here that for general
  $L\subset T^*N$ there are three potentially different orientation
  choices that seems important in this discussion. We will for
  simplicity also assume that $L$ is orientable. The first two are
  only dependent on $L \subset T^*L$ and not how it is embedded into
  $T^*N$. As explained in Appendix~\ref{cha:coor} different finite
  dimensional approximations may give different results.
  \begin{itemize}
  \item {As noticed by Viterbo in \cite{MR1617648} the finite
      dimensional approximations for $L \subset T^*L$ the zero-section
      do in fact reproduce the loop space homology. I.e. no
      transgression of the second Stiefel-Whitney class.}
  \item {The symplectic cohomology of $T^*L$ reproduces the loop space
      homology twisted by the canonical transgression of the second
      Stiefel Whitney class to $\pi_1(\Lambda L)$ - as described by
      Abouzaid in \cite{abou2}.}
  \item {The finite dimensional approximation of the action provided
      by $DT^*L \subset T^*N$ gives a potentially third result, which
      is loop space homology twisted by the transgression of the
      \emph{relative} second Stiefel Whitney class of $L\to N$
      transgressed to $\pi_1(\Lambda L)$. This is the homology of the
      target spectrum in Equation~\eqref{eq:8sdafuc} above, and is
      also the homology of $\FL$ appearing in this paper, and the one
      used in the main argument.}
  \end{itemize}
  The reason for these differences is that the orientation in
  symplectic homology does not generally agree with the one defined by
  finite dimensional approximations, which for this purpose can be
  thought of as depending on a choice of Lagrangian foliation. This
  is explained in more details in Appendix~\ref{cha:coor}.
\end{Remark}

Section~\ref{cha:pseudo} and Section~\ref{cha:paracon} defines
parametrized pseudo-gradients and parametrized Conley index theory,
which can be thought of as a version of parametrized Morse theory on
non-compact manifolds, which do not suffer from the defect that one
actually needs the functions to be Morse. And for the finite
dimensional approximations we thus get parametrized Conley indices
over the base $N$. The fact that we do not need Morse means that the
fibers are defined even in the non-generic cases. 

In Section~\ref{cha:fib-susp} we relate the finite dimensional
approximations for different $r$'s. The relation can heuristically be
thought of as the idea that the fibers of the parametrized Conley
index for $r+1$ is given by a reduced suspension of the fibers of the
parametrized Conley index for $r$. This is not completely true because
the fibers may behave slightly bad due to the non-genericity that we
cannot avoid in a parametrized setting. However, in the end these
problems will disappear when considering $\FL$ because it is fibrant.

In Section~\ref{cha:gen-spec} we use the suspension result from
Section~\ref{cha:fib-susp} to construct fibered spectra over $N$ whose
homology (i.e. the homology of the total space) will be the
(twisted) Floer homology associated to the action $A$. The homologies
of the fibers will \emph{not} form local systems - i.e. these
fibered spectra will \emph{not} be \emph{fibrant}. In fact, we will
only see this special feature for $\FL$ because it is defined
carefully as a certain limit (se below). This fibrancy can in fact be
generalized to fibered spectra computing the (twisted) symplectic
cohomology of any compact exact sub-Liouville domain in $T^*N$, but we
will not consider this here.

In Section~\ref{theex} we define $\FL$ as a limit when $s\to \infty$
of the fibered spectra defined in Section~\ref{cha:gen-spec} given by a
family of Hamiltonians $H^s,s\in [0,\infty[$. This is a refinement of
the idea that symplectic cohomology of $DT^*L$ as a symplectic
neighborhood of $L$ can be defined by the limit of Floer homology for
some Hamiltonians that go to infinity on the complement of $DT^*L$ and
stays $0$ on $L$. To get all the structure we need - most importantly
fibrancy - we have to be very careful in defining these Hamiltonians,
and this is done in Appendix~\ref{cha:pd-ham}.

In Section~\ref{homoindex} we construct parallel transport in
parametrized Conley indices. We need this to properly define the
continuation maps increasing $s$ (i.e. the maps over which 
we took the limit) in Section~\ref{theex} as fibered maps (called
ex-maps in the paper) over $N$. We also need this at several points
later to prove fibrancy and construct other structures on $\FL$.

In Section~\ref{cha:fib-fl} and Section~\ref{serre} we use this
parallel transport to construct (stable) homotopy lifts in $\FL$ and
prove that the homology of the fibers form a local system and that we
have an associated Serre spectral sequence as mentioned above.

Section~\ref{cha:inc-const} through Section~\ref{cha:serprod} concerns
structures such as products and a careful analysis of the constant
loops on $L$, which is then related to the spectral sequence. The
product is a generalization of the Chas-Sullivan product (and thus
probably related to the pair of pants products). This product on the
spectral sequence can be considered a generalization of the product
constructed in \cite{MR2039760}. In fact, when $L=N$ is the zero
section we recover their spectral sequence.

All these structures are needed for the main argument, which is
carried out in Section~\ref{cha:proof-theor-refh}. Since this is great
motivation for all the structure and is readable without understanding
much of the paper we invite the reader to skip ahead and to read
the main argument before continuing.

\textbf{Acknowledgments.} I would like to thank John Rognes for many
informative discussions on various subjects involved in the
construction, and Mohammed Abouzaid for many discussions and
suggestions which lead me to rewrite the paper in a more transparent
fashion. Of course any non-transparency and errors are due solely to
me.


\chapter{Finite Dimensional Approximation of Action} \label{finred}

In this section we will define smooth functions on finite dimensional
manifolds approximating the action in Equation~\eqref{action} on
the space of paths restricted to:
\begin{enumerate}
\item{The free loop space $\Lambda T^*N$ of \emph{closed} paths in
    $T^*N$.}
\item{The path space $\Omega_qT^*N$ of paths in $T^*N$ starting and
    ending in the fiber $T^*_q N$ over a point $q\in N$.}
\end{enumerate}
In Contrast to the infinite dimensional cases each approximation in
case 2 will simply be the restriction of an approximation in case 1,
and so the approximations will be compatible with each other. Indeed,
this is because to approximate the action on closed loops we in some
sense insert discontinuities, which are similar to the one of
restricting start and end point to a fiber. These approximations
will of course depend on a Hamiltonian $H$ as the action does, but
they will also depends on a sub-division $\alpha$ of the unit
interval. In fact they will smoothly depend on the pair
$(H,\alpha)$. We will only need to consider a very restricted class of
Hamiltonians. Indeed, we will only consider Hamiltonians linear at
infinity with a slight positive slope.

We fix a Riemannian structure on $N$ once and for all. Throughout the
paper let $\epsin$ denote some fixed constant such that $2\epsin$ is
smaller than the injective radius of $N$ (Recall that we assume $N$ to
be closed).
\begin{Definition}
  \label{def:5}
  The subspace $\HS \subset C^\infty(T^*N,\R)$ of \textbf{admissible}
  Hamiltonians is defined by; $H\in \HS$ iff
  \begin{align} \label{Hdef}
    H(q,p)= \mu \norm{p} + c \qquad \textrm{when} \qquad
    \norm{p}\geq 1
  \end{align}
  with $c\in \R$ and $\mu \in ]0,\epsin]$.
\end{Definition}
The slope $\mu$ is significant since this implies that the Hamiltonian
flow outside $DT^*N$ never starts and ends in the same fiber. Indeed,
the Hamiltonian flow co-insides with the geodesic flow with constant
speed equal to this slope at infinity (see e.g. \cite{hejeh} Section
3).

We let $C_1^H$ denote the supremum of the norm of the covariant
derivative (gradient) of $H$ on the compact set $DT^*N$. Similar let
$C^H_2$ denote the supremum of the norm of the second order covariant
derivative on $DT^*N$. These are both continuous as functions of
$H\in \HS$. Note that $C_1^H$ is a bound on the gradient on all of
$T^*N$ because of the linearity assumption outside of $DT^*N$.

The Riemannian structure on $N$ provides a canonical vector bundle
isomorphism between the tangent bundle $TN$ and the cotangent bundle 
$T^*N$. We will without further warning suppress this isomorphism from
the notation. However, we use the notation $T^*N$ to denote the
canonical symplectic manifold of dimension $2d$, and we use the
notation $TN$ to denote the vector bundle over $N$. The Riemannian
structure on $N$ induces a Riemannian structure on $T^*N$ compatible
with the symplectic structure (see e.g. section 3 in \cite{hejeh}). We
also get induced  Riemannian structures on the manifolds defined
below, by considering them sub-manifolds of either $N^r$ or
$(T^*N)^r=T^*(N^r)$. We define the finite dimensional manifold
approximation to free loops in $N$ by 
\begin{align*}
  \Lambda_r N = \{ \arq = (q_j)_{j\in \Z_r} \in N^r \mid
  \dist(q_j,q_{j+1}) < \epsin \}.
\end{align*}
This is a manifold of dimension $rd$ and we see that its cotangent
bundle is 
\begin{align*}
  T^*\Lambda_r N = \{ (q_j,p_j)_{j\in \Z_r} \in (T^*N)^r \mid 
  \dist(q_j,q_{j+1}) < \epsin \}.
\end{align*}
We will use the notation $\arz=(\arq,\arp)$ for a point in this
space. Compatible with this we use the notation $z_j=(q_j,p_j)$ for a
single ``coordinate'' in $T^*N$. 

We define the evaluations $\ev_i \colon T^*\Lambda_r N \to N$ by
$\ev_i(\arz)=q_i$. We then define
\begin{align*}
  T^\bullet \Omega_{q,r} N = \ev_0^{-1}(q).
\end{align*}
This is the sub-manifold where $q_0$ is fixed. It is the 
cotangent bundle above restricted to the $(r-1)d$-dimensional
sub-manifold
\begin{align*}
  \Omega_{q,r} N = \{ \arq \in \Lambda_r N \mid q_0=q\} \subset
  \Lambda_r N.
\end{align*}
Obviously it is not the cotangent bundle because we have the extra
cotangent vector $p_0$ - hence the $\bullet$ notation.

The approximations will depend on a subdivision of the interval
$I=[0,1]$ into $r$ sub-intervals. These are represented by points in
the simplex $\Delta^{r-1}$ defined by
\begin{align*}
  \Delta^{r-1} = \{\alpha=(\alpha_0,\dots,\alpha_{r-1}) \in \R^{r}\mid
  \alpha_j \geq 0\forall j, \sum_j \alpha_j =1 \}.
\end{align*}
The correspondence is such that the length of the $j^\textrm{th}$
interval is $\alpha_{j-1}$. We define the length $l$ of a subdivision
$\alpha$ as the longest interval in the subdivision given by
$l(\alpha)=\max_j \alpha_j$.

Assume that $\alpha\in \Delta^{r-1}$ is a subdivision and $\gamma
\colon I \to T^*\Lambda_rN$ is a smooth path starting and ending in
the same fiber $T^*_qN$. Also assume that $\norm{\gamma'(t)} l(\alpha)
< \epsin$ for all $t\in I$. We may then define the $\alpha$-dissection
of $\gamma$ to be the point
\begin{align} \label{alpdis}
\arz=(z_j)_{j\in \Z_r} \in T^*\Lambda_r N \qquad \textrm{given by}
\qquad z_j= \gamma(\sum_{i=0}^{j-1} \alpha_i) 
\end{align}
for $j=0,\dots,r-1$. The assumption on the length makes the distance
$\dist(z_j,z_{j-1})$ less than $\delta_0$ so that this is
well-defined. Note that this is well-defined even though the end
points of $\gamma$ are far apart because their projection to $N$
is the same. Indeed, this shows that
\begin{align*}
  \dist(q_{r-1},q_0) =& \dist(q_{r-1},\pi(\gamma(0))) =
  \dist(q_{r-1},\pi(\gamma(1))) \leq \\
  \leq& \dist(z_{r-1},\gamma(1)) = \dist(\gamma(1-\alpha_{r-1}),\gamma(1))< \delta_0.
\end{align*}

we will be $\alpha$-dissecting time-1 flow curves of the Hamiltonian
flow of $H$ so we better assume that
\begin{align} \label{flowbound}
  l(\alpha)C_1^H < \epsin.
\end{align}
The converse of dissecting a time-1 flow curve is to start with an
``approximated curve'' $\arz \in T^*\Lambda_r N$ and create a piece
wise flow curve by defining the $j^{\textrm{th}}$ piece
\begin{align} \label{eq:16}
  \gamma_j\colon[0,\alpha_j]\to T^*N
\end{align}
to be the Hamiltonian flow curve
$\gamma_j(t)=\varphi^{H}_t(z_j)=\varphi^H_t(q_j,p_j)$ for all $j\in
\Z_r$. These may not fit together to a continuous curve, and the
construction of the approximations of the action is all about fixing
this \emph{without} creating additional critical points, and retaining
that the set of connecting orbits is compact.

Let $P_q^{q'} \colon T_{q}N \to T_{q'}N$ denote the parallel transport
induced by the Riemannian structure of tangent vectors along the
unique geodesic when $\dist(q,q')<2\epsin$.

Because of our assumption in Equation~\eqref{flowbound} each of the
flow pieces $\gamma_j$ are shorter than $\epsin$. So we may for
all $j\in \Z_r$ define
\begin{figure}[ht]
  \centering
  \includegraphics{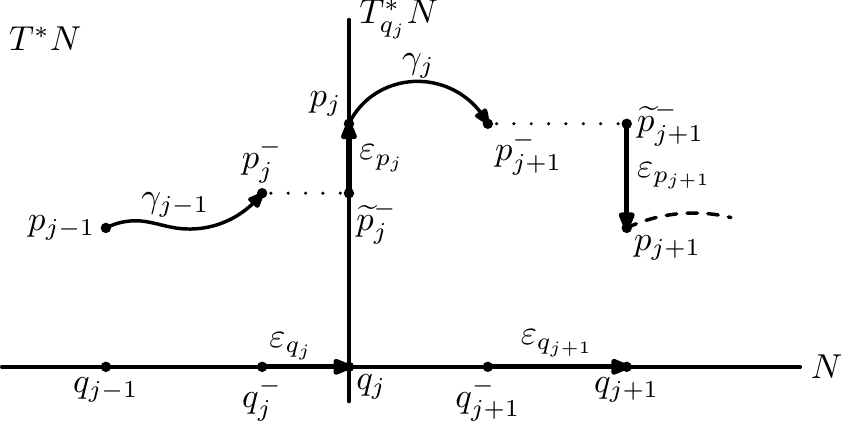}
  \caption{Broken ``geodesics''.}
  \label{fig:pq}
\end{figure}
\begin{align*}
  (q_j^-,p_j^-)&=\gamma_{j-1}(\alpha_j), &
  \tp_j^- &= P_{q_j^-}^{q_j}(p^-) \in T_{q_j}^*N,    \\
  \epsilon_{q_j}&=\exp^{-1}_{q_j^-}(q_j) \in T_{q_j^-} N, &
  \te_{q_j}&=P^{q_{j-1}}_{q_j^-}(\epsilon_{q_j}) \in T_{q_{j-1}}N,  \\
  \epsilon_{p_j}& =p_j-\tp_j^-, \quad \textrm{and}  &
  P &= \max_{j\in \Z_r} \norm{p_j}.
\end{align*}
These, with the exception of the two last in the right column, are
visualized in figure \ref{fig:pq}. Even though all of these depends on
$\alpha, H$ and $\arz$ we suppress this from the notation.

We can now define the finite dimensional approximations. Since we
will only be interested in these when the assumptions of
Lemma~\ref{new42} are satisfied we only define them in such cases. The
assumption in Equation~\eqref{flowbound} may be obtained by simply
assuming that $\delta<\epsin$ in the definition below. However, to get
a buffer around the critical points we assume that $\delta<
\epsin/5$. So all the flow curves $\gamma_j$ has length less than
$\epsin/5$.

\begin{Definition}
  \label{def:Ardef}
  For any subdivision $\alpha\in \Delta^{r-1}$ and any Hamiltonian
  $H\in \HS$ such that
  \begin{align*}
    l(\alpha)(C^H_1+C_2^H) \leq \delta,
  \end{align*}
  with $\delta$ the constant from Lemma~\ref{new42}, we define 
  \begin{align*}
    S_r \colon T^* \Lambda_r N \to \R    
  \end{align*}
  by
  \begin{align} \label{eq:29}
    S_r(\arz) = \sum_{j\in \Z_r} \pare*{\int_{\gamma_j}(\lambda-Hdt) +
    p_j^-\epsilon_{q_j}}.
  \end{align}
\end{Definition}

The integration term in each summand of Equation~\eqref{eq:29} is the
action of the flow-curve piece. The second term is the cotangent
vector $p_j^-$ evaluated on the tangent vector $\epsilon_{q_j}$ and
can be viewed as a way to compensate for the fact that we really
should integrate $\lambda$ over a closed curve. In fact this term is
$\lambda$ integrated over the horizontal geodesic from $(q_j^-,p_j^-)$
to $(q_j,\tp_j^-)$ (which is the doted line in
Figure~\ref{fig:pq}). The integration of $\lambda$ over the
unique line segment connecting $(q_j,\tp_j^-)$ to $(q_j,p_j)$ (the
vertical arrows in Figure~\ref{fig:pq}) is zero. So if we so wished we
could rewrite $S_r$ as the integration of $\lambda$ over a closed
curve depending on $\arz$ plus $\sum_j (\alpha_j H(z_j))$.

By definition we see that: if for each $j\in\Z_r$ the curve $\gamma_j$
ends where $\gamma_{j+1}$ starts then $\arz$ is an $\alpha$-dissection
of a closed periodic orbit for the Hamiltonian flow of $H$ starting at
$z_0\in T^*N$ and the value of $S_r$ is the action of this curve. The
following lemma and corollary is taken from \cite{hejeh}.

\begin{Lemma} \label{new42}
  There are constants $K>0$ and $0<\delta<\epsin/5$ independent of
  $r$, $H\in \HS$ and $\alpha\in \Delta^{r-1}$ such that the
  assumption
  \begin{align} \label{appas}
    l(\alpha)(C_1^H+C_2^H)\leq \delta
  \end{align}
  in Definition~\ref{def:Ardef} implies the estimates
  \begin{align} \label{pjap}
    &\norm{\nabla_{p_j} S_{r} - \te_{q_{j+1}}} \leq 
    \tfrac{1}{4}\norm{\te_{q_{j+1}}}
    =\tfrac{1}{4}\norm{\epsilon_{q_{j+1}}} 
    \qquad \textrm{and} \\
  \label{qjap}
    &\norm{\nabla_{q_j} S_{r} + \epsilon_{p_j}} \leq K \max(1,P)
    (\norm{\epsilon_{q_j}} +
    \norm{\epsilon_{q_{j+1}}})
  \end{align}
  for all $j\in \Z_r$.
\end{Lemma}

Here $\nabla_{q_j} S_r \oplus \nabla_{p_j} S_r =\nabla_{(q_j,p_j)}
S_r = \nabla_{z_j} S_r$ is the gradient with respect to the
$j^\textrm{th}$ factor split into horizontal and vertical vectors.

\begin{Corollary} \label{cor:crit}
  The critical points of $S_r$ are precisely the $\alpha$-dissections
  of the 1-periodic Hamiltonian flow curves of $H$, and the associated
  critical value is the action of this orbit.
\end{Corollary}

Using Lemma~\ref{new42} we can deduce a similar corollary for the
approximations restricted to $T^\bullet \Omega_{q,r} N$, which we
denote by
\begin{align*}
  S_{q,r} = S_{r\mid T^\bullet \Omega_{q,r} N}.
\end{align*}
The proof of this corollary for the restriction is almost identical to
the proof of the above corollary.

\begin{Corollary} \label{cor:critq}
  The critical points of $S_{q,r}$ are precisely the
  $\alpha$-dissections of the time-1 Hamiltonian flow curves of $H$
  which starts and ends in the fiber $T^*_qN$, and the associated
  critical value is the action of this flow curve.
\end{Corollary}

\begin{bevis}
  The gradient estimates in Lemma~\ref{new42} are the same for the
  restriction, except that there is no $\nabla_{q_0} S_{q,r}$
  since $q_0$ is constantly equal to $q$. So if $\arz \in T^\bullet
  \Omega_{q,r} N$ is a critical point we may first use
  Equation~\eqref{pjap} to conclude that $\epsilon_{q_j}=0$ for all
  $j\in \Z_r$. Then we use Equation~\eqref{qjap} to conclude that all
  $\epsilon_{p_j}=0$ for $j\in \Z_r-\{0\}$. This goes both ways so
  this is a necessary and sufficient condition to be a critical
  point. So all the $\gamma_j$ ends where $\gamma_{j+1}$ begins except
  for $\gamma_{r-1}$ which ends in the fiber over $q=q_0$. On the other
  hand any time-1 flow curve starting and ending over $q$ may be
  $\alpha$-dissected and produce a critical point. The critical value
  is the action by the definition of $S_r$ and the fact that all
  $\epsilon_{q_j}$ are 0.
\end{bevis}


\chapter{Pseudo-gradients and
  Compactness/Completeness}\label{cha:pseudo}

In this section we define pseudo-gradients $X_r$ and $Y_{q,r}$ for the
approximations $S_r$ and $S_{q,r}$ from Section~\ref{finred}. These
will smoothly dependent on all parameters ($H$, $\alpha$, $q$). We
will then prove two important compactness/completeness results for
these pseudo-gradients:
\begin{itemize}
\item[C1] {Their flows exists for all times (future \emph{and} past).}
\item[C2] {The rate of change of $S_r$ when flowing with them is bounded
    from below by a positives constant outside a compact set.}
\end{itemize}
These statements are short versions of Lemma~\ref{lem:infflow} and
Lemma~\ref{lem:comp} below, and the combination of these are very
useful. The fact that the pseudo-gradients $Y_{q,r}$ are the
restrictions of a globally defined vector field $Y_r$ is very
important, and this $Y_r$ is an example of what we will define as a
parametrized pseudo-gradient (Definition~\ref{def:6}).

Recall that a pseudo-gradient for a smooth function $f\colon M \to \R$
is a smooth vector field such that $X(f)\geq 0$ and $X=0$ only at
critical points for $f$.

To define the pseudo-gradients we fix once and for all a smooth bump
function $\chi \colon [0,\epsin] \to I$ such that
\begin{align*}
  \chi(t) = \left\{
    \begin{array}{ll}
      1 & t < \epsin/5 \\
      0 & t > \epsin/4.
    \end{array}
  \right..
\end{align*}
Using this we define the smooth bump functions $\chi_r \colon
T^*\Lambda_r N \to I$ by
\begin{align*}
  \chi_r(\arz) = \prod_{j=1}^r \chi(\norm{\epsilon_{q_j}}).
\end{align*}
Note that since $\epsilon_{q_j}$ depends on $H$ and $\alpha$ so will
this function, but smoothly. Assuming $S_r$ is defined we define the
following vector fields.

\begin{Definition} \label{Xdef}
  Let $X_r$ be the smooth vector field on $T^*\Lambda_r N$ given by
  \begin{align*}
    (X_r)_\arz = (\chi_r(\arz)\nabla_{\arq} S_r, \nabla_{\arp} S_r)
  \end{align*}
  using the splitting associated to coordinates $\arz=(\arq,\arp)$.

  Similarly we define the smooth vector field $Y_r$ on $T^*
  \Lambda_r N$ by
  \begin{align*}
    (Y_r)_\arz = (0,\chi_r(\arz)\nabla_{\arq'} S_r, \nabla_{\arp} S_r),
  \end{align*}
  where $\arq'=(q_1,\dots,q_{r-1})$ does not have the $q_0$
  factor.
\end{Definition}

The coordinates in the definition of $Y_r$ are those given by the
splitting of the tangent bundle associated to
$\arz=(q_0,\arq',\arp)$. In both cases we are as in Lemma~\ref{new42}
implicitly using the Riemannian structure to properly define the
horizontal directions associated to $\arq$.

\begin{Remark}
  \label{rem:param}
  The vector field $Y_r$ is vertical with respect to the submersion
  \begin{align*}
    \ev_0 \colon T^*\Lambda_r N \to N.
  \end{align*}
  Indeed, $Y_r$ is simply the orthogonal projection of $X_r$ to the
  vertical vectors (the kernel of $D\ev_0$). These vertical vectors
  are canonically identified with the tangent spaces of $T^\bullet
  \Omega_{q,r} N$, and so we may define 
  \begin{align*}
    Y_{q,r} = Y_{r \mid \ev_0^{-1}(q)}
  \end{align*}
  as a vector field on $T^\bullet \Omega_{q,r} N$.
\end{Remark}

This leads us to the definition of parametrized pseudo-gradients. Let
$M$ be a smooth compact manifold possibly with boundary, corners,
etc. A submersion of smooth manifolds $\pi \colon M' \to M$ is a
smooth map with surjective differential such that $\partial M' =
\pi^{-1}(\partial M)$ - making each fiber a manifold of constant
dimension \emph{without} boundary.

\begin{Definition}
  \label{def:6}
  A \textbf{parametrized pseudo-gradient} for $f\colon M' \to \R$ with
  respect to a submersion $\pi \colon M'\to M$ is a vector field on
  $M'$, which lies in the kernel of $D\pi$, and which in each fiber
  over $q\in M$ is a pseudo-gradient for $f$ restricted to that
  fiber.
\end{Definition}

\begin{Remark}
  \label{rem:comp}
  Since we assume that the base is compact, we only need to check C1
  and C2 fiber-wise. This will come in handy later when we in fact
  will deal with compact smooth families of both $X_r$ and $Y_r$.
\end{Remark}

We will call a critical point of $f_{\mid \pi^{-1}(q)}$ for any $q\in
M$ a \textbf{fiber-critical point}. We will call the value of any
fiber-critical point a \textbf{fiber-critical value}. Note that this
is for all fibers simultaneously, i.e. we do not have an analogue of
Sard's theorem for these. Similarly we define \textbf{fiber-regular
  value}.

\begin{Lemma} \label{pseudo}
  The vector fields $X_r$ and $Y_{q,r}$ are pseudo-gradients for $S_r$
  and $S_{q,r}$ respectively. In particular $Y_r$ is a parametrized
  pseudo-gradient for $S_r$ with respect to $\ev_0 \colon
  T^*\Lambda_rN \to N$.
\end{Lemma}

\begin{bevis}
  The first gradient estimate \eqref{pjap} in Lemma~\ref{new42} tells
  us that $\nabla_{\arp} S_r$ is non-zero on the set where any of the
  $\epsilon_{q_j}$'s are non-zero. So we conclude that on the set
  where we multiply the $\arq$-component of $\nabla S_r$ with $0$ in
  the definition of both $X_r$ and $Y_r$ the $\arp$ components of
  either gradients are non-zero.
\end{bevis}

\begin{Lemma}
  \label{lem:infflow}
  The flows of $\pm X_r$ and $\pm Y_r$ are defined for all time.
\end{Lemma}

\begin{bevis}
  These vector fields are constructed such that the positive or
  negative flows never reaches the ``boundary'' where
  $\dist(q_j,q_{j+1})=\epsin$. Indeed, this distance function is
  preserved by the flows when $\chi$ is zero, and this is the case
  if we assume $\dist(q_j,q_{j+1}) \geq 9\epsin/20$. Indeed, by the
  triangle inequality and this assumption we see that
  \begin{align*}
     \norm{\epsilon_{q_j}} \geq \dist(q_j,q_{j+1}) -\epsin/5 \geq
     9\epsin/20 -\epsin/5 = \epsin /4.
  \end{align*}

  The gradient estimates in Equation~\eqref{pjap} and the fact that
  $\nabla_{p_j} S_r=\nabla_{p_j} S_{q,r}$ are the $p_j$ components of
  both $Y_r$ and $X_r$ proves that
  \begin{align*}
    \absv{\pm Y_r(\norm{p_j})} = \absv{ \pm X_r(\norm{p_j})} \leq
    \norm{\nabla_{p_j} S_r} \norm{(\nabla \norm{p_j})} \leq
    \norm{\nabla_{p_j} S_r} \leq 5\epsin/4
  \end{align*}
  when $p_j \neq 0$. This implies that the rate of change of
  $\norm{p_j}$ is less than $5\epsin/4$ and the lemma follows.
\end{bevis}

\begin{Lemma} \label{lem:comp}
  There is an $\epsilon>0$ such that $X_r(S_r) \geq Y_r(S_r) \geq
  \epsilon$ on the compliment of the compact set
  \begin{align*}
    B_r = \{\arz \in T^*\Lambda_r N \mid \dist(q_j,q_{j+1})\leq \epsin/2
    \quad \vee \quad P \leq 2 \}.
  \end{align*}
\end{Lemma}

\begin{bevis}
  By definition $X_r(S_r) \geq Y_r(S_r)$ so we only need to consider
  $Y_r$. Also by definition of $Y_r$ we see that $Y_r(S_r) \geq
  \norm{Y_r}^2$. So getting a lower bound on $\norm{Y_r}$ suffices.

  We will thus assume that $\arz$ is a point in the set where
  $\norm{Y_r}\leq \epsilon$ intersected with the compliment of $B_r$,
  and arrive at a contradiction for $\epsilon$ small
  enough. However, for notational convenience we will often suppress
  the $\arz$ in the following.

  Using the estimate \eqref{pjap} of $\nabla_{p_j} S_r$ in
  Lemma~\ref{new42} we conclude that
  \begin{align} \label{eboundp}
    \norm{\epsilon_{q_{j+1}}} \leq
    \tfrac{4}{3}\norm{\nabla_{p_j} S_r} \qquad j\in 0,\dots,r-1.
  \end{align}
  Since the $\arp$-part of $Y_r$ is the gradient $\nabla_{\arp}
  S_r$ the above inequalities imply
  \begin{align} \label{qebound22}
    \sum_{j=1}^r \norm{\epsilon_{q_j}}^2 \leq \sum_{j=0}^{r-1} 2
    \norm{\nabla_{p_j}S_r}^2 \leq 2 \norm{Y_r}^2 \leq 2\epsilon^2
  \end{align}
  So with $2\epsilon^2 \leq \epsin/5$ we see that
  \begin{align} \label{eq:18}
    \norm{\epsilon_{q_j}} \leq \epsin/5 \qquad \forall
    j=1,\dots,r,
  \end{align}
  which implies
  \begin{align*}
    \dist(q_j,q_{j+1}) \leq 2\epsin/5 < \epsin/2 \qquad \forall
    j=1,\dots,r.
  \end{align*}
  Since we assumed that we are in the compliment of $B_r$ this implies
  the key observation that
  \begin{align*}
    P\geq 2.
  \end{align*}

  Equation~\eqref{eq:18} also implies that $\chi_r=1$ in
  Definition~\ref{Xdef}, which means that $Y_r$ is the actual
  fiber-wise gradient and we thus get
  \begin{align} \label{gradbound}
    \norm{Y_r}^2 = \sum_{j\in\Z_r} \norm{\nabla_{p_j} S_r}^2 +
    \smashoperator{\sum_{j\in\Z_r-\{0\}}} \norm{\nabla_{q_j} S_r}^2.
  \end{align}
  The inequality $\norm{x-y} \leq \norm{z}$ implies that $\norm{y}^2
  \leq (\norm{x}+\norm{z})^2$. Using (in order) this on the gradient
  estimates in Equation~\eqref{qjap}, the inequality $(c+d+e)^2\leq
  3c^2+3d^2+3e^3$, Equation~\eqref{eboundp} for all $j\in \Z_r$, and
  Equation~\eqref{gradbound} we see that
  \begin{align*}
    \sum_{j=1}^{r-1} \norm{\epsilon_{p_j}}^2 & \leq
    \sum_{j=1}^{r-1} \pare*{
      \norm{\nabla_{q_j} S_r} +
      KP(\norm{\epsilon_{q_j}}+\norm{\epsilon_{q_{j+1}}})}^2 \leq \\
    &  \leq \sum_{j=1}^{r-1} \pare*{
      3\norm{\nabla_{q_j} S_r}^2 + 3K^2P^2(\norm{\epsilon_{q_j}}^2 
      +\norm{\epsilon_{q_{j+1}}}^2)} \leq \\
    &  \leq \sum_{j=1}^{r-1} \pare*{
      3\norm{\nabla_{q_j} S_r}^2 + 6K^2P^2(\norm{\nabla_{p_j} S_r}^2 
      +\norm{\nabla_{p_{j+1}}S_r}^2)} \leq \\
    & \leq \max(3,12K^2P^2) \norm{Y_r}^2 \leq \max(3,12K^2P^2)\epsilon^2
  \end{align*}
  So for $\epsilon^2 < (4\max(3,12K^2)r)^{-1}$ we have
  \begin{align} \label{sepbound}
    \sum_{j=1}^{r-1} \norm{\epsilon_{p_j}}^2 < P^2/4r, \qquad
    \textrm{which implies} \qquad \sum_{j=1}^{r-1} \norm{\epsilon_{p_j}} <
    P/2.
  \end{align}
  This is another key observation and Lemma~\ref{lem:help6} below thus
  implies that $\norm{p_j}\geq P/2$ for all $j\in \Z_r$.

  In particular we see that $\norm{p_j}\geq 1$ for all $j\in
  \Z_r$, and on this set we see that $S_r$ only depends on the
  constants $\mu$ and $c$ (Definition~\ref{def:5}) and the
  sub-division $\alpha$. So in fact the gradient of $S_r$ only depends
  on $\mu$ and $\alpha$. Analyzing the Hamiltonians $H(q,p)=\mu
  \norm{p}$ on $T^*N-DT^*N$ and the associated functions $S_r$ (for
  any $\alpha$) we have
  \begin{align} \label{scale}
    S_r(\arq,t\arp)=tS_r(\arq,\arp),
  \end{align}
  as long as $(\arq,\arp)$ has all $\norm{p_j}\geq 1$ and $t\geq
  1$. This implies that the norm of the gradient of $S_r$ increases
  when $\arp$ is multiplied by $t$. Indeed, we have
  \begin{align} \label{eq:6}
    \norm{\nabla_{\arq} S_r(\arq,t\arp)} = t\norm{\nabla_{\arq}
      S_r(\arq,\arp)}
  \end{align}
  and
  \begin{align} \label{eq:5}
    \norm{\nabla_{t\arp} S_r(\arq,t\arp)} = \norm{\nabla_{\arp}
      S_r(\arq,\arp)}
  \end{align}
  on that same set. 

  The flow defined by these Hamiltonian has no time-1 flow curves
  starting and ending in the same fiber outside of $DT^*N$. Indeed,
  there is no geodesic starting and ending at the same point on $N$
  with length $\mu\in ]0,\epsin]$. This means that on the compact set
  given by the equations 
  \begin{itemize}
  \item{$1 \leq \min_j \norm{p_j} \leq \max_j \norm{p_j}\leq 2$ and}
  \item{$\dist(q_j,q_{j+1})\leq \epsin/2$}
  \end{itemize}
  the norm of $\nabla S_r$ has a lower bound. We thus finally make
  $\epsilon$ smaller than this lower bound.

  Using Equation~\eqref{eq:5} and Equation~\eqref{eq:6} this is a
  lower bound on the gradient of $S_r$ on the non-compact set given by
  the equations 
  \begin{itemize}
  \item{$1\leq P/2 \leq \min_j \norm{p_j} \leq \max_j \norm{p_j}= P$ and}
  \item{$\dist(q_j,q_{j+1})\leq  \epsin/2$,}
  \end{itemize}
  which combined with Lemma~\ref{lem:help6} provides the contradiction.
\end{bevis}

\begin{Lemma}
  \label{lem:help6}
  If $P\geq 2$ and
  \begin{align*}
    \sum_{j\in \Z_r -\{0\}} \norm{\epsilon_{p_j}} \leq P/2
  \end{align*}
  then $\norm{p_j} \geq P/2$ for all $j\in Z_r$.
\end{Lemma}

\begin{bevis}
  This relies on the fact that if either $p_j$ or $\tp_{j+1}^-$ lies
  outside $DT^*N$ then so does the other and
  \begin{align} \label{eq:outside}
    \norm{\tp_{j+1}^-}=\norm{p_j}.
  \end{align}
  Indeed, these two are related by a composition of a parallel
  transport map (which preserves norms) and the Hamiltonian flow,
  which preserves the norm when outside $DT^*N$.

  Since $0\in \Z_r$ is special (omitted from the sum above) it
  is convenient in the following to identify $\Z_r$ with
  $\{0,\dots,r-1\}$ and use the usual ordering on this set. We may
  pick $i$ in this set such that $\norm{p_i}=P\geq 2$. Now
  \emph{assume for contradiction} that we have a $0\leq j<i$ such that
  $\norm{p_j}<P/2$ we may assume WLOG that $j$ is the largest such
  $j<i$ with $\norm{p_j}<P/2$. Then $\norm{p_{j+1}}\geq P/2\geq 1$ and
  using Equation~\eqref{eq:outside} we see that 
  \begin{align*}
    \norm{\epsilon_{p_{j+1}}} = \norm{p_{j+1}-\tp_{j+1}^-} >
    \norm{p_{j+1}}-P/2. 
  \end{align*}
  Using the same equation for all remaining $p_k$'s with $j<k\leq
  i$, which also lies outside $DT^*N$, we get 
  \begin{align*}
    \norm{\epsilon_{p_k}} = \norm{p_k-\tp_k^-} \geq \norm{p_k} -
    \norm{\tp_k^-} = \norm{p_k}-\norm{p_{k-1}}.
  \end{align*}
  Adding these for all $j<k\leq i$ we get that
  \begin{align*}
    \sum_{k=j+1}^i \norm{\epsilon_{p_k}} > \norm{p_i}- P/2 = P-P/2 \geq
    P/2, 
  \end{align*}
  which contradicts Equation~\eqref{sepbound}. The case with $i<j\leq
  r-1$ is similar summing over $k$'s such that $i\leq k <j$.
\end{bevis}

\chapter{Parametrized Conley Indices}\label{cha:paracon}

In this section we define the Conley index (CI) associated to a smooth
function $f$ and a pseudo-gradient $X$ satisfying the two compactness
conditions C1 and C2 from Section~\ref{cha:pseudo}. These are based
compact Hausdorff spaces $I_a^b(f,X)$ depending on a choice of interval
$[a,b]$ and one additional contractible choice. These will come with
a canonical inclusion
\begin{align} \label{eq:212}
  I_a^b(f,X) \subset f^{-1}([a,b]) / f^{-1}(a),
\end{align}
which we will see is a homotopy equivalence. The only reason we do not
use this ``level-set quotient'' as the actual definition is that later
it will be very convenient to use the compact replacements that we
describe in this section. We will also define a category of
parametrized based spaces (called ex-spaces), and define the
parametrized CI (PCI) as an object in such a category. A PCI will
encode the fiber-wise CI, and Lemma~\ref{lem:11} describes how these
are related for different 
bases. In particular Corollary~\ref{cor:4} relates the PCI of
$(S_r,Y_r)$ (whose fibers will be $I_a^b(S_{q,r},Y_{q,r})$) with the
\emph{un}parametrized CI of $(S_r,X_r)$. For simplicity we 
will \emph{always} assume that $b$ is an upper bound on all
fiber-critical values, and hence a fiber-regular value. However, since
we will not be assuming such genericity for $a$ the fiber CI will be a
slight generalization of the usual notion of CIs. See
e.g. \cite{MR797044} for a more conventional treatment of CIs.

Since the unparametrized case is a special case of the parametrized
case (the base is a single point) we formulate everything
parametrized. We will need some basics about parametrized homotopy
theory. As usual a map is a continuous function. In the following $W$
is any compact Hausdorff space.

\begin{Definition} \label{exspace}
  The category $\Com^+_W$ of compact \textbf{ex-spaces} over $W$ has
  \begin{itemize}
  \item {objects (or ex-spaces): $A=(A,p_A,s_A)$, with $A$ a compact
      Hausdorff space, $p_A \colon A \to W$ a map, and $s_A \colon W
      \to A$ a map such that $p_A \circ s_A = \id_W$}
  \item{morphisms (or ex-maps): $f\colon A \to B$ a
      map such that $p_B \circ f = p_A$ and $s_B = f \circ
      s_A$.}
  \end{itemize}
\end{Definition}

We call $p_A$ and $s_A$ the projection of $A$ and the section of $A$
respectively. For each $q\in W$ the \textbf{fiber} $A_{\mid q} =
p_A^{-1}(q)\subset A$ is naturally a based space with base point
$s_A(q)$. So intuitively this category is the category of based
compact Hausdorff spaces continuously parametrized by $W$.

\begin{Definition}
  \label{def:ex-hom}
  For any $A,B \in \Com_W^+$ an \textbf{ex-homotopy} $G \colon A\times I
  \to B$ is a map such that $G_{\mid A\times \{t\}}$ is an ex-map for
  all $t\in I$.  
\end{Definition}

In this definition we used the canonical identification of $A\times
\{t\}$ with $A$. As for spaces this provides the notions of
ex-homotopy equivalences, ex-homotopy type, etc. In particular
restricting any ex-homotopy equivalence $h \colon A \to B$ to the
fibers over $x\in W$ we get a based homotopy equivalence of the fibers
$h_{\mid x} \colon A_{\mid x} \to B_{\mid x}$.

For any pair of spaces $(A,B)$ with a map $\pi\colon A\to W$ we define
the \textbf{fiber-wise quotient} $A:_W^\pi B \in \Com_W^+$ by the push
out diagram
\begin{align} \label{quotpush}
  \xymatrix{
    B \ar[d]^{\pi} \ar[r] & A \ar[d] \\
    W \ar[r] & A:_W^\pi B.
  }
\end{align}
The projection is the push out of the map $\pi_{\mid A}$ and the
identity on $W$, i.e. the unique $\pi'$ in the diagram
\begin{align} \label{quotpush2}
  \xymatrix{
    B \ar[d]^{\pi} \ar[r] & A \ar[d] \ar@/^/[ddr]^{\pi_{\mid A}} \\
    W \ar[r]^-{s'} \ar@/_/[rrd]^{\id} & A:_W^\pi B \ar[rd]^{\pi'} \\
    & & W.
  }
\end{align}
The section is the $s'$ in the diagram. When unambiguous we omit the
projection $\pi$ and write $A:_W B$. The fiber of this construction at
any point $q\in W$ is canonically identified with the based space
given by the quotient $p^{-1}(q)/(p^{-1}(q)\cap B)$. We notice that if
$A$ is compact Hausdorff and $B$ is closed in $A$ then this is an
ex-space in $\Com_W^+$. To make the section more concrete  one can add a
$W$ to both $A$ and $B$. I.e. $A :_W B \cong (A\sqcup W):_W (B\sqcup
W)$. However, this is not actually necessary by the usual definition
of push out.

The fact that we will mostly use these notions for compact Hausdorff
spaces makes many things easier, however for a general treatment we
refer to \cite{MR2271789}. We will at some points use non-compact
ex-spaces. However, the properties we will use for these will be very
basic.

\begin{figure}[ht]
  \centering
  \includegraphics{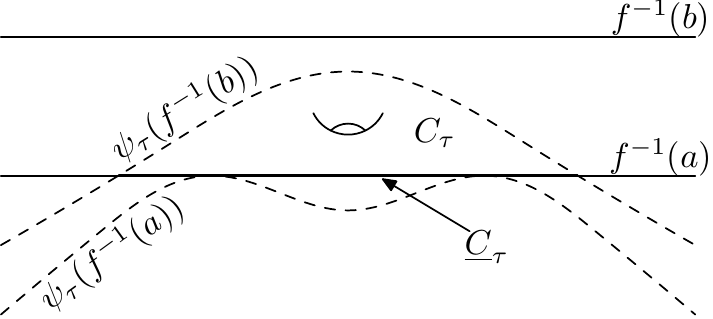}
  \caption{$C_\tau$ and $\uC_\tau$.}
  \label{ctaus}
\end{figure}

Let $f\colon M' \to \R$ be a smooth function and let $Y \colon M' \to
TM'$ be a parametrized pseudo-gradient (Definition~\ref{def:6}) with
respect to a submersion
\begin{align*}
  \pi \colon M' \to M.
\end{align*}
We assume that $Y$ satisfies C1 and C2 from
Section~\ref{cha:pseudo}. So there is an $\epsilon>0$ such that
$Y(f)>\epsilon$ outside some compact set $K \subset M'$. Let $\psi_t$ 
denote the flow of $-Y$, which exists for all time. We may define
\begin{align} \label{eq:15}
  C_\tau = \psi_\tau(f^{-1}([a,b])) \cap f^{-1}([a,b])
\end{align}
and
\begin{align} \label{eq:15b}
  \uC_\tau = \psi_\tau(f^{-1}([a,b])) \cap f^{-1}(a) \subset C_\tau.
\end{align}
These sets are illustrated in Figure~\ref{ctaus}.

\begin{Lemma}
  \label{lem:7}
  The subspace $C_\tau \subset f^{-1}([a,b])$ is compact for
  $\tau\geq \epsilon^{-1}(b-a)$.
\end{Lemma}

\begin{bevis}
  Let $t_0=\epsilon^{-1}(b-a)$. Then define
  \begin{align*}
    K' = \smashoperator{\bigcup_{t \in [0,t_0]}} \psi_t (K)
  \end{align*}
  This is compact, and any point $x\in f^{-1}([a,b])$ will either flow
  into and thus end in this set after time $t_0$ or satisfy 
  \begin{align*}
    \pd{t} (f\circ \gamma) (t) = -Y_{\gamma(t)} (f) < -\epsilon. 
  \end{align*}
  In which case we get $f(\gamma(t_0))-f(\gamma(0)) <
  -\epsilon(b-a)\epsilon^{-1}$. This implies that $f(\gamma(t_0))
  < a$. Ultimately proving that
  \begin{align*}
    C_\tau \subset K'
  \end{align*}
  for all $\tau \geq t_0$. Since $C_\tau$ is closed and $K'$ is
  compact we have proved the lemma.
\end{bevis}

This leads us to define the parametrized Conley index (PCI) as
\begin{align} \label{eq:17}
  I_a^b(f,Y,\tau)_M = I_a^b(f,Y)_M = C_\tau :_M \uC_\tau,
\end{align}
which for $\tau\geq t_0$ is compact and thus lies in $\Com_M^+$. The
choice of $\tau$ is a contractible choice, and we will assume that
$C_\tau$ and thus the PCI is compact if nothing else is mentioned. We
will comment on the choice of $\tau$ when appropriate. The following
lemma justifies this suppression.

\begin{Lemma}
  \label{lem:8}
  The subspace $C_\tau$ is a subspace in $C_{\tau'}$ when $\tau\geq
  \tau'\geq 0$, and the inclusion induces an ex-homotopy equivalence
  \begin{align*}
    I_a^b(f,Y,\tau)_M \to I_a^b(f,Y,\tau')_M
  \end{align*}
  even when the indices are non-compact.
\end{Lemma}

\begin{bevis}
  Since $Y$ is a pseudo-gradient it is clear that
  \begin{align*}
    \psi_{\tau'}(f^{-1}([\infty,b])) \subset
    \psi_\tau(f^{-1}([\infty,b]))
  \end{align*}
  which is preserved when intersecting with $f^{-1}([a,b])$. So
  $C_{\tau} \subset C_{\tau'}$.

  Furthermore, we may define what we will call a \textbf{stop-flow}
  $\underline{\psi}_t$ on $f^{-1}([a,b])$ by the formula
  \begin{align*}
    \underline{\psi}_t(x) = \psi_{s(x,t)}(x),
  \end{align*}
  where $s(x,t) = t$ if $f(\psi_t(x))\geq a$ and
  otherwise we put $s(x,t) = \max \{ s \in [0,t] \mid
  f(\psi_t(x)) \geq a \}$. This is the flow of $\psi_s$ restricted to
  $f^{-1}([a,b])$ and stopped (abruptly) when reaching the set
  $f^{-1}(a)$. This is continuous, and since it preserves the set
  $f^{-1}(a)$, and the submersion to $M$ it defines ex-maps
  \begin{align*}
    \underline{\psi}_{\tau-\tau'} \colon I_a^b(f,Y,\tau)_M \to
    I_a^b(f,Y,\tau')_M,
  \end{align*}
  which is an ex-homotopy inverse to the inclusions. Indeed, in both
  cases the ex-homotopy to the identity is simply given by these maps
  $\underline{\psi}_t$ for $t\in [0,\tau'-\tau]$.
\end{bevis}

Let $\pi' \colon M \to M''$ be another submersion and let $X$ be
another parametrized pseudo-gradient for $f$ this time with respect to
$\pi' \circ \pi \colon M' \to M''$ satisfying C1 and C2. Then we have
two parametrized Conley indices $I_a^b(f,Y)_M$ and $I_a^b(f,X)_{M''}$.

\begin{Lemma}
  \label{lem:11}
  There is a \emph{canonical} (up to contractible choices) ex-homotopy
  equivalence
  \begin{align*}
    I_a^b(f,X)_{M''} \simeq_{M''} I_a^b(f,Y)_M :_{M''} M.
  \end{align*}
\end{Lemma}

The fiber-wise quotient of the section for any map $h\colon W\to W'$
is called the push forward of ex-spaces and is a functor denoted
\begin{align*}
  h_! \colon \Com_W^+ \to \Com_{W'}^+.
\end{align*}
The relation to sheaf-shriek push forward is subtle.

\begin{Corollary}
  \label{cor:4}
  The Conley index $I_a^b(S_r,X_r)$ is homotopy equivalent to
  \begin{align*}
    (I_a^b(f,Y_r)_N)/N,
  \end{align*}
  where $(\sm)/N$ is the functor from ex-spaces over $N$ to based
  spaces given by collapsing the section.
\end{Corollary}

\begin{bevis:lem:11}
  This follows from Lemma~\ref{lem:8} and the fact that by definition
  the level-set quotients satisfies this - i.e.
  \begin{align*}
    \pare*{f^{-1}([a,b]) :_M f^{-1}(a)} :_{M''} M \cong f^{-1}([a,b]) :_{M''}
    f^{-1}(a). 
  \end{align*}
  This ``canonical'' ex-homotopy equivalences are thus defined as the
  maps induced (when taking quotients over $M''$) by the composition
  of
  \begin{itemize}
  \item {the inclusion $C_\tau \to f^{-1}([a,b])$ and }
  \item {the flow of $-X$ taking the image into a set $C'_{\tau'}$
      defining $I_a^b(f,X)_{M''}$.}
  \end{itemize}
  Here $C_\tau$ is defined using $-Y$ and $C'_{\tau'}$ is defined
  using $-X$.
\end{bevis:lem:11}


\chapter{The Suspension Maps} \label{cha:fib-susp}

In this section we will for a given Hamiltonian $H$ relate the
Conley indices for $S_r$ and $S_{q,r}$ defined in
Section~\ref{cha:paracon} for different choices of $r$. These turn out
to be related by twisted suspensions, which essentially is why we can
use these to define spectra. More specifically assume that $H$ is a
Hamiltonian and $\alpha \in \Delta^{r-1}$ is a subdivision such that
$S_r$ is defined. Then the top face map
\begin{align} \label{eq:24}
  d_r \colon \Delta^{r-1} \subset \Delta^r
\end{align}
given by appending a $0$ defines $S_{r+1}$ (and $Y_{r+1}$ and $X_{r+1}$)
using the same Hamiltonian $H$. Then for each $q\in N$ we construct a
map
\begin{align*}
  \sigma_{q,r} \colon \Sigma^d I_a^b(S_r,Y_r)_q \to
  I_a^b(S_{r+1},Y_{r+1})_q.
\end{align*}
Here $\Sigma^d(-)$ denotes the $d$-fold reduced suspension of based
topological spaces. We construct these such that they are defined
continuously over $N$, and in the generic case they will induce
homotopy equivalences of unparametrized indices (when $a$ is regular
for $S_r$). Here continuous in $N$ means that we in fact define a map
of ex-spaces 
\begin{align} \label{firstsigmar}
  \sigma_r \colon I_a^b(S_r,Y_r)_N^{TN} \to I_a^b(S_{r+1},Y_{r+1})_N,
\end{align}
where $I_a^b(S_r,Y_r)_N^{TN}$ denotes the twisted fiber-wise
reduced suspension of the index $I_a^b(S_r,Y_r)_N$ using the one-point
compactification of the fibers of $TN$. We start by defining this
properly. We will refer to this procedure as $TN$-suspending.

The construction of the twisted fiber-wise reduced suspension does not
depend on the fact that $TN$ is the tangent bundle of $N$. So assume
that $\pi\colon \beta \to W$ is any metric vector bundle of dimension
$l$ over a compact Hausdorff space $W$. We have the associated disc
bundle.
\begin{align*}
  D\beta = \{v\in \beta \mid \norm{v} \leq 1\}
\end{align*}
and sphere bundle
\begin{align*}
  S\beta = \{v\in \beta \mid \norm{v} = 1\}.
\end{align*}
We then define
\begin{align*}
  S^\beta = (D\beta):_W (S\beta) \in \Com_W^+,
\end{align*}
where $:_W$ is the fiber-wise quotient defined in
Equation~\eqref{quotpush}. This may be identified with the sphere
bundle of $\beta\times \R \to W$, with its canonical section. Or
it may be identified with the fiber-wise 1-point compactification of
$\beta$. 

For any two ex-space $A$ and $B$ over $W$ we also define the
fiber-wise product over $W$ as the usual pull back
\begin{align*}
  A \times_W B = \{(x,y) \in A\times B \mid p_A(x)=p_B(y)\}.
\end{align*}
This has the obvious section and projection and it is the categorical
product. We also define the fiber-wise wedge product as
\begin{align*}
  A \vee_W B = A \bigcup_{W} B
\end{align*}
using the sections $s_A$ and $s_B$. This is the categorical coproduct.
We may then define the fiber-wise smash product as
\begin{align*}
  A \wedge_W B = (A\times_W B) :_W (A\vee_W B).
\end{align*}
These are symmetric monoidal products. The smash product $\wedge_W$
distributes over the wedge $\vee_W$ and defines a symmetric bimonoidal
structure on the category of ex-space over $W$. Restricting to a fiber 
$q\in W$ defines a symmetric bimonoidal functor to based spaces (with
the usual $\wedge$ and $\vee$ as products - so we are really just
saying that these structures on ex-spaces behaves similarly).

Any ex-space $A$ over $W$ may now be fiber-wise reduced suspended
by the sphere bundle $S^\beta$. That is we define the
\textbf{$\beta$-suspension}
\begin{align} \label{eq:relthom}
  A^\beta = A \wedge_W S^\beta.
\end{align}
This notation suggests that we can consider $S$ an ex-space, and
indeed this is $S^0_W$ appearing below. A little more concretely one
may identify $A^\beta$ with
\begin{align*}
  D(p_A^*\beta) :_W \pare*{ D(\beta)  \cup S(p_A^*\beta)},
\end{align*}
Here $D(\beta)$ is identified with the restriction of $D(p_A^*\beta)$
to the image of the section $s_A \colon W \to A$.

\begin{Remark}
  \label{rem:NN}
  The symmetric bimonoidal structure on ex-spaces over $W$ has
  \begin{itemize}
  \item {zero: $W_W = (W,\id_W,\id_W)$ (to be read $W$ over $W$) and}
  \item {unit: $S^0_W$.}
  \end{itemize}
  Here
  \begin{align*}
    S^n_W = (S^n \times W,\pi_W,*\times \id_W),
  \end{align*}
  where $\pi_W \colon S^n \times W \to W$ is the projection and
  $*\times\id_W\colon W\to S^n\times W$ is the constant map to the
  base point of $S^n$ times the identity.

  The reason we denote the zero by $W_W$ and not simply $W$ is because
  this would be very ambiguous with the usual notation for Thom
  spaces. I.e. for $\beta \to W$ we would define its Thom space to be
  \begin{align*}
    W^\beta = D\beta/S\beta \cong ((S^0_W)^\beta)/W = S^\beta/W,
  \end{align*}
  which is not $W_W^\beta (\cong W_W)$. So the subscript emphasizes that
  we are in the category over $W$ such that we know that $(-)^\beta$
  means the $\beta$-suspension.
\end{Remark}

So $I_a^b(S_r,Y_r)^{TN}_N$ is defined as above using the metric 
bundle $TN\to N$ and a choice of the parametrized Conley index from
Section~\ref{cha:paracon}. So we may identify the restriction of this
to the fiber over $q\in N$ with the usual reduced suspension of the
Conley index of that fiber:
\begin{align*}
  \Sigma^d I_a^b(S_{q,r},Y_{q,r}) = I_a^b(S_{q,r},Y_{q,r}) \wedge S^d.
\end{align*}
However, not canonical since it depends on a choice of isomorphism
$T_qN\cong \R^d$. 

Assume now that $(H,\alpha) \in \HS \times \Delta^{r-1}$ are given and
defines $S_r$. Then we may define the stabilized subdivision $\alpha'
\in \Delta^r$ by $\alpha'_j=\alpha_j$ for $j=0,\dots, r-1$ and
$\alpha_r=0$. Indeed, this is the usual top face map $d_r$ from
Equation~\eqref{eq:24} evaluated on $\alpha$. Using $(H,\alpha')$ we
then define the function $S_{r+1}$ and parametrized pseudo-gradient
$Y_{r+1}$.

We will now define $\sigma_r$ in a few steps. Let $\uC_\tau \subset
C_\tau \subset S_r^{-1}([a,b])$ be as defined in
Equation~\eqref{eq:15} and Equation~\eqref{eq:15b} using the function
$S_r$ and parametrized pseudo-gradient $Y_r$. Let $\tau$ be such that
$C_\tau$ and hence
\begin{align*}
  I_a^b(S_r,Y_r)_N = C_\tau :_N (C_\tau \cap S_r^{-1}(a))
\end{align*}
is compact. Let $E \to C_\tau$ be the vector bundle given by $\ev_0^*
TN$ then we define the map
\begin{align} \label{suspdef}
  f=f_r \colon E \to S_{r+1}^{-1}([a,b])
\end{align}
given by
\begin{align} \label{suspind}
  f(z_0,\dots,z_{r-1},v) = (z_0,\dots,z_{r-1},(q_0,v)).
\end{align}
Here $\arz=(z_0,\dots,z_{r-1}) \in T^*\Lambda_r N$ and $v\in
T_{q_0}N$, and the new ``coordinate'' in $f(\arz)$ is simply given by
$(q_0,v)$. So if we denote $\arz'=f(\arz)$ we may write $z_j'=z_j$ for
$j=0,\dots,r-1$ and $z'_r = (q_0,v)$. This notation is a little
abusive since we originally assumed that the indices $j$ are in 
$\Z_r$ and $\Z_{r+1}$, and now we are identifying these with
$\{0,\dots,r-1\}$ and $\{0,\dots,r\}$ respectively.

Notice that this in fact does land in $S_{r+1}^{-1}([a,b])$ because
the diagram
\begin{align} \label{comrdia}
  \xymatrix{
    E \ar[d] \ar[r]^-{f} & T^*\Lambda_{r+1} N \ar[d]^{S_{r+1}} \\
    C_\tau \ar[r]^-{S_r} & \R
  }
\end{align}
commutes. Indeed, the two new terms in the sum are both zero
independent of $v$. This is because $\epsilon_{q_r}=0$ and the flow
curve $\gamma_r$ is the constant curve parametrized by the point $0\in
\R$ so the integral vanishes. In particular we see that $S_{r+1}$ is
constant on the image of $f$ restricted to any fibers of the vector
bundle.

The image of $f$ is not compact, and it does not induce a map of the
quotients that we need it to. However, composing with the flow
$\psi'_{\tau'}$ of $-Y_{r+1}$ and intersecting with
$s_{r+1}^{-1}([a,b])$ will by Lemma~\ref{lem:7} be
compact if $\tau'$ is large enough. That is, if we define the compact
set
\begin{align*}
  C'_{\tau'} = \psi'_{\tau'} (S_{r+1}^{-1}([a,b])) \cap
  S_{r+1}^{-1}([a,b]) 
\end{align*}
as in Equation~\eqref{eq:15}, but for $S_{r+1}$. Then we get
\begin{align*}
  (\psi'_{\tau'} \circ f) (E) \subset C'_{\tau'} \cup
  S_{r+1}^{-1}(]-\infty,a]). 
\end{align*}
We may then (using compactness) find a $R>>0$ such that $S_RE$ (the
sphere bundle of $E$ with radius $R$) is sent by this map to
$S_{r+1}^{-1}(]-\infty,a])$. This means that it induces a map of
quotients
\begin{align*}
  D_RE / (S_RE \cup (D_RE)_{\mid \uC_\tau}) \to C'_{\tau'} \cup
  S_{r+1}^{-1}([a,b]) / S_{r+1}^{-1}(]-\infty,a]). 
\end{align*}
However, all maps are over $N$, so it actually induces a map of the
quotients over $N$. The source of this ex-map is identified as
\begin{align*}
  D_RE :_N (S_RE \cup (D_RE)_{\mid \uC_\tau}) &\cong
  I_a^b(S_r,Y_r)_N^{TN}
\end{align*}
and the target is
\begin{align*}   
  C'_{\tau'} \cup S_{r+1}^{-1}([a,b]) :_N S_{r+1}^{-1}(]-\infty,a])
  \cong C'_{\tau'} :_N \uC'_{\tau'} = I_a^b(S_{r+1},Y_{r+1})_N. 
\end{align*}
So we have defined the wanted map
\begin{align*}
  \sigma_r \colon I_a^b(S_r,Y_r)_N^{TN} \to I_a^b(S_{r+1},Y_{r+1})_N,
\end{align*}
depending on the contractible choices of $\tau$ and $\tau'$.

Since this construction generalizes the construction in \cite{hejeh}
to a fiber-wise construction we immediately get the following
corollary.

\begin{Corollary} \label{cor:ths}
  If $a$ is a regular value for $S_r$ (and thus also
  $S_{r+1}$) then the map
  \begin{align*}
    \sigma_r/N \colon (I_a^b(S_r,Y_r)^{TN}_N)/N \to
    I_a^b(S_{r+1},X_{r+1})
  \end{align*}
  is a homotopy equivalence.
\end{Corollary}

Note that by assumption $b$ is always assumed regular. One could
try to prove more generally that this corollary is true for any $a$,
and even prove a fiber-wise version. However, we will not need this,
and it is thus convenient to avoid the technicalities of this by
simply making sure later that when this corollary is needed $a$ is,
indeed, regular for $S_r \colon T^* \Lambda_r N \to \R$. 

\begin{bevis}
  Corollary~\ref{cor:4} tells us that using the functor $(-)/N$ we
  recover the unparametrized indices $I_a^b(S_r,X_r)$ and
  $I_a^b(S_{r+1},X_{r+1})$. These where in Lemma~4.6 of \cite{hejeh}
  proved to be related by a $TN$-suspension relative to $\uC_\tau$,
  and inspecting the above defined map we see that it induces this
  homotopy equivalence.
\end{bevis}


\chapter{The Generating Function Spectra}\label{cha:gen-spec} 

In this section we assume that we are given:
\begin{enumerate}
\item{a Hamiltonian $H\colon T^*N \to \R$ in $\HS$,} 
\item{a subdivision $\alpha \in \Delta^{r-1}$ such that $S_r$ from
    Section~\ref{finred} is defined using $H$ and $\alpha$, and}
\item{a compact interval $[a,b]\subset \R$ such that
    \begin{itemize}
    \item{ $a$ is not the action of \emph{closed} 1-periodic orbit
        associated to $H$, and}
    \item{ $b$ is a strict upper bound on the action of the time-1
        flow curves starting and ending in the same fiber of $T^*N \to
        N$.}
    \end{itemize} }
\end{enumerate}
After some additional choices (which we suppress from the notation) we
will associate to this a fibered spectrum
\begin{align*}
  \GF_a^b(H,\alpha)  
\end{align*}
over $N$. Here fibered\footnote{Not to be confused with
  \emph{fibrant}} simply means that to each $q\in N$ we have a 
spectrum, and these glue together continuously in a sense which we
will make precise. This fibered spectrum captures the Floer homology
of the action integral restricted to those curves with action above
$a$, but encode what takes place over each fiber. We warn that 
as defined in this section these are \emph{not} parametrized spectra
in the sense of \cite{MR2271789} - unless $\beta$ in the following is
trivial. Appendix~\ref{cha:funcrep}, however, describes a natural
functor to a more standard notion of parametrized spectra.

We will then define the notion of homology of such fibered spectra,
and prove that for $\GF_a^b(H,\alpha)$ this does not depend on
anything but $H$ and $a$, and this will recover the notion of
generating function homology. In fact homology forgets the fibered
structure, and only depends on the unparametrized index. If one wants
information about fibers one should restrict to a fiber and
\emph{then} take homology.

Motivated by the construction in Section~\ref{cha:fib-susp} we
define the following very naive version of a category of spectra over
any compact Hausdorff space $W$. Let $\beta \to W$ be any metric
vector bundle. We will use that $(-)^\beta \colon \Com_W^+ \to \Com_W^+$
(defined in Section~\ref{cha:fib-susp}) is a functor. Indeed,
we may extend the construction to ex-maps $f\colon A\to B$ by smashing
with the identity on $S^B$. I.e. we get 
\begin{align*}
  f^{\beta} = f \wedge_W \id_{S^\beta} \colon A^\beta \to B^\beta.
\end{align*}
We note that there is an obvious natural isomorphism of functors
between $((-)^{\beta_1})^{\beta_2}$ and
$(-)^{\beta_1\oplus\beta_2}$. We will use this implicitly on
iterations of $(-)^\beta$ and use the notation $(-)^{k\beta}$. For
the sake of notation we will (when there is no ambiguity) suppress
these exponents on ex-maps. So we will often simply write $f$ instead
of $f^{k\beta}$. Sometimes we will write $f^\oplus$ to emphasize that
some (possible) exponent has been suppressed.

\begin{Definition}
  \label{def:betaspec}
    A $\beta$-spectrum $\As$ is a sequence of pairs
    $\As=(A_r,\sigma_r)_{r\in \N_0}$ with each $A_r$ an ex-space over
    $W$ and each $\sigma_r$ an ex-map
    \begin{align*}
      \sigma_r=\sigma_r^{\As} \colon (A_r)^\beta \to A_{r+1}
    \end{align*}
    over $W$.
\end{Definition}
We call the $\sigma_r$'s \textbf{structure maps}. We define
$\sigma_{r,r'}=\sigma_{r,r'}^{\As}$ as the composition
\begin{align*}
  \sigma_{r'-1}\circ \sigma_{r'-2} \circ \cdots \circ
  \sigma_r \colon A_r^{(r'-r)\beta} \to A_{r'}.
\end{align*}
So that $\sigma_{r,r+1}=\sigma_r$ and $\sigma_{r,r}=\id$. Here we
have suppressed some $\beta$-exponents. E.g. we have
\begin{align*}
  \sigma_{r,r+2} = \sigma_{r+1} \circ (\sigma_r)^{\beta}
\end{align*}
as the map
\begin{align*}
  A_r^{2\beta} \xrightarrow{\phantom{a}(\sigma_r)^\beta} A_{r+1}^\beta
  \xrightarrow{\sigma_{r+1}} A_{r+2}.
\end{align*}
We will later define morphisms of $\beta$-spectra and we denote the
category of these by $\Sp_W^\beta$.

If we take a $\beta$-spectrum and restrict it level-wise to the
fiber over a single point $q\in W$ then we recover a naive
spectrum (this identification depends on a trivialization
of $\beta_q$). That is, we get a sequence of based spaces $A_r$ and a
sequence of maps $\sigma_r \colon \Sigma^l(A_r) \to
A_{r+1}$, where $\Sigma^l(-)$ denotes reduced suspension. The
``restriction'' of the suppressed natural isomorphism of functors from
$((-)^\beta)^\beta$ to $(-)^{2\beta}$ is then simply the natural
isomorphisms of functors from $((-)\wedge S^l) \wedge S^l$ to
$(-) \wedge S^{2l}$.

Now let $H$ be any Hamiltonian as in Section~\ref{finred} and $\alpha
\in \Delta^{r_0-1}$ with $\alpha_{r_0} \neq 0$. We will think of this
as a choice of subdivision for all $r\geq r_0$ by using the top
face inclusion $\Delta^{r-1} \subset \Delta^r$ defined by appending
any sub-division by a zero. This was already introduced when we
defined the ex-map $\sigma_r$ in Section~\ref{cha:fib-susp}. Using
these maps we define the infinite simplex $\Delta^\infty$ by taking
the topological union 
\begin{align} \label{infsimp}
  \Delta^\infty = \cup_{r\in \N} \Delta^{r-1}
\end{align}
with the weak topology. for any point $\alpha\in \Delta^\infty$ we
define $r_\alpha \in \N$ by
\begin{align}\label{rofinf}
  r_\alpha = \inf \{ r \in \N \mid \alpha \in \Delta^{r-1}\}.
\end{align}
Now let $\alpha \in \Delta^\infty$ be such that $S_{r_\alpha}$ is
defined using $H$ and $\alpha$. This implies that $S_r$ is defined for
all $r\geq r_\alpha$, which means that the Conley indices
$I_a^b(s_r,Y_r)_N$ are defined. Also all the suspension ex-maps
$\sigma_r$ in Section~\ref{cha:fib-susp} between these indices for
different $r$'s are defined.  

We thus define the $TN$-spectrum
\begin{align} \label{gfspec}
  \GF_a^b (H,\alpha)
\end{align}
by $\GF_a^b(H,\alpha)_r=N_N$ (from Remark~\ref{rem:NN}) for
$r<r_\alpha$. For each $r\geq r_\alpha$ we define
\begin{align*}
  \GF_a^b(H,\alpha)_r=I_a^b(S_r,Y_r)_N.
\end{align*}
We define $\sigma_r$ by the canonical identification $N_N^{TN} \cong
N_N$ when $r<r_\alpha-1$, and we define $\sigma_{r_\alpha-1}$ by the
section (or the fact that $N_N$ is initial in $\Com_N^+$). For
$r>r_\alpha$ we use the ex-maps
\begin{align*}
  \sigma_{r-1} \colon \GF_a^b(H,\alpha)^{TN}_{r-1} \to
  \GF_a^b(H,\alpha)_r
\end{align*}
constructed in Section~\ref{cha:fib-susp} as the structure maps.

\begin{Remark}
  \label{rem:2}
  It may seem odd that this fibered
  spectrum is trivial on the first levels. However, this is a general
  fact about spectra. It does not really matter up to equivalence
  if one replaces the first finite number of spaces by the trivial
  space. In fact it is \emph{intuitively} convenient to think of
  spectra as a sort of colimit as $r\to \infty$ (which is not defined
  in the category of spaces). This is reflected in the following
  definition of homology.
\end{Remark}

Let $\F$ be either $\Q$ or $\F_p$ for some prime $p$. Assume that
$\beta$ is $\F$ orientable and that we have chosen an
orientation. Using the functor in Appendix~\ref{cha:funcrep} we could
define this for any general homology theory, but to make things as
transparent as possible we wont. Recall that $N$ is viewed as a
subspace of $A_r$ using its section. So in the following
$H_*(A_r,N;\F)$ denotes the usual relative singular homology groups.

\begin{Definition}
  \label{def:3}
  The homology of a $\beta$-spectrum $\As$ with coefficients in
  $\F$ is defined to be the colimit
  \begin{align*}
    H_*(\As;\F) = \colim_{r\to \infty} H_{*+lr} (A_r,N;\F)
  \end{align*}
  using the maps $H_*(A_r,N;\F) \to H_{*+l}(A_{r+1},N;\F)$ given by
  the composition of the maps:
  \begin{itemize}
  \item{the Thom-map $H_*(A_r,N;\F) \to H_{*+l}(A_r^\beta,N;\F)$, and}
  \item{The map induced by the structure map $\sigma_r \colon
      A_r^\beta \to A_{r+1}$}
  \end{itemize}
\end{Definition}

The \textbf{Thom-map} will when the sections are cofibrant be a Thom
isomorphism. This is, indeed, the case for $\GF_a^b(H,\alpha)$ since
$a$ is regular (see the proof below). However, we want to talk about
the homology of fibers as well. So we do not restrict the definition
to this case, and in general we define the Thom-map by the following
steps.
\begin{itemize}
\item {First we use the Thom-isomorphism 
    \begin{align*}
      H_*(A_r,N;\F) \cong H_{*+l}(Dp^*\beta,Sp^*\beta\cup D\beta;\F),
    \end{align*}
    of pairs, which depends on the orientation on $\beta$. In the
    proof of Proposition~\ref{Spectral} a different perspective on
    this map is used in the case where $\beta$ is trivial.}
\item {Then we compose with the map induced by the fiber-wise collapse
    map $(Dp^*\beta,Sp^*\beta\cup D\beta) \to (A_r^\beta,N)$
    defining $A_r^\beta$ in Equation~\eqref{eq:relthom}.}
\end{itemize}

Although we have not defined morphisms of $TN$-spectra yet we note
that homology is a functor.

\begin{Lemma} \label{lem:12}
  The homology $H_*(\GF_a^b(H,\alpha);\F)$ is isomorphic to
  \begin{align*}
    H_{*+lr_\alpha}(I_a^b(S_{r_\alpha},Y_{r_\alpha})_N,N;\F) \cong
    \widetilde{H}_{*+lr_\alpha}(I_a^b(S_{r_\alpha},X_{r_\alpha});\F).
  \end{align*}
\end{Lemma}

Note that the notion of generating function homology (associated to
the action on \emph{closed} loops) discussed in e.g. \cite{MR1403954}
is recovered in this lemma. However, we will later only use the
explicit computations in \cite{hejeh}, and not relate it to actual
(twisted) Floer homology.

There is a notion of an ex-space being well-sectioned, which is more
than the section being a standard cofibration of spaces. This notion
does \emph{not} appear in the following proof.

\begin{bevis}
  If we can argue that all the sections are cofibrations then
  Corollary~\ref{cor:4} tells us that
  $H_{*+lr_\alpha}(I_a^b(S_r,Y_r)_N,N;\F) \cong
  \widetilde{H}_{*+lr_\alpha}(I_a^b(S_r,X_r);\F)$. Combining this with
  Corollary~\ref{cor:ths} we see that in this case the colimit in
  Definition~\ref{def:3} defining the homology is the colimit of the
  diagram 
  \begin{align*}
    0 \cdots 0 \to
    H_*(I_a^b(S_{r_\alpha},Y_{r_\alpha})_N,N;\F) \xrightarrow{\cong} 
    H_*(I_a^b(S_{r_\alpha+1},Y_{r_\alpha+1})_N,N;\F) \xrightarrow{\cong} \cdots
  \end{align*}
  proving the corollary.

  The sections \emph{are} cofibrations. Indeed, since $a$ is a regular
  value for $S_r$ the inclusion $\uC_\tau \subset C_\tau$ (defined in
  Equation~\eqref{eq:15} and Equation~\eqref{eq:15b}) is a
  cofibration, and for any cofibrant pair $(A,B)$ with a map to
  $N$ the section of $A:_N B$ is cofibrant.
\end{bevis}

\begin{Corollary}
  \label{cor:1}
  The homology of $\GF_a^b(H,\alpha)$ does not depend on $\alpha$.
\end{Corollary}

\begin{bevis}
  In the proof above we related the homologies of the parametrized
  indices to the homologies of the unparametrized indices, and these
  satisfies homotopy invariance. So the corollary follows since the
  fact that $a$ and $b$ are regular values for $S_r$ does not depend
  on the sub-division $\alpha$.
\end{bevis}


\chapter{The $TN$-Spectrum $\FL$ Associated to $L$ Exact
  in $T^*N$} \label{theex} 

Let $j \colon L \to T^*N$ be as in the introduction. In this section we
construct the associated $TN$-spectrum $\FL$ as a ``colimit'' of
certain generating function spectra. We put colimit in quotes since we
have not yet discussed morphisms of $TN$-spectra. We will instead
construct $\FL$ explicitly. The homology of this $TN$-spectrum will be
the homology of the free loop space of $L$ with possibly twisted
coefficients. It will in fact represent the symplectic homology of
$DT^*L$, but with possibly different twisted coefficients - we
discussed these local coefficients in Remark~\ref{rem:coef} and
Appendix~\ref{cha:coor}. The actual construction will as the
generating function spectra depend on a lot of choices. Indeed, we
will construct a smooth family
\begin{align*}
  (H^s,\alpha^s) \qquad \textrm{for} \qquad s\in [0,\infty[,
\end{align*}
where this family consists of a Hamiltonian and a subdivision. Then
for each $s$ we may pick $r$ large enough to define associated finite
dimensional approximations $S_r^s$, and define associated
(un)parametrized pseudo-gradients $X_r^S$ and $Y_r^s$ - all as in
Section~\ref{finred} and Section~\ref{cha:pseudo}. This will
define a family of generating function spectra, and we will explicitly
describe the ``colimit'' of these as $s$ tends to infinity. Since we
will be tweaking these choices in later sections, and work with more
than on set of choices simultaneously, we introduce the terminology
that a set of choices defines an \textbf{instance} of $\FL$. However,
we note that $\FL$ is in the homotopy category of $TN$-spectra defined
up to unique isomorphism, and we will discuss this in
Appendix~\ref{cha:hominv}.

By choosing an appropriate extension (and rescaling) of $j$ and an
appropriate Riemannian structure on $L$ we may assume that we have a
symplectic embedding of the closed unit disc bundle $D T^*L$ into
$D_{1/2} T^*N$. By abuse of notation we will also denote this
embedding $j$.

Let $\lambda_N$ and $\lambda_L$ denote the canonical Liouville forms
on $T^*N$ and $T^*L$ respectively. Similarly we will let $\norm{\sm}_N
\colon T^*N \to \R$ and $\norm{\sm}_L\colon T^*L \to \R$ denote the
norm of the cotangent vectors in $N$ and $L$ respectively.

Because $j$ is exact symplectic we may choose a function $F\colon
DT^*L\to \R$ such that $dF= j^*\lambda_N - \lambda_L$, and this is
unique up to a constant because $L$ is connected. Since $DT^*L$ is
compact $F$ and its derivatives are bounded. Because of this we will
assume that $F$ takes the value 0, but is non-negative, and hence
$\norm{F}=\norm{F}_\infty$ is a bound on $\absv{F(z_1)-F(z_2)}$.

If $\gamma \colon I \to DT^*L$ is any smooth path we can now consider
its action in $T^*L$ (using the 1-form $\lambda_L$) or the action of
$j\circ \gamma$ in $T^*N$ (using $\lambda_N$). We denote these two
different action functionals by $A_L$ and $A_N$ respectively. They
satisfy the formula
\begin{align} \label{eq:1}
  A_N(j\circ \gamma) = A_L(\gamma) + (F(\gamma(1))-F(\gamma(0))).
\end{align}
The extra terms vanish on closed loops. So on closed loops the two
action integrals agree.

To create $\FL$ we have to construct some families of Hamiltonians
associated to $j$. We will describe some smooth functions and some of
their properties in the following, and we postpone the actual
construction of these to Appendix~\ref{cha:pd-ham}. Indeed, we do so
because we need them to satisfy a lot of additional properties later,
which we do not want to mention here.

\begin{figure}[ht]
  \centering
  \includegraphics{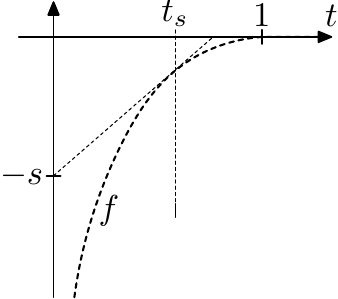} $\qquad$
  \includegraphics{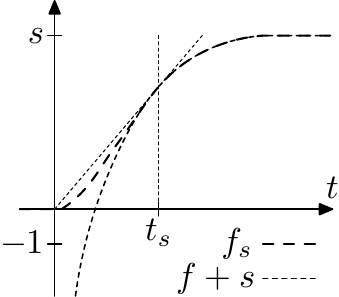}
  \caption{The function $f$ on the left, and $f_s$ and $f+s$ on the right.}
  \label{fig:fs}
\end{figure}

First we fix a smooth concave function $f\colon \R_+ \to \R$ (see
figure~\ref{fig:fs}) such that:
\begin{itemize}
\item { $f(t) \to -\infty$ when $t\to 0$.}
\item { $f(t) = 0$ when $t\in [1,\infty[$.}
\item { $f''(t) < 0$ when $t\in ]0,1[$ - i.e. strictly concave on
    $]0,1[$.}
\end{itemize}
For each $s>0$ we define $t_s$ to be the unique point $t_s\in ]0,1[$
such that the tangent of $f$ at $t_s$ intersects the $2^\textrm{nd}$
axis at $-s$ (see figure~\ref{fig:fs}).
By assumptions on $f$ this $t_s$ will go monotonely to 0 as $s$ goes
to infinity.

For any smooth function $g \colon \R_{\geq 0} \to \R$ we define the
\textbf{Action set} of $g$ as the set of $-s$ such that there is a
tangent to $g$ with slope the length of a \emph{closed} geodesic on
$L$ intersecting $(0,-s)$. We also define this for $f$ even though
$f$ is not defined at 0. The action set of $f$ consists of non-negative
real numbers. If we define the Hamiltonian $H(q,p)=f(\norm{p}_L)$ away
from the zero-section on $T^*L$ the action of 1-periodic orbits would
coincide with this action set of $f$ - hence the name. Indeed, this is
a standard calculation (see e.g. \cite{hejeh} Section 3). Of course we
cannot use this Hamiltonian as it is not extendable to $L\subset
T^*L$. It will thus be convenient to translate the function $f$ and
cap it of depending on a parameter $s$.

So let $f_s\colon \R \to \R$ for $s\in [0,\infty[$ be the smooth 
family of smooth increasing functions constructed in
Appendix~\ref{cha:pd-ham} using $f$. These satisfy
\begin{itemize}
\item[f2)]{For $t$ close to $0$ we have $f_s(t) = ct^2$ for some
    $c=c(s)>0$.}
\item[f3)]{For $t\geq 1$ we have $f_s(t) = s$.}
\item[f4)]{For $s>0$ the restriction of $f_s(t)$ to $]0,1[$ is
    strictly increasing.}
\end{itemize}
These forces $f_0=0$. In this section we particularly need:
\begin{itemize}
\item[f5)]{The tangent to $f_s$ at any $t\in [0,t_s]$
    intersects the $2^\textrm{nd}$ axis in $]-1,0]$.}
\item[f6)]{For $s\geq 5$ we have $f_s(t)=f(t)+s$ when $t\geq t_s$ -
    making the tangents of $f_s$ at these points intersect the
    $2^\textrm{nd}$ axis in $[0,s]$.}
\end{itemize}
In the latter 5 is a rather arbitrary choice since we in this section
really only care about the limit $s\to \infty$. The function $f_s$
(for $s\geq 5$) is illustrated in Figure~\ref{fig:fs}. Let $H_\infty
\colon T^*N \to \R$ be a function in $\HS_\infty\subset \HS$ from
Appendix~\ref{cha:pd-ham}. In particular $H_\infty$ satisfies
\begin{itemize}
\item {$H_\infty$ is zero on $D_{2/3}T^*N$.}
\item {Any time-1 Hamiltonian flow curve for $H_\infty$ starting and
    ending in the same fiber $T^*_qN$ is constant.}
\end{itemize}
Now define the Hamiltonians $H^s \colon T^*N \to \R$ for \emph{all}
$s\geq 0$ by
\begin{align*}
  H^s(z) = \left\{ 
    \begin{array}{ll}
      H_\infty(z) + s & z\notin j(DT^*L) \\
      f_s(\norm{j^{-1}(z)}_L) & z\in j(DT^*L).
    \end{array}
  \right.
\end{align*}
\begin{figure}[ht]
  \centering
  \includegraphics{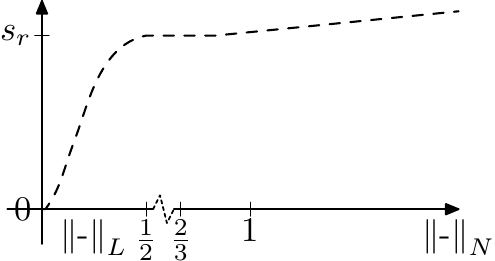}
  \caption{Hamiltonian $H^{s_r}$}
  \label{fig:Hamsr}
\end{figure}

The reason for introducing this $H_\infty$ is simply to get a slight
positive slope at infinity - so that $H^s$ is in $\HS$. These are
smooth because each function $f_s(\norm{j^{-1}(z)}_L)$ is smoothly
extendable by the constant function $s$ on a neighborhood of
$j(DT^*L)$ and this agrees with $H_\infty(z)+s$. Since $f_s$ is smooth
in $s$ we in fact have that the adjoint function
\begin{align*}
  H^\sm \colon T^*N \times [0,\infty[ \to \R
\end{align*}
is smooth - making $H^s$ a smooth family of smooth Hamiltonians.

We wish to explicitly construct the $TN$-spectrum $\FL$ by a limiting
process as $s\to \infty$ of certain generating function spectra
associated to $H^s$ defined in Section~\ref{cha:gen-spec}. First we
consider the critical set associated to the action on \emph{closed}
loops (corresponding to the unparametrized Conley indices).

\begin{Lemma}
  \label{lem:closact}
  For any $s\in [0,\infty]$ the action of the closed 1-periodic orbits
  associated to the Hamiltonian $H^s$ with negative action coincides
  with the action set of $f$ translated down by $s$ and intersected
  with $]-\infty,0[$.
\end{Lemma}

\begin{bevis}
  By Equation~\eqref{eq:1} the action $A_L$ and $A_N$ on closed loops
  agree in $DT^*L$ and $T^*N$. By the explicit construction of $H^s$
  (and $H_\infty$ in Appendix~\ref{cha:pd-ham}) any closed 1-periodic
  orbit in $T^*N$ is either:
  \begin{itemize}
  \item {constant in $D_{2/3}T^*N-j(DT^*L)$, or}
  \item {a 1-periodic orbit in $j(DT^*L)$.}
  \end{itemize}
  The constant curves has critical value $-s$, and this corresponds to
  the unique horizontal tangent to $f_s$ intersecting the
  $2^{\textrm{nd}}$ axis at $s$. The orbits inside $j(DT^*L)$ can be
  calculate using the same precise formula that defines the action set
  of $f_s$. we conclude that there is a 1-1 correspondence between the
  action of 1-periodic orbits and the action set of $f_s$.

  For the purpose of this section we now assume that $s\in
  [5,\infty]$. We postpone the rest of the proof to
  Appendix~\ref{cha:pd-ham}. Indeed, we will not need this for small
  $s$ in this section.

  The action set of $f+s$ is the action set of $f$ shifted down by
  $s$. So we need to see that the negative part of the action
  set of $f_s$ is the same as the negative part of the action set of
  $f+s$. This is indeed the case by construction of the family
  $f_s$ for $s\geq 5$. Indeed, f5 and f6 shows that we have been 
  careful to cap of $f+s$ to get $f_s$ without distorting the
  negative part of the action set.
\end{bevis}

The above lemma makes it rather easy to construct continuation maps if
we did not care about the fiber-wise structure of the Conley
indices. However, we need the fiber-wise structure, and thus we need
similar control of fiber-critical values. So now we consider these.

\begin{Lemma} \label{actions}
  For any $s\geq 5$ all time-1 Hamiltonian flow curves of $H^s$
  starting and ending in the same fiber has action less than
  $\norm{F}+1$, and if the action is less than or equal to $-\norm{F}$
  the flow curve will lie in the set $T^*N - j(D_{t_s}T^*L)$.
\end{Lemma}

The most important part of this construction is that we can control
the \emph{open} flow lines with action less than the fixed number
$-\norm{F}$ in such a way that their action is simply translated
downwards when $s$ is increased. This we fully prove in
Lemma~\ref{lem:10}, but first we need this lemma.

\begin{bevis}
  On the set $j(DT^*L)$ the Hamiltonian flow is constructed such
  that it preserves $\norm{j^{-1}(\cdot)}_L$. So either all of a
  flow-curve is inside $DT^*L$ (or $D_{t_s}T^*L$) or none of it
  is. Given any time-1 Hamiltonian flow curve $\gamma \in
  C^\infty(I,DT^*L)$ we consider the two action integrals in
  Equation~\eqref{eq:1}. The action $A_L$ of $\gamma$ in $DT^*L$ is
  given by the usual intersection formula, which works on open paths
  in $DT^*L$, because here $H^s$ is a function of $\norm{-}_L$. So the
  action $A_N$ of $j\circ \gamma$ in $T^*N$ is given by this
  \emph{plus} the extra factors $F(\gamma(1))-F(\gamma(0))$.

  The first statement then follows by using the above on any time-1
  flow curve $\gamma$ in $j(DT^*L)$ to see that
  \begin{align*}
    A_N(\gamma) \leq 1+F(\gamma(1))-F(\gamma(0)) \leq \norm{F} +1.
  \end{align*}
  Here $1$ comes as minus the smallest possible intersection
  from f5 (combined with f6 this is a global bound on the
  intersection with the $2^{\textrm{dn}}$ axis). By construction of
  $H_\infty$ any flow-curve outside of $j(DT^*L)$ which starts and
  ends in the same fiber is constant and have action $-s\leq -5$.
  
  The second statement follows from using it on any flow curve
  $\gamma$ inside $j(D_{t_s}T^*L)$ to get:
  \begin{align*}
    A_N(\gamma) \geq 0 +  F(\gamma(1))-F(\gamma(0)) \geq -\norm{F}
  \end{align*}
  as needed. Here we used the upper bound 0 from f5.
\end{bevis}

The interval of action (energy) $[a,b]$ that we will be interested in
as $s$ tends to infinity are determined by the above lemma. That is we
will now assume that $a$ and $b$ are fixed such that
\begin{align*}
  a < -\norm{F} \qquad \textrm{and} \qquad b \geq \norm{F}+1.
\end{align*}
Indeed, this makes $b$ an upper bound on all fiber-critical values
of the action integrals as we have been assuming in many constructions
thus far. The reason for the bound on $a$ is that the critical values
below or equal to $a$ behaves in a very concrete manor (even
fiber-wise), and Lemma~\ref{lem:10} below is the precise description
of this behavior. This idea and Lemma~\ref{lem:10} will be essential
in constructing the fiber-wise continuation maps and proving the
fibrancy property.

We need a suitable family of sub-divisions $\alpha^s, s\in
[0,\infty[$. This means that we need to choose a smooth map
\begin{align*}
  \alpha^\sm \colon [0,\infty[ \to \Delta^\infty
\end{align*}
such that the finite dimensional approximation $S_{r_{\alpha^s}}^s$
(see Equation~\eqref{rofinf}) are defined using $(H^s,\alpha^s)$ for
all $s\in [0,\infty[$. Here smoothness is well-defined because the top
face inclusions $\Delta^{r-1} \subset \Delta^r$ are smooth and we used
the weak topology on $\Delta^\infty$. Such families are easy to
construct because all we need $\alpha^s$ to satisfy is
\begin{align*}
  l(\alpha^s) (C_1^{H^s} + C_2^{H^s}) \leq \delta,
\end{align*}
with $\delta$ as in Definition~\ref{def:Ardef}. We may, in fact,
choose this such that $\alpha^s$ is constantly equal to the unique
sub-division in $\Delta^0$ for $s$ close to $0$. Indeed,
this follows from the fact that $f_0=0$ and property H4 of
$H_\infty$ described in Appendix~\ref{cha:pd-ham}.

\begin{Lemma}
  \label{lem:sr}
  There exists a strictly increasing sequence $(s_r)_{r\in \N} \in
  [0,\infty[$ converging to infinity such that 
  \begin{itemize}
  \item { $S_r^s$ is defined for all $s\in [0,s_r]$, and}
  \item { The action associated to $H^{s_r}$ has no 1-periodic orbit
      with action $a$ - i.e. $a$ is regular for $S_r^{s_r}$.}
  \end{itemize}
\end{Lemma}

\begin{bevis}
  By construction $r_{\alpha^s}$ is defined and finite for all $s\in
  [0,s']$ and by compactness of $[0,s']$ it is bounded on such
  intervals (here we use the topology on $\Delta^\infty$). So we may
  indeed pick a sequence $s_r$ (for all $r\geq 1$) such that $S_r^s$ is
  defined for $s\in [0,s_r]$. Making each $s_r$ smaller preserves this
  property and since $r_{\alpha^s}$ starts out being constantly equal
  to $1$ at $s=0$ we can always rechoose the sequence to make it strictly
  increasing. By the shift by $s$ in Lemma~\ref{lem:closact} and
  Sard's theorem we see that since $a$ is fixed and negative the subset
  of $\Jinf$ where the last bullet point is satisfied is open and
  dense. Hence, we may perturb any sequence slightly downwards to
  satisfy this. 
\end{bevis}

We now pick such a sequence $s_r \in [0,\infty[$ and call them
\textbf{jump points}. As in Section~\ref{finred} we may thus define the
fiber-wise approximations 
\begin{align*}
  S^s_{q,r} \colon T^\bullet \Omega_{q,r} N \to \R
\end{align*}
as the restriction of $S_r^s$ to $T^\bullet \Omega_{q,r} N$ for each
$s\in \Jr$. We use the notation
\begin{align*}
  S_r^\sm \colon T^*\Lambda_r N \times \Jr \to \R
\end{align*}
for the function defined by $S_r^\sm(\arz,s)=S_r^s(\arz)$. Similarly
for $S_{q,r}^\sm$. These functions are smooth because $H^s$ and
$\alpha^s$ are smooth in $s$.

We may also define the pseudo-gradients $X_r^s$ and parametrized
pseudo-gradients $Y_{q,r}^s$ as in Section~\ref{cha:pseudo}. In fact
we define $X_r^\sm$ to be the unique vector field on $T^*\Lambda_r N
\times \Jr$ which restricts (and lies in the correct sub-space) to
$X_r^s$ on $T^*\Lambda_r N\times \{s\}$ for each $s\in \Jr$. This
makes $X_r^\sm$ a parametrized pseudo-gradient with respect to the
projection to the second factor 
\begin{align*}
  T^*\Lambda_r N \times \Jr \to \Jr.
\end{align*}
Note that $X^\sm_r$ is smooth because of how it is defined in
Section~\ref{finred} and because $S_r^\sm$ is smooth.

Similar we may define $Y_r^\sm$ making it a parametrized
pseudo-gradient with respect to the projection
\begin{align} \label{projection}
  \pi = \ev_0 \times \id \colon T^*\Lambda N \times \Jr \to N \times \Jr.
\end{align}

\begin{Lemma}
  \label{lem:10}
  If $\arz \in T^*\Lambda_r N$ is a critical point for $S_{q,r}^{s'}$
  with critical value less than $-\norm{F}$, for some $s' \in ]5,s_r]$,
  then
  \begin{align*}
    (\pd{s} S^\sm_{q,r}(\arz)) (s') = -1.
  \end{align*}
\end{Lemma}

Intuitively this lemma says that the critical values of $S_{q,r}^s$
below $-\norm{F}$ is when $s\geq 5$ is increased simply translated
downwards with the same speed. Indeed, this should be considered a
fiber-wise version of Lemma~\ref{lem:closact}.

\begin{bevis}
  By Corollary~\ref{cor:critq} $\arz$ is the $\alpha^{s'}$-dissection
  of a time-1 flow curve for $H^{s'}$ with action less than
  $-\norm{F}$. By Lemma~\ref{actions} this flow curve lies in
  the set $T^*N-j(D_{t_{s'}}T^*L)$. By f6 above and how we constructed
  $H^s$ using $f_s$ this means that there is an $\epsilon>0$ such
  that: for any $z\in T^*N$ sufficiently close to the time-1 flow
  curve we have
  \begin{align*}
    H^s(z) = H^{s'}(z) + (s-s') \qquad \qquad s\in
    [s'-\epsilon,s'+\epsilon].
  \end{align*}
  This implies that this time-1 flow curve is a time-1 flow curve for
  all $s\in [s'-\epsilon,s'+\epsilon]$ with the action translated by
  $s'-s$.

  This does not mean that $\arz$ is a critical point for all $s$ in
  this interval. Indeed, the subdivision $\alpha^s$ is not necessarily
  constant so the dissection could change. So for each
  $s\in[s'-\epsilon,s'+\epsilon]$ we define $\arz_s\in T^*\Omega_{q,r}
  N$ to be the $\alpha^s$-dissection of this common time-1 flow
  curve. This is smooth in $s$ because $\alpha^s$ and the flow curve
  is. By Corollary~\ref{cor:critq} $\arz_s$ is now a critical point
  for $S_{q,r}^s$ when $s\in [s'-\epsilon,s+\epsilon]$.

  The corollary also tells us that the critical value is the same as
  the action, which shows that
  \begin{align*}
    S_{q,r}^s(\arz_s) = S_{q,r}^{s'}(\arz) + (s'-s).
  \end{align*}
  This proves that the directional derivative of $S_{q,r}^\sm$ along
  the tangent of the curve $s\mapsto (\arz_s,s)$ is $-1$, and the fact
  that $\arz$ is a critical point for $S_{q,r}^{s'}$ implies that the
  directional derivative along the curve $s\mapsto (\arz,s)$ at the
  point $(\arz,s')$ is the same. Indeed, the tangent vectors of the
  two curves at the point $(\arz,s')\in T^\bullet \Omega_{q,r} N
  \times \Jr$ differ by a tangent which is 0 on the second factor.
\end{bevis}

The above lemma is the reason for the choice of $a$ to be less
than $-\norm{F}$. Indeed, this gives us the control we need when $s$
tends to infinity. The general idea of why such a lemma provides
ex-continuation-maps is described intuitively in the introduction of
Section~\ref{homoindex}, and the proof of existence is postponed to 
the proof of Lemma~\ref{lem:1}. The following paragraph is the
construction using Lemma~\ref{lem:1} of quotient ex-map $q_r$. These
are the usual quotient maps of Conley indices, but lifted to ex-maps
over $N$. So they are the fiber-wise continuation ex-maps that we are
after - modulo the fact that we also need the $TN$-suspensions from
Section~\ref{cha:fib-susp}.

Firstly we pick $r_0>0$ such that $s_{r_0}>5$, and consider any $r>
r_0$. Put $W=N$ and define $G\colon W\times I \to  N \times \Jr$ by
\begin{align*}
  G(q,t) = (q,(1-t)s_{r-1} + ts_r).
\end{align*}
Then Lemma~\ref{lem:10} tells us that the adjoint of $G$, say $g
\colon W\to C^\infty(I,N\times \Jr)$, is space-like with respect to
$S_r^\sm$ (see Definition~\ref{def:4}), and so Lemma~\ref{lem:1}
provides us with an ex-map (contractible choice)
\begin{align} \label{eq:20}
  q_{r-1} = p_0^1 \colon I_a^b(S_r^{s_{r-1}},Y_r^{s_{r-1}})_N \to
  I_a^b(S_r^{s_r},Y_r^{s_r})_N
\end{align}
over $N$. The map induced on the unparametrized indices by $q_{r-1}$ is
simply the usual quotient map of Conley indices combined with homotopy
invariance when changing $H$ and $\alpha$ (this uses the assumption
that $a$ is a regular value for the action). However, to constructed
this in a parametrized setting while changing the functions and
pseudo-gradients is delicate - simply because $a$ is not always
fiber-regular and we may thus have some non-generic fibers. This is
why we need the notion of parallel transport in Lemma~\ref{lem:1}.

We now define the $TN$-spectrum $\FL$ by
\begin{align*}
  \FL_r = \left\{ \begin{array}{ll} 
      N_N & r< r_0 \\
      \GF_a^b(H^{s_r},\alpha^{s_r})_r = I_a^b(S_r^{s_r},Y_r^{s_r})_N &
      r \geq r_0
    \end{array}
  \right.
\end{align*}
The first $r_0$ structure maps are trivial ($N_N$ is initial in
$\Com_N^+$). For $r\geq r_0$ we define $\sigma^\FL_r \colon \FL_r^{TN}
\to \FL_{r+1}$ as the composition of the two maps:
\begin{itemize}
\item {The generating functions structure map 
    \begin{align*}
      \sigma_r \colon \GF_a^b(H^{s_r},\alpha^{s_r})_r^{TN} \to
      \GF_a^b(H^{s_r},\alpha^{s_r})_{r+1},
    \end{align*}
    which were also the ex-maps $\sigma_r$ defined in
    Section~\ref{cha:fib-susp}. }
\item{The ex-map from Equation~\eqref{eq:20} above
    \begin{align*}
      q_r \colon \GF_a^b(H^{s_r},\alpha^{s_r})_{r+1} \to
      \GF_a^b(H^{s_{r+1}},\alpha^{s_{r+1}})_{r+1}.
    \end{align*} }
\end{itemize}

Although we will not need it in the proof of the main theorem - we
will in Appendix~\ref{cha:hominv} discuss the fact that for all
possible choices made we get canonically (up to homotopy) ex-homotopy
equivalent instances of $\FL$ - meaning that they define an object in
the ex-homotopy category up to unique isomorphism. The only thing that we
will need for the main theorem is the following computation of the
homology given some homological assumptions on $L \subset T^*N$.

\begin{Corollary}
  \label{cor:2}
  Let $\F$ be $\Q$ or $\F_p$ for $p$ prime. If either
  \begin{itemize}
  \item{$N$ and $L$ are orientable and $j' \colon L\to N$ is relative
      spin, or}
  \item{$\F=\F_2$}
  \end{itemize}
  then the homology $H_*(\FL,\F)$ is isomorphic to the homology of
  the free loop space $H_*(\Lambda L,\F)$ with the usual
  component-wise regrading by the Maslov class.
\end{Corollary}

\begin{bevis}
  Since we assumed that $a$ is regular for $S_r^{s_r}$ for all
  $r\geq r_0$ we know from Lemma~\ref{lem:12} that the homology of
  $\GF_a^b(H^{s_r},\alpha^{s_r})$ is the reduced homology of the
  unparametrized index $I_a^b(S_r^{s_r},X_r^{s_r})$. Since $a$ is
  regular the usual homotopy invariance of such Conley indices
  applies. This was used in \cite{hejeh} to prove that this Conley
  index is stably (i.e. after some standard reduced suspensions)
  homotopy equivalent to some Thom-suspension of
  \begin{align*}
    \Lambda^\mu L,
  \end{align*}
  where $\Lambda^\mu L$ is the loop space of loops in $L$ with length
  less than $\mu$ and $\mu$ is the slope of the tangent of $f+s$
  intersecting the second axis at $-a$. These identifications were
  also shown to be compatible with the inclusion of loops when the
  slope $\mu$ were increased to $\mu'$. So that the structure maps on
  reduced homology from $(\FL_r)/N$ to $(\FL_{r+1})/N$ is the reduced
  homology (with the coefficients associated to the compatible Thom
  suspensions) of the usual inclusion from $\Lambda^\mu L$ to
  $\Lambda^{\mu'} L$.

  This provides the isomorphism for $\F=\F_2$ because here the
  Thom isomorphism always work - except that we should note that the
  Maslov index shift comes because the dimension of this bundle
  depends over each component on the Maslov index. For general $\F$ we
  may use the assumptions in the corollary. Indeed, the stable bundle
  over $\Lambda L$ constructed in \cite{hejeh} is identified with
  $j^*TN - TL +\eta$, where $\eta$ is the Maslov bundle introduced
  there. This is orientable when $\eta$, $TN$ and $TL$ are all
  orientable. So by the assumptions in the corollary we may focus
  solely on the first Stiefel Whitney class of the Maslov bundle
  $\eta$.

  The Maslov bundle $\eta$ was constructed in \cite{hejeh} by using a
  canonical lift of the map $f\colon L \to B\Or$, which classifies the
  virtual bundle $TN-TL=j^*TN-TL$, to a map $f' \colon L \to \U/\Or$,
  and then looping this to get a map
  \begin{align} \label{eq:2}
    \Lambda L \to \Lambda (\U/\Or) \to \Omega (\U/\Or) \simeq \Z \times
    B\Or. 
  \end{align}
  here the middle map is the canonical projection using that $\U/\Or$
  is a based loop space (in fact an infinite loop space), and the
  latter map is Bott periodicity as described in \cite{MR0163331}. The
  question of orientability of $\eta$ comes down to whether or not
  there is an element $\gamma \in \pi_1(\Lambda L)$ which is mapped to
  the non-trivial class in any of the components of $\Omega(\U/\Or)$,
  and we see that this is equivalent to the existence of a class
  $\alpha \in \pi_2(L)$ mapping to the non-trivial element in
  $\pi_2(\U/\Or) \cong \Z_2$. A short analysis of the long exact
  sequence of homotopy groups for $O\to U \to U/O$ shows that it must
  end with
  \begin{center}
    \begin{tabular}{c | c c c c c c c}
      4 & 
      $\xrightarrow{}$ & $0$
      & $\xrightarrow{}$ & $0$
      & $\xrightarrow{}$ & $0$ & $\xrightarrow{}$ \\
      3 & 
      $\xrightarrow{}$ & $\Z$
      & $\xrightarrow{\cdot 2}$ & $\Z$
      & $\twoheadrightarrow$ & $\Z_2$ & $\xrightarrow{}$ \\
      2 & 
      $\xrightarrow{}$ & $0$
      & $\xrightarrow{}$ & $0$
      & $\xrightarrow{}$ & $\Z_2$ & $\xrightarrow{g}$ \\
      1 & 
      $\xrightarrow{g}$ & $\Z_2$
      & $\xrightarrow{}$ & $\Z$
      & $\xrightarrow{\cdot 2}$ & $\Z$ & $\twoheadrightarrow$ \\
      0 & 
      $\twoheadrightarrow$ & $\Z_2$
      & $\xrightarrow{}$ & $0$
      & $\xrightarrow{}$ & $0$ & $\xrightarrow{}$ \\
      \hline
      n & & $\pi_n(\Or)$ & & $\pi_n(\U)$ & & $\pi_n(\U/\Or)$ \\
    \end{tabular}
  \end{center}
  In particular $g$ is an isomorphism. So orientability of $\eta$ is
  equivalent to whether or not the second Stiefel-Whitney class
  of the bundle classified by $f$ vanishes on $\pi_2(L)$. This is less
  than relative spin.
\end{bevis}


\chapter{Parallel Transport in Parametrized Conley
  Indices}\label{homoindex}

In this section we introduce the notion of parallel transport in
parametrized Conley indices, which is the notion we use to construct
the ex-maps $q_r$ from Equation~\eqref{eq:20}. Heuristically these are
a ``continuous combination'' of homotopy invariance and the natural
quotient maps
\begin{align*}
  I_a^b(f,Y)_M \to I_{a'}^b(f,Y)_M,
\end{align*}
where $a<a'<b$. We could also simultaneously consider varying $b$ and
the natural inclusion/quotient maps:
\begin{align*}
  I_a^b(f,Y)_M \to I_{a'}^{b'}(f,Y)_M,
\end{align*}
However, as we assume that $b$ is always an upper bound on all
fiber-critical values this will not be needed. So in the following we
assume that $b$ is always a fiber-regular value. We will illustrate
the idea by describing a simple case, and describe what assumptions we
need and what parallel transport then is (for the most general
description see Lemma~\ref{lem:1}). So assume that
\begin{itemize}
\item {The base $M$ is the unit interval $I$,}
\item {we have a smooth function $f\colon M'\to \R$, where $M'$ is
    the product of a manifold without boundary and $I$, and}
\item {we have a parametrized pseudo-gradient $Y$ on $M'$, with
    respect to the projection to one factor $M' \to I$, which
    satisfies the two properties C1  and C2 described in the
    introduction to Section~\ref{cha:pseudo}.}
\end{itemize}
This is equivalent to having a smooth homotopy $(f_s,Y_s),s\in I$ of
smooth functions and pseudo-gradients, satisfying C1 and C2 for all
$s\in I$. The property we need to be able to define parallel transport
forward in $I$ is: any critical points $x$ for $f_s$ (i.e. $(x,s)\in
M'$ is a fiber-critical point for the fiber over $s\in I$) with
critical value $a$ satisfies
\begin{align} \label{eq:13}
  \pare*{\pd{s} f}(x,s) < 0.
\end{align}
The parallel transport (or homotopy lift if the reader prefer) are
then maps
\begin{align*}
  p_s^{s'} \colon I_a^b(f_s) \to  I_a^b(f_{s'})
\end{align*}
for $0\leq s\leq s' \leq 1$ continuous in $s$ and $s'$ such that
$p_s^s=\id$. The existences of these maps in both directions when $a$
is regular for all fibers will describe mutual homotopy inverses and
thus recovers the notion of homotopy invariance of Conley indices.

The condition in Equation~\eqref{eq:13} is rather intuitive: if we
encounter a critical point with value $a$ when increasing $s$, it 
must cross downwards through $a$ such that going forward corresponds to a
usual collapse map of indices. So in generic cases this is homotopic
to a composition of homotopy invariance when $a$ is regular and
collapse maps for the crossings. The trouble is non-generic cases and
the fact that we need these maps to exist in compact families -
i.e. as ex-maps.

For the general setup we assume that $(f,Y)$ is a smooth function and
parametrized pseudo-gradient with respect to the submersion
\begin{align}\label{eq:14}
  \pi \colon M' \to M,
\end{align}
of smooth manifolds as in Definition~\ref{def:6}. We as usual denote
the flow of $-Y$ by $\psi_t$. In addition we are as mentioned above
assuming the properties:
\begin{itemize}
\item[C1)]{The flow $\psi_t$ of $-Y$ exists for all time (future and
    past).}
\item[C2)]{Outside a compact set $Y(f)$ is bounded from below by a
    positive constant.}
\end{itemize}
We will in this section use the abbreviation
\begin{align*}
  I(f)_M = I_a^b(f,Y)_M.
\end{align*}

To describe the general property we need to be able to construct the
parallel transport map we need a few notions. Let $\gamma \colon I \to
M$ be a smooth path. We may pull back everything using
$\gamma$. I.e. we get the pull back diagram
\begin{align*}
  \xymatrix{
    M_I' = \gamma^* M' \ar[r] \ar[d]^{\pi_I} & M' \ar[d]^{\pi} \\ 
    I \ar[r]^{\gamma} & M.
  }
\end{align*}
where the two vertical maps are smooth submersions. The function and
parametrized pseudo-gradient $(f,Y)$ pulls back to a function and
parametrized (with respect to $\pi_I$) pseudo-gradient 
$(f \circ \gamma , \gamma^*Y )$. The manifold $M'_I$ is not in general
a product, but there is still a canonical way to ask whether or not
Equation~\eqref{eq:13} is satisfied. Indeed, we may pick a lift $Z
\colon M_I' \to TM_I'$ of $\pd{s}$ in the sense that
\begin{align*}
  \xymatrix{
    M'_I \ar[r]^-{Z} \ar[d]^{\pi_I} & TM'_I \ar[d]^{D\pi_I} \\
    I \ar[r]^-{\tpd{s}} & TI
  }
\end{align*}
commutes. Then $Z_x(f \circ \gamma)$ for $x\in M_I'$ does not depend
on the choice of lift provided $x$ is a fiber-critical point. So we
may define $(\pd{s})^x(f\circ \gamma)$ as this unique value when $x$
is a fiber-critical point. We can thus ask that
\begin{align} \label{parlifteq}
  (\pd{s})^x(f\circ \gamma) <0
\end{align}
for all fiber-critical points $x$ with critical value $a$.

\begin{Definition}
  \label{def:2}
  Any smooth curve $\gamma \colon I \to M$ is called
  \textbf{space-like} (with respect to $f$, $\pi$ and $a$) if the pull
  back above satisfies Equation~\eqref{parlifteq} for the
  fiber-critical points with critical value $a$. We denote the space
  of these paths by $\Gamma^a(f,\pi)$ and topologize it as a sub-space of
  $C^\infty(I,M)$ with the weak topology.
\end{Definition}

In the following $W$ is compact Hausdorff. In particular any map
$W \to \Gamma^a(f)$ have all derivatives continuous as maps from
$W\times I$.

\begin{Definition}
  \label{def:4}
  The adjoint $G\colon W\times I \to M$ of any map $W \to
  C^\infty(I,M)$ is called a \textbf{path-smooth}
  homotopy. Furthermore if the image of the adjoint is in
  $\Gamma^a(M,f)$ we call it a \textbf{space-like homotopy}.
\end{Definition}

For such homotopies we use the notation
\begin{align*}
  G_s = G(\sm,s) \colon W \to M
\end{align*}
for $s\in I$, and $C^\infty(I,M) \ni G(w) = G(w,\sm) \colon I \to M$
for $w \in W$.

Notice that since the flow of $Y$ is fiber-wise the properties C1 and
C2 is preserved under the pull back. In fact the Conley indices pulls
back as ex-spaces. Indeed, if $h\colon W' \to W$ is any map of compact
Hausdorff spaces then we have an induced pull back functor
\begin{align} \label{eq:8puba}
  h^* \colon \Com_W^+ \to \Com_{W'}^+
\end{align}
of ex-spaces. Indeed, this is simply defined as the usual pull back of
spaces, which has canonical sections and projections.

\begin{Lemma} \label{lem:1}
  Let $I(f)_M$ be a parametrized Conley index with all the assumptions
  above. Let $G\colon W \times I \to M$ be any space-like
  homotopy. Then there exist a \emph{contractible} choice of
  a continuous family of parallel transport ex-maps
  \begin{align*}
    p_s^{s'} \colon G_s^*(I(f)_M) \to G_{s'}^*(I(f)_M)
  \end{align*}
  over $W$ for $0\leq s\leq s' \leq 1$ such that
  $p_s^{s''}=p_{s'}^{s''}\circ p_s^{s'}$ and $p_s^s$ is the identity.
\end{Lemma}

By contractible choice we simply mean that all the choices we make to
construct these ex-maps form a contractible space, and that the family
of course depends continuously on these choices.

\begin{bevis}
  First we consider the case where $W$ is a point, and thus $G$ is a
  single smooth curve $G \colon I \to N$. As above we can use $G$ to
  pull everything back and have everything parametrized by $I$, and
  thus in this first case we may assume that $M=I$, and the path is
  simply the identity on $I$.
  
  Recall the subspaces
  \begin{align*}
    C_\tau = \psi_{\tau}(f^{-1}([a,b])) \cap f^{-1}([a,b])    
  \end{align*}
  and
  \begin{align*}
    \uC_\tau = C_\tau \cap f^{-1}(a)
  \end{align*}
  defined in Section~\ref{cha:paracon}. Assume that $\tau$ is large
  enough for this to be compact by Lemma~\ref{lem:7}. This is thus one
  of the possible definitions of $I(f)_I = C_\tau :_I \uC_\tau$. We
  see that
  \begin{align*}
    C_\tau=\{x\in M'\mid f(x)\geq a,f(\psi_{-\tau}(x)) \leq b\},
  \end{align*}
  which implies that the boundary satisfies
  \begin{align}\label{eq:21}
    \partial C_\tau \subset f^{-1}(a) \cup (f \circ
    \psi_{-\tau})^{-1}(b).
  \end{align}
  So let $g= f\circ \psi_{-\tau}$. Notice that $Y$ is also a
  parametrized pseudo-gradient for $g$ and $b$ is a fiber-regular
  value for $g$. 

  We define $D_s=C_\tau \cap \pi^{-1}(s)$ and $\uD_s = \uC_{\tau} \cap
  \pi^{-1}(s)$ for each $s\in I$ and we thus have
  \begin{align*}
    C_\tau = \bigcup_{s\in I} D_s, \qquad \textrm{and} \qquad
    \uC_\tau  = \bigcup_{s\in I} \uD_s.
  \end{align*}
  Moreover, we have that the fiber of $I(f)_I$ at $s\in I$ is
  $D_s/\uD_s$. The goal is thus to construct maps
  \begin{align*}
    p_s^{s'} \colon D_s/\uD_s  \to D_{s'} / \uD_{s'}
  \end{align*}
  satisfying the requirements.

  The idea is to flow the set $D_s, s\in I$ to the fiber over
  $s'\in [s,1]$, in such a way that we get maps induced on the
  quotients $D_s/\uD_s \to D_{s'}/\uD_{s'}$. In the following we
  make this idea precise.

  Start by choosing any vector field $Z \colon M' \to TM'$ lifting
  $\pd{s}$. To accommodate the general situation later we will only
  assume that $Z$ is continuous. The fact that the curve is space-like
  means that we have that 
  \begin{align*}
    Z_x(f)=(\tpd{s})^x f < 0
  \end{align*}
  for all $x\in M'$ such that $f(x)=a$ and $x$ is a critical point for
  $f_s = f_{\mid \pi^{-1}(s)}$.

  We will alter this $Z$ to another lift $Z'$ which has better
  properties. Indeed we will find a $c>>0$ such that the flow of the
  vector field $Z'=Z - cY \colon M' \to TM'$ defines a map from one
  fiber of the parametrized Conley index to the other. Notice that
  $Z'$ is indeed another lift of $\pd{s}$ since $Y$ is vertical. So
  the flow of $Z'$ for time $s'-s$ takes the fiber $M_s'=\pi^{-1}(s)$
  to the fiber $M_{s'}'$ (when it is defined). The point is that we
  may pick $c>>0$ such that
  \begin{itemize}
  \item {$Z'_x(f)<0$ for all $x\in C_\tau$ with $f(x)=a$, and}
  \item {$Z'_x(g)<0$ for all $x\in C_\tau$ with $g(x)=b$.}
  \end{itemize}
  Indeed, this is true by compactness of $C_\tau$ and in each separate
  case because:
  \begin{itemize}
  \item{When $f(x)=a$ we have either $Y(f)<0$ or by the space-like
      assumption $Y(f)=0$ and $Z(f) <0$.}
  \item {When $g=b$ then $Y(f)<0$ since $b$ is fiber-regular.}
  \end{itemize}
  
  Now let $s\in I-\{1\}$ be given. Then since $D_s\subset C_\tau$ is 
  compact there is an $\epsilon>0$ such that the flow (which is $C^1$)
  of $Z'$ is defined on $D_s$ for time $[0,\epsilon]$. We may shrink
  this $\epsilon$ and get that any point in a neighborhood of $g=b$ in
  $D_s$ is flown to the set where $g\leq b$. Similar for $f=a$. So
  by possibly shrinking $\epsilon$ again we may assume that the flow of
  $Z'$ for time $t \in [0,\epsilon]$ maps $D_s$ to
  \begin{align*}
    g^{-1}(]-\infty,b]) \cup f^{-1}(]-\infty,a]) = D_{s+t} \cup
    f^{-1}(]-\infty,a])
  \end{align*}
  and similarly $\uD_s$ to
  \begin{align*}
    f^{-1}(]-\infty,a]),
  \end{align*}
  and hence induces a map of quotients
  \begin{align*}
    p_s^{s+t} \colon D_s / \uD_s \to (D_{s+t} \cup
    f^{-1}(]-\infty,a])) / f^{-1}(]-\infty,a]) \cong D_{s+t} /
    \uD_{s+t}. 
  \end{align*}
  Obviously $p_s^s$ is the identity, and dividing the interval
  $[0,\epsilon]$ into smaller bits and composing the resulting maps
  yields the same map - so we may uniquely extend this to all $s'\geq
  s$ and the maps will then satisfy
  \begin{align*}
    p_{s'}^{s''} \circ p_s^{s'} = p_s^{s''}.
  \end{align*}
  There are no problems at $1\in I$ since we can argue backwards in
  time and find the $\epsilon$ above such that the maps exists locally
  for any $s\in [1-\epsilon,1]$ and any $1\geq s'\geq s$. 

  For any map $g\colon W \to \Gamma^a(f)$ the same argument
  extends. The pull back $M''$ defined topological by the pull back
  diagram 
  \begin{align*}
    \xymatrix{
      M'' \ar[r] \ar[d] & M' \ar[d] \\
      W\times I \ar[r] & M
    }
  \end{align*}
  is when restricted to $\{w\} \times I$ for each $w\in W$ naturally a
  manifold $M''_w$. Their tangent bundles define a vector bundle $E
  \to W\times I$ with a surjective map to the vector bundle $W\times TI
  \to W\times I$. We may lift $\pd{s}$ to a continuous section $Z$,
  and also the pull back of $Y$ defines a continuous section in
  $E$. By compactness we may now find $c>>0$ as above working for all
  $w\in W$ simultaneously. The flow of $Z'$ is then defined
  fiber-wise, but glues together to continuous ex-maps over $W$.

  Contractibility follows from the fact that lifts $Z'$ satisfying the
  needed equations above is a contractible space. Also one could
  consider $\tau$ as part of the chosen structure, but this is also
  contractible, and even though the conditions on $Z'$ depends
  on $\tau$ the combined set of choices is still contractible.
\end{bevis}


\chapter{3S Fibrancy of $\FL$}\label{cha:fib-fl}

In this section we define the notion of a $\beta$-spectrum over a
smooth compact manifold $M$ being 3S-fibrant. As we will see in
Section~\ref{serre} this will imply that the homology of the fibers is
a local coefficient system and that there is a Serre spectral
sequences, with second page isomorphic to the homology of $N$ 
with coefficients in this local coefficient system. This will converge
to the homology of the $\beta$-spectrum. It is convenient at this
point to introduce the version of morphisms we will use in the
categories $\Sp_W^\beta$. The fact that we restricted our attention to
spectra with levels given by compact Hausdorff ex-spaces makes
everything a little easier than they would be in a more general
setting.

The morphisms of $\beta$-spectra we will defined can be considered
analogous to the notion of maps of colimits over
$\N$ of finitely generated Abelian groups. Indeed, even though the
homology groups of a $\beta$-spectrum is not always finitely generated
at each level (e.g. the Hawaiian earring example) the morphisms of
our spectra will behave as if they were - due to compactness. In
particular this means that for any level, say $r$, of a
$\beta$-spectrum the morphisms has to be defined on the image
$\sigma_{r,r'}$ for some structure map with sufficiently large
$r'$. Hence the following definition, which really only defines a map
in the homotopy category.

A map of $\beta$-spectra $f=[(f_r,h_r,k_r)] \colon \As \to \Bs$, with
$\As,\Bs \in \Sp_W^\beta$, is an equivalence class of triples of
\begin{itemize}
\item{a sequence of increasing integers $k_r\geq r$,}
\item{a sequence of ex-maps
    \begin{align*}
      f_r \colon \As_r^{(k_r-r)\beta} \to \Bs_{k_r}
    \end{align*}
    in $\Com_W^+$, and}
\item{a sequence of ex-homotopies
    \begin{align*}
      h_r \colon \As_r^{(k_{r+1}-r)\beta}\times I \to
      \Bs_{k_{r+1}}
    \end{align*}
    between the maps $f_{r+1} \circ \sigma_r^\As$ and
    $\sigma^{\Bs}_{k_r,k_{r+1}} \circ f_r$ (here we have suppressed some
    $\beta$-suspensions on ex-maps).}
\end{itemize}
The equivalence relation is generated by:
\begin{itemize}
\item {Ex-homotopy equivalence on both sequences $f_r$ and $h_r$
    through ex-maps such that they satisfy the above.}
\item {Replacing $k_r$ with $k'_r \geq k_r$ and thus; compatibly
    $\beta$-suspending $f_r$ (and $h_r$) with $(k_r'-k_r)\beta$ and
    composing with $\sigma^\Bs_{k_r,k_r'}$.}
\end{itemize}
For $W$ a point this defines a morphisms in the usual homotopy
category of spectra (only for spectra with compact Hausdorff
levels - if not one has to be a little more careful).

Without the ex-homotopies $h_r$ this would not make much sense as maps
of colimits - i.e. homology would not be a functor. We could simply
have asked that the maps were ex-homotopic and not included the 
ex-homotopies $h_r$ in the structure. However, it is well-known that
this produces an extra unwanted equivalence on the morphisms, and
although this has little bearing on the main argument we include it so
that we really do get morphisms in a more conventional homotopy
category of $\beta$-spectra.

Let $g\colon W' \to W$ be any map between compact Hausdorff
spaces. The pull back $g^*S^\beta$ of $S^\beta$ is naturally  
identified with $S^{g^*\beta}$. Also the smash product over $W$ is
pulled back to the smash product over $W'$. So $g^*$ naturally defines
a functor from $\beta$-spectra over $W$ to $(g^*\beta)$-spectra over
$W'$.

Let $M$ be a closed neighborhood retract of some open manifold. Let
$\As$ be a $\beta$-spectrum in $\Sp_M^\beta$, and let 
\begin{align*}
  G \colon W \times I \to M
\end{align*}
with $W$ compact Hausdorff. Then we say that $G$ has a
\textbf{stable lift} to $\As$ if there is a map of $G^*\beta$-spectra 
\begin{align} \label{parmap}
  \overline{G} \colon (G_0^*\As) \times I \to G^*\As
\end{align}
over $W\times I$, where
\begin{itemize}
\item {$G_0 = G_{\mid W\times \{0\}}$,}
\item {we have identified $G_0^*\beta \times I$ with $G^*\beta$ in a
    way compatible with the canonical choice over $0\in I$, and}
\item { $\overline{G}$ is the canonical identification (or identity)
    over $0\in I$.}
\end{itemize}
The existence of such a stable lift has nothing to do with the
choice in the second bullet point. Indeed, any two choices are related
by an automorphism, which induce an automorphism of $G^*\As$ - which
is seen to be the identity because we only defined morphisms over
$W\times I$ up to (ex-)homotopy. Also this choice is a contractible
choice - since we required it to be the
canonical choice over $0\in I$. In fact, these identifications will
in the following construction come for free since we are using the
parallel transport from Lemma~\ref{lem:1}. Indeed, using the lemma on
a parametrized Conley index which happens to be the sphere bundle
$D\beta/S\beta$ (e.g. having a single non-degenerate critical point in
each fiber) produces maps that are ex-homotopic to standard
parallel transport in $\beta$ descended to the sphere bundle
$D\beta/S\beta$.

\begin{Definition} \label{3S}
  An ex-spectrum $\As\in \Sp_M^\beta$ is \textbf{3S fibrant}
  (smoothly and stably Serre fibrant) if it has stable lifts of all
  path-smooth homotopies $G \colon W\times I\to M$ with $W$ a compact 
  CW complex.
\end{Definition}

The $TN$-spectrum $\FL$ has this property mainly because of the
following lemma and Lemma~\ref{lem:10}. Let $H^s,\alpha^s$ for each 
$s\in \Jinf$ and $a,b,(s_r)_{r\geq r_0}$ be as in
Section~\ref{homoindex} defining the families of smooth finite
dimensional approximations $S_r^s$ (and $S_{q,r}^s$) and subsequently
defining an instance of $\FL$.

\begin{Lemma}
  \label{lem:3}
  For any $s\in \Jr$ let $\arz$ be a fiber-critical point for $S_r^s$
  with respect to $\ev_0 \colon T^*\Lambda_r N \to N$. I.e. a critical
  point for $S_{r,q_0}^s$. Then  
  \begin{align*}
    \norm{(\nabla S_r^s)_\arz} \leq 2.
  \end{align*}
\end{Lemma}

The proof of this is a simple reference to the gradient approximation
from Lemma~\ref{new42} of the finite dimensional
approximations. However, the action on the infinite dimensional space
of paths starting and ending in the same fibers behaves in the exact
same way. Indeed, the gradient of an open time-1 path for the
Hamiltonian flow only depends on the end-points and if we calculate 
the gradient keeping in mind that we are restricted to paths that
start and end in the same fiber within $DT^*N$ we get the exact same 
gradient and bound as in the proof below.

\begin{bevis}
  Corollary~\ref{cor:critq} tells us that the fiber critical points
  over $q\in N$ is precisely those for which $\epsilon_{q_j}=0$ for
  all $j$ and $\epsilon_{p_j}=0$ for all $j\neq 0$. The gradient
  estimates in Lemma~\ref{new42} then tells us that
  \begin{align*}
    \norm{\nabla S_r^s} = \norm{ \nabla_{q_0} S_r^s} =
    \norm{\epsilon_{p_0}} \leq 2.
  \end{align*}
  The latter because all our Hamiltonians are such that $S_r^s$ only
  has fiber-critical points with all $z_j$ in $DT^*N$. In particular
  $p_0 \in DT^*N$ and $p_0^- \in DT^*N$.
\end{bevis}

To prove the 3S fibrancy structure of $\FL$ we need to
understand how to create ex-homotopies by ``varying'' the jump points
$s_r$ where we Thom-suspend using the map from
Section~\ref{cha:fib-susp} and combine this with the parallel 
transports from Section~\ref{homoindex}. The best way to understand
this concept is proving the following concrete lemma.

\begin{Lemma}
  \label{lem:exhomjump}
  If we replace $s_{r+1}$ with $s_{r+1}' \in [s_r,s_{r+1}]$ in the
  definition of $\FL$ we get the same composed structure map
  $\sigma^\FL_{r,r+2}$ up to ex-homotopy.
\end{Lemma}

When we defined $\FL$ we did not actually allow $s_r=s_{r+1}$. However,
this was only to exclude the possibility of having to do parallel
transport along constant paths (see definition of $q_r$ from
Equation~\eqref{eq:20}) and because it will be convenient later when
considering products. It is convenient to allow the case of
$s_r=s_{r+1}$ in this section and the argument in the following proof
makes it very clear how this can be handled in more generality.

\begin{bevis}
  The definition of the composed structure map $\sigma_{r,r+2}^{\FL}$
  depends on contractible choices. We will make these choices now such
  that the similar structure maps, say $\sigma_{r,r+2}^{s_{r+1}'}$,
  defined using $s_{r+1}'$ instead of $s_{r+1}$ is continuous in
  $s_{r+1}'$. Indeed, let
  \begin{align*}
    P_s^{s'} \colon I_a^b(S_{r+2}^s,Y_{r+2}^s)_N \to 
    I_a^b(S_{r+2}^{s'},Y_{r+2}^{s'})_N
  \end{align*}
  be parallel transport maps for $s\leq s' \in [s_r,s_{r+2}]$ defined
  by the obvious space-like homotopy $N\times I \to N\times
  [s_r,s_{r+2}]$ (similar to those defining $q_{r+1}$ in
  Equation~\eqref{eq:20}) and let
  \begin{align*}
    p_s^{s'} \colon I_a^b(S_{r+1}^s,Y_{r+1}^s)_N \to 
    I_a^b(S_{r+1}^{s'},Y_{r+1}^{s'})_N
  \end{align*}
  be defined similarly for $s\leq s' \in [s_r,s_{r+1}]$. These are
  illustrated in Figure~\ref{fig:parsusp}. 
  \begin{figure}[ht]
    \centering
    \includegraphics{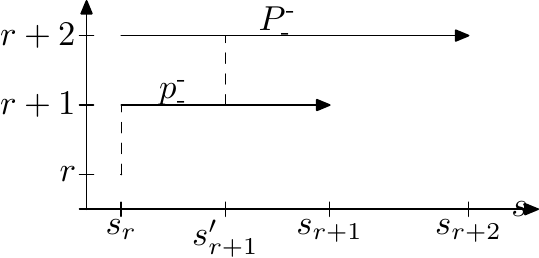}
    \caption{The horizontal lines represents parallel transport
      forward in $s$ and the vertical dotted line represents when we
      choose to use the suspension from Section~\ref{cha:fib-susp}.}
    \label{fig:parsusp}
  \end{figure}

  Then the structure maps may be defined by
  \begin{align*}
    \sigma_{r,r+2}^{s_{r+1}'} = P_{s_{r+1}'}^{s_{r+2}} \circ
    \sigma^{s_{r+1}'}_{r+1} \circ p_{s_{r}}^{s'_{r+1}} \circ \sigma_r,
  \end{align*}
  where $\sigma^{s_{r+1}'}_{r+1}$ is the suspension ex-map
  \begin{align*}
    I_a^b(S_{r+1}^{s_{r+1}'},Y_{r+1}^{s_{r+1}'})_N \to
    I_a^b(S_{r+2}^{s_{r+1}'},Y_{r+2}^{s_{r+1}'})_N,
  \end{align*}
  which are easily defined continuous in $s_{r+1}'$. Indeed, the
  suspension maps from Section~\ref{cha:fib-susp} can be defined
  continuously over any compact smooth family of Hamiltonians and
  sub-divisions.
\end{bevis}

\begin{Proposition}
  \label{prop:1}
  Any instance of $\FL$ is 3S-fibrant.
\end{Proposition}

The intuitive idea of the following proof is: Lemma~\ref{lem:3} tells
us that the movement of the fiber-critical points for fixed $s$ is
bounded by $2$ when following a unit length path in
$N$. Lemma~\ref{lem:10} tells us that the relevant critical points
moves down with speed 1 when $s$ increases by unit length. This means
that if we lift any smooth path in $N$ to $N\times [s_r,\infty]$ such
that the speed in the second factor is more than two times the speed
in the first factor 
then we get a space-like path. This is also the reason for the
suggestive name space-like, indeed, we are lifting the path to the
cone where the rate of change of these critical points is guaranteed to
be negative by Lemma~\ref{lem:10}.

\begin{bevis}
  Let $G \colon W\times I \to M$ be any path-smooth homotopy, and let
  $r\in \N_0$ be given. Then let $c>>0$ be a bound on the norm of the
  first derivative of all the paths $G(w) \in C^\infty(I,M)$. Pick an
  increasing sequence $k_r\geq r$ such that
  \begin{align} \label{eq:19}
    s_{k_r} - s_r > 2c.
  \end{align}
  We claim that there is a stable lift $\overline{G}=(f_r,h_r,k_r)$ as
  in Equation~\eqref{parmap} with this sequence $k_r$.
  
  The structure maps of the generating function spectrum
  $\GF(H^{s_r},\alpha^{s_r})$ (the ex-map from
  Section~\ref{cha:fib-susp}) are
  \begin{align} \label{eq:22}
    \sigma_{r,k_r} \colon \FL_r^{\beta(k_r-r)} =
      I_a^b(S_r^{s_r},Y_r^{s_r})_N^{(k_r-r)\beta} \to
      I_a^b(S_{k_r}^{s_r},Y_{k_r}^{s_r})_N. 
  \end{align}
  This does not map to the correct ex-space. Indeed, we need to map to
  the indices defined at $s_{k_r}$ (not $s_r$). But moreover we need
  it to be defined on the pull back using $G_0$ and ``lift'' the
  homotopy $G$ in the sense of Equation~\eqref{parmap}.

  We will generalize the construction of the structure maps using the
  bound $s_{k_r}-s_r>2c$ to construct space-like homotopies, which are
  equal to restrictions (and reparametrization) of $G$ in the $N$
  factor. Indeed, for each $r\geq r_0$ define
  \begin{align*}
    F=F^r \colon W \times I \times I \to N \times [s_r,s_{k_r}]
  \end{align*}
  by
  \begin{align*}
    F (w,t,\tau) = (G(w,\tau t) , (s_{k_r}-s_r) \tau + s_r).
  \end{align*}
  Then this is space-like (using the $\tau$-factor) with respect to 
  $S_{k_r}^\sm$ and $a$. Indeed, with $s=(s_{k_r}-s_r)\tau +s_r$ we
  see that at a fiber critical point (i.e. critical for $S_{q,r}^s$)
  we have
  \begin{align*}
    (\pd{\tau} F ) (w,t,\tau) \leq  2 c t + (s_{k_r}-s_r) (\pd{s}
    S_{k_r}^{\sm})(s).
  \end{align*}
  Indeed, this uses Lemma~\ref{lem:3}, the assumed bound on $G$, and
  the chain rule. Then using Lemma~\ref{lem:10} (remembering that
  $a<-\norm{F}$ and $s_r\geq 5$) and Equation~\eqref{eq:19} we may
  bound this to be strictly less than $2ct - 2c$ which is non-positive
  for all $t\in I$.

  Lemma~\ref{lem:1} thus gives us parallel transport ex-maps
  \begin{align*}
    p_0^1 &\colon F_0^*I_a^b(S_{k_r}^\sm,Y_{k_r}^\sm)_N =
    G_0^*I_a^b(S_{k_r}^{s_r},Y_{k_r}^{s_r})_N \times I \to \\
    & \to F_1^*I_a^b(S_{k_r}^\sm,Y_{k_r}^\sm)_N =
    G^* I_a^b(S_{k_r}^{s_{k_r}},Y_{k_r}^{s_{k_r}})_N,
  \end{align*}
  over $W\times I$, which when pre-composed with the map from
  Equation~\eqref{eq:22} defines the wanted lift $f_r$ for each
  $r$. However, we also need to construct the homotopies $h_r$ and
  argue that the morphisms over $0\in I$ is (homotopic to) the
  identity.

  Firstly we notice that the restriction of each $f_r$ to $0\in I$ is
  ex-homotopic to the structure maps (pulled back by $G_0$). Indeed,
  the homotopy for $t=0$ is a homotopy which is constant in the $N$ 
  factor. This means that the map is defined precisely as in the lemma
  above - except that we are considering more than two suspensions and
  the longer interval $[s_r,s_{k_r}]$. The natural generalization of
  the above lemma to several jump points $s'_{r+1} \leq \dots \leq
  s'_{k_r-1}$ with each $s_m' \in [s_r,s_{m}]$ shows that $f_r$
  restricted to $0\in I$ is defined precisely as $G_0^*\sigma_{r,k_r}$
  if we had $s_r=s_{r+1}'\cdots=s'_{k_r-1}$. Intuitively this means that
  we get an ex-homotopy between the two maps by sliding down the jump
  points in the interval to all be equal to $s_r$.
  
  A generalization of the same idea can be used to construct
  $h_r$. For this it is convenient to introduce some notions:
  \begin{itemize}
  \item {We say that two homotopies $G',G'' \colon W\times I \to N$
      are concatenable if $G'_1=G''_0$.}
  \item{When two space-like homotopies $G',G''$ are concatenable we
      define the concatenated parallel transport
      \begin{align*}
        p_0^s = \left\{ \begin{array}{cc}
            (p')_0^{2s} & s \in [1,1/2] \\
            (p'')_0^{2s-1} \circ (p')_0^1 & s\in [1/2,1]
          \end{array} \right.
      \end{align*}
      where $p'$ and $p''$ are parallel transports associated to $G'$
      and $G''$ respectively.} 
  \item{We may generalize this to any number of concatenable
      homotopies, and in fact in general get parallel transport for
      any piecewise space-like homotopy.}
  \item{The choice of such parallel transports are contractible even
      when considering different ways of considering the paths
      piece-wise space-like. Indeed, the choice of a parallel
      transport is contractible on each piece, and introducing a
      redundant division point we may consider the usual parallel
      transport as defined by the concatenation above.}
  \end{itemize}
  We now realize that the two compositions $f_{r+1} \circ \sigma_r$
  and $\sigma_{k_r,k_{r+1}} \circ f_r$ are in fact both defined in the 
  following general way:
  \begin{itemize}
  \item{We have given piece-wise space-like homotopies in $N\times
      [5,\infty[$, which are space-like because their derivative lie
      in the cone given by $\norm{\gamma'_N(t)}\leq 2
      \gamma'_s(t)$. Here the sub-script denotes the components
      in $TN$ and $T[5,\infty[$ respectively.}  
  \item{We then parallel transport along these homotopies and at
      several given values of $s$ we stop the parallel transport 
      and suspend to a higher level and continue the parallel
      transport on this higher level.}
  \end{itemize}
  The generalization of the above lemma tells us that we can ignore
  where we do the suspensions. Indeed, we may homotope them all to
  happen at the beginning - i.e. at $s_r$. The space-like homotopies
  used to define these two different maps are homotopies between the
  same two maps - i.e. the maps $W\times I \to N$ defined using $G_0$
  (constant in the $I$ factor) and all of $G$. Figure~\ref{Ghoms}
  \begin{figure}[ht]
    \centering
    \includegraphics{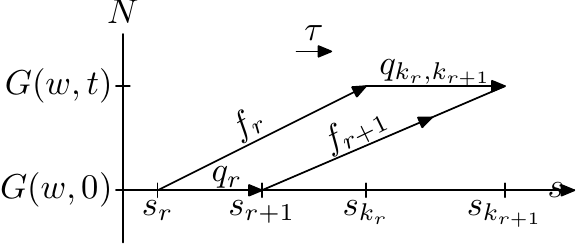}
    \caption{Paths in $N\times [0,\infty[$ defined by homotopies for a
      fixed $(w,t)\in W\times I$. The paths are labeled by which
      ex-map they induce, and $q_{k_r,k_{r+1}}=q_{k_{r+1}-1}\circ
      \cdots \circ q_r$.} 
    \label{Ghoms}
  \end{figure}
  illustrates the two different paths taken for some fixed $(w,t)\in
  W\times I$ (by construction we see that all these paths are constant
  in $N$ when $t=0$). We can interpolate between these two homotopies
  and stay within the space-like cone, and we may use parallel
  transport on the family of such an interpolation to construct the
  ex-homotopies $h_r$.
\end{bevis}


\chapter{Serre's Spectral Sequence and 3S fibrancy} \label{serre}

In section~\ref{cha:fib-fl} we introduced the notion of a
$\beta$-spectrum $\As\in \Sp_N^\beta$ being 3S fibrant over $N$. In
this section we describe the Serre spectral sequence converging to the
homology of such a $\beta$-spectrum. In particular we show that page 2
is the homology of $N$ with coefficients in the local coefficient
system given by the homology of the fibers. Some may want to call this
an Atiyah Hirzebruch spectral sequence. However, since we are
calculating standard singular homology with coefficients $\F$ of a
fibered spectrum, and not thinking of this as homology of the base
with coefficient in a twisted homology theory - we have kept the
Serre terminology.

Firstly - using Appendix~\ref{cha:funcrep} - we will assume that
$\beta$ comes with a trivialization $\beta \cong \R^l \times
N$. Secondly, the following lemma tells us that we need not worry that
we can only lift path-smooth homotopies.

\begin{Lemma}
  \label{lem:4}
  Let $W$ be compact, any homotopy $H\colon W\times I \to N$ is
  homotopic rel $W\times \partial I$ to a homotopy which is
  path-smooth. If the original homotopy was constant on a
  subspace $A\subset W$ we may assume that the new path-smooth
  replacement is also constant on $A$.
\end{Lemma}

\begin{bevis}
  One can compose smooth paths and get smooth paths if one
  is willing to reparametrize them to make them constant close to the
  endpoints. So by subdividing the homotopy using compactness we
  may assume that each path $h(x,\sm)$ is contained within an $\epsin$
  neighborhood of $h(x,0)$ for all $k\in K$. Now the family of 
  paths $\gamma(q,q',t)$ given by the unique shortest geodesic between
  two points $q$ and $q'$ in $N$ with distance less than $\epsilon_0$
  are smooth in all three variables. This means that all the
  derivatives with respect to $t$ (defined locally in charts) are
  continuous as functions of $(q,q',t)$. So let $\pi_W
  \colon W\times I \to W$ and $\pi_I \colon W\times I \to I$ be the
  projections then the homotopy defined by $\gamma(h^0\circ
  \pi_W,h^1\circ \pi_W,\pi_I) \colon W \times I \to N$ is adjoint to a
  map $W \to C^\infty(I,N)$ and is homotopic to $h$ rel
  $W\times \partial I$.
  
  The last statement in the lemma follows by construction.
\end{bevis}

As in Appendix~\ref{cha:funcrep} we let $l\in \N$ denote the trivial
$l$-dimensional bundle over any space.

\begin{Lemma}
  \label{lem:2asd}
  All fibers of a 3S fibrant $l$-spectrum $\As \in \Sp_N^l$ are
  homology equivalent and their homologies define a local 
  coefficient system on $N$.
\end{Lemma}

\begin{bevis}
  This is standard - except that from the definition of 3S fibrant we
  only have a lifting property for smooth paths and path-smooth
  homotopies of these smooth paths. However, we simply note that
  defining the fundamental groupoid using smooth paths and path-smooth
  homotopies yields the same groupoid as the continuous construction -
  essentially because of Lemma~\ref{lem:4}. 
\end{bevis}

We will need to get a good hold on products and suspensions on the
level of chains. So we need the Eilenberg-Zilber operators. Let $A$
and $B$ be spaces. Define the Eilenberg-Zilber operator
\begin{align*}
  P^{n,m} \colon C_n(A)\otimes C_m(B) \to C_{n+m}(A\times B)
\end{align*}
by subdividing $\Delta^n\times \Delta^m$ in the following way: let
$\Delta^n=[x_0,\dots,x_n]$ and let $\Delta^m=[y_0,\dots,y_m]$ then 
each subset $S \subset \{1,\dots,n+m\}$, with $\absv{S}=n$ defines a
sequence of cross-products of 0-simplices in $\Delta^n\times
\Delta^m$ by starting with $(x_0,y_0)$ and then if 1 is in $S$ the
next point is $(x_1,y_0)$ but if 1 is not in $S$ it is $(x_0,y_1)$
continuing like this the $k^{\textrm{th}}$ point is
$(x_{\absv{S\cap\{1,\dots,k\}}},y_{k-\absv{S\cap\{1,\dots,k\}}})$,
this defines a non-degenerate linear $n+m$ simplex in $\Delta^n
\times \Delta^m$, and these subdivide the product into $n$ choose
$n+m$ simplices. We then define $P^{n,m}$ on the tensor product of
generators $\alpha\colon \Delta^n \to A$ and $\beta\colon \Delta^m
\to B$ as the sum over $\alpha\times\beta \colon \Delta^n\times
\Delta^m \to A\times B$ pre-composed with each of these times the sign
on each that corresponds to preserving orientation.

The Eilenberg-Zilber operators are strictly associative since we may
index similar non-degenerate simplices in the triple product
$\Delta^n\times \Delta^m \times \Delta^k$ by two disjoint subsets
$S_1,S_2 \subset \{1,\dots,n+m+k\}$ of size $\absv{S_1}=n$ and
$\absv{S_2}=m$, and the sign is still the orientation preserving
sign. So this defines  
\begin{align*}
  P^{n,m,k} \colon C_n(A)\otimes C_m(B) \otimes C_k(C) \to
  C_{n+m+k}(A\times B\times C),
\end{align*}
and similar for more than 3 factors.

\begin{Remark}
  \label{rem:5sa}
  By keeping track of the orientations one may check that the
  Eilenberg-Zilber operator satisfies the derivation property
  \begin{align} \label{deriv}
    \partial P^{n,m}(\alpha\otimes \beta) = P^{n-1,m}(\partial \alpha
    \otimes \beta) + (-1)^{nm}P^{n,m-1}(\alpha \otimes \partial
    \beta),
  \end{align}
  which implies that they induce maps on homology
  \begin{align*}
    H_*(A)\otimes H_*(B) \to H_*(A\times B),
  \end{align*}
  which we will use to define products, and also to systematically
  treat the suspensions. 
\end{Remark}

\begin{Remark}
  \label{rem:6}
  We will also need the fact that the Eilenberg-Zilber operators are
  strictly (graded) commutative. That is: the diagram
  \begin{align}
    \label{eq:3}
    \xymatrix{
      C_n(A)\otimes C_m(B) \ar[r]^{P^{n,m}} \ar[d]^{(-1)^{nm}\tau} &
      C_{n+m}(A\times B) \ar[d]^{T_*} \\
      C_m(B)\otimes C_n(A) \ar[r]^{P^{m,n}} &
      C_{m+n}(B\times A)      
    }
  \end{align}
  commutes, where $\tau$ and $T$ are the obvious twists. The sign is
  obvious by considering that the twist on the product $\Delta^n
  \times \Delta^m$ is not always orientation preserving.
\end{Remark}

All of the above can be done with coefficients in any ring, and we may
now formulate the following stable (in the fiber direction) version of
the Serre spectral sequence.

\begin{Proposition} \label{Spectral}
  Let $\F$ be any coefficient ring. For a 3S fibrant $\As\in
  \Sp_N^\beta$, there is a 1. and 4. quadrant spectral sequence
  $\{E_{n,m}^r,d_r\}$ such that
  \begin{itemize}
  \item{The spectral sequence strongly converges to a filtered
      quotient of $H_*(\As;\F)$.}
  \item{$E_{n,m}^2 \cong H_n(N;H_m(\As_{\mid\bullet};\F))$, where
      $H_*(\As_{\mid\bullet};\F)$ is the graded local coefficient
      system defined by $\As$ from Lemma~\ref{lem:2asd}. Since $N$ is
      $d$ dimensional this implies that the spectral sequence
      collapses on the $(d+1)^{\textrm{st}}$ page.}
  \end{itemize}
\end{Proposition}

\begin{Remark}
  The proof follows several ideas from \cite{MR2039760},
  \cite{MR1867354}, \cite{HatchSpec}, and \cite{MR1793722}.
\end{Remark}

\begin{bevis}
  For notational purposes we will only consider $\F=\Z$. However, the
  general case is similar.

  We start by using the Eilenberg-Zilber operators to define the
  stabilization maps, used in Definition~\ref{def:3}, on the
  level of chains. Define the chain map on singular chains
  \begin{align}\label{eq:1023}
    \Sigma_* =\Sigma_{*,r} \colon C_*(\As_r,N) \to C_{*+l}(\As_{r+1},N)
  \end{align}
  as the composition of
  \begin{itemize}
  \item{The chain map
      \begin{align*}
        C_*(\As_r,N) \to C_{*+l}(\As_r \times S^l,N\times S^l\cup \As_r\times
        \{s_0\}), 
      \end{align*}
      given by sending $\alpha$ to $P^{*,l}(\alpha\otimes \beta)$
      where $\beta \in C_l(S^l)$ is some fixed representative of the 
      generator of $H_l(S^l,\{s_0\})$. This is a chain map due to the
      derivation property~\eqref{deriv} and the fact that $\beta$ is
      closed.}
  \item{The chain map induced by the quotient map 
      \begin{align*}
        (\As_r \times S^l,N\times S^l\cup \As_r\times \{s_0\}) \to
        (\As^\beta_r,N).
      \end{align*}}
  \item{The chain map induced by the structure map $\sigma_r$ of $\As$.}
  \end{itemize}
  Now define 
  \begin{align*}
    C_* = C_*(\As) = \colim_{r\to\infty}(C_{*+lr}(\As_r,N))
  \end{align*}
  in the category of chain complexes, using the chain maps
  $\Sigma_*$. Since limits of chain complexes commute with
  homology we have $H_*(C_*) \cong H_*(\As)$.

  We call a singular $m$-simplices $\alpha \colon \Delta^m \to \As_r$
  $n$-degenerate over $N$ if the projection $p\circ \alpha$ is a
  singular simplex in $N$ which is a degeneration of a simplex of
  dimension less than or equal to $n$. I.e. there is a commutative
  diagram
  \begin{align*}
    \xymatrix{
      \Delta^m \ar[r] \ar[d]^L & \As_r \ar[d]^{p_{\As_r}}\\
      \Delta^n \ar[r] & N,
    }
  \end{align*}
  where $L$ is a linear degeneration. So $L$ is given by an order
  preserving map $\{0,\dots,m\}$ to $\{0,\dots,n\}$. Let
  $F_nC_m(\As_r,N)$ be the free Abelian group generated by the
  $n$-degenerate $m$-simplices in $\As_r$ quotiented by the free
  subgroup generated by those in $N$ (or more precisely in the image
  of the section of $\As_r$). Since the degeneracy property is
  preserved when taking boundary we see that $F_nC_*(\As_r,N)$ is a
  sub chain complex. We also see that the chain map $\Sigma_*$
  preserves
  this property since the suspension is in fiber direction. Indeed in
  general if $\alpha \in C_n(\As_r)$ and $\beta\in C_m(Y)$ then
  $P^{n,m}(\alpha,\beta)$ is $n$ degenerate when projected to
  $\As_r$. So we may define the $n^{\textrm{th}}$ filtration of $C_*$
  by
  \begin{align*}
    F_nC_* = \colim_{r\to\infty} (F_nC_{*+lr}(\As_r,N)) \subset C_*.
  \end{align*}

  Since $F_nC_* = 0$ for $n<0$ and $\cup_n F_nC_* = C_*$ we get a
  filtration of $H_*(C_*)$ as $0\subset H_*(F_0C_*)\subset
  H_*(F_1C_*)\subset \cdots$, giving a spectral sequence
  $(E_{n,m}^*,d_*)$ converging as described in the first bullet
  point. Here the initial exact triangle $A^1 \to A^1 \to E^1
  \to A^1$ is given by $A^1_{n,m}=H_{n+m}(F_nC_*)$ and
  $E^1_{n,m}=H_{n+m}(F_nC_*/F_{n-1}C_*)$, with the usual maps.

  To get the second bullet point we need to inspect the first page and
  its differential $d_1 \colon E_{*,*}^1 \to E^1_{*-1,*}$. First
  notice that any $n$-degenerate simplex, which is not
  $(n-1)$-degenerate has a unique non-degenerate $n$ simplex in $N$ over
  which it is non-degenerate. This implies that the quotient
  $F_nC_*/F_{n-1}C_*$ is given by a direct sum as chain complexes over
  all the non-degenerate $n$ simplices in $N$. Let $\alpha\colon
  \Delta^n \to N$ be such a non-degenerate simplex, and let
  $F_s^{\alpha}C_*$ denote the subspace in $F_s C_*$ spanned by those
  $s$-degenerate simplices which if projected to $N$ is a composition
  of face and degeneracies of $\alpha$. Then the direct summand of
  the quotient $F_nC_*/F_{n-1}C_*$ corresponding to $\alpha$ may be
  canonically identified with $F^\alpha_n C_*/F_{n-1}^{\alpha}C_*$.

  Let $h \colon \Delta^n \times I \to N$ be a homotopy from the
  constant simplex mapping to $q=\alpha(1,0,\dots,0)$ and $\alpha$ 
  relative to $\Delta^0=\{(1,0,\dots,0)\}\subset \Delta^n$. We may
  assume by Lemma~\ref{lem:4} that $h$ is path-smooth and thus
  invoke Definition~\ref{3S} and get a map (when evaluating at $1\in
  I$) of ex-spectra
  \begin{align*}
    \Delta^n \times \As_{\mid q} \to \alpha^*\As
  \end{align*}
  over $\Delta^n$. By using the homotopy in the opposite direction we
  may get a map of ex-spectra in the other direction
  \begin{align*}
    \alpha^*\As \to \Delta^n \times \As_{\mid q}.
  \end{align*}
  By usual arguments (using a null-homotopy of concatenated homotopies)
  these two maps are homotopy inverses as $(\alpha^*\beta)$-spectra
  over $\Delta^n$ (since we are only considering maps up to homotopy
  this means it is an isomorphism in the category). So they induce a
  homology isomorphism, but that is not precisely what we are after:
  looking at the universal property for pullbacks:
  \begin{align*}
    \xymatrix{
      \Delta^{N} \ar@/^/[drr] \ar@/_/[ddr] \ar[dr] & \\
      & \alpha^*\As \ar[r] \ar[d] & \As \ar[d] \\
      & \Delta^n \ar[r] & N
    }
  \end{align*}
  we realize that the chain complexes $F^\alpha_sC_*$ is isomorphic to
  the chain complex given by the span of those simplices in
  $\alpha^*\As$ which projects to actual simplices in the simplicial
  structure of $\Delta^n$ which are $s$-degenerate. However, the maps
  also induce a chain homotopy equivalence of such restricted chain
  complexes simply because the maps \emph{and homotopies} are over
  $\Delta^n$. This means that we may identify the homology of the
  quotient complex $F_n^\alpha C_*/F_{n-1}C_*$ with
  $H_n(\Delta^n,\partial\Delta^n; H_*(\As_q))$. In fact we may do this
  such that $d_1$ is compatible with the usual boundary operator on
  $\alpha$. We have thus argued that $(E_{*,m}^1,d_1)$ is chain
  homotopy equivalent to the \emph{normed} chain complex
  $C_*(N;H_m(\As_{\mid \bullet}))$, which establishes the second bullet
  point.
\end{bevis}


\chapter{The Inclusion of Constant Loops} \label{cha:inc-const}

The identification of the homology in Corollary~\ref{cor:2} tells us
that there is an inclusion of constant loops $H_*(L,\F) \to
H_*(\FL,\F)$. The goal of this section is to prove that this inclusion
defines a map of spectral sequences from the usual Serre spectral
sequence for a fibrant replacement of the map $L \to N$ into the Serre
spectral sequence associated to $\FL$ in Proposition~\ref{Spectral}. 

The main part of the section is to construct a map of pairs
\begin{align}
  \label{eq:9}
  i_r \colon (D V, S V) \to (\FL_r , N)
\end{align}
for each $r$ sufficiently large (depending on the instance of $\FL$).
Here $V \to L$ will (under the assumptions of Corollary~\ref{cor:2})
be an $\F$-oriented $rd$-dimensional vector bundle and the induced map
on homology
\begin{align*}
  H_*(L;\F) & \cong H_{*+rd}(DV,SV;\F) \xrightarrow{i_*} \\
   & \to H_{*+rd}(\FL_r,N;\F) \to H_*(\FL;\F) \cong H_{"*"}(\Lambda
   L;\F)
\end{align*}
will be the inclusion. Note that the latter isomorphism is graded
component wise by the Maslov index. However, the grading on the
component containing constant loops is, indeed, 0. Of course since we
want to involve the Serre spectral sequence associated to the fibrant
replacement of $L\to N$ we need $i_r$ to somehow relate to this
map. Indeed, we will construct $i_r$ such that restricting it to
the zero-section $L\subset DV$ and projected to $N$ (using projection
of $\FL_r$) it will be homotopic to the composition $p\colon L\to T^*N
\to N$.

\begin{Lemma} \label{lem:14as}
  A map $i_r$ as in Equation~\eqref{eq:9} with the above properties
  will induce a map from the Serre spectral sequence associated to the
  fibrant replacement of $p\colon L\to N$ into the Serre spectral
  sequence associated to $\FL$, which then realizes the filtered
  quotient of the inclusion of constant loops.
\end{Lemma}

\begin{bevis}
  Define $A=DV:^P_N SV$ by using the extension $P \colon DV
  \xrightarrow{i_r} \FL_r \to N$ of $p$. This is an ex-space over $N$
  but it is \emph{not} in general a sphere bundle, and not in general
  fibrant in any way. The map from $A \to \FL_r$ induced by $i_r$ is
  by definition of $A$ an ex-map over $N$. So if we fibrantly replace
  both sides as ex-spaces over $N$ we get maps of the Serre spectral
  sequences (defined relative to the image of the sections because we
  are talking about ex-spaces). Since the spectrum $\FL$ is already
  3S-fibrant there is a map from the target spectral sequence and into
  the spectral sequence associated to $\FL$ from
  Proposition~\ref{Spectral}. So all we need to do to prove the lemma
  is to identify the Serre spectral sequence associated to the fibrant
  replacement of $A$.

  When constructing $A$ and taking fibrant replacement we would get
  ex-homotopy equivalent results if we replaced the map to $N$ by a
  homotopy equivalent map. So since $P$ restricted to the zero section
  is $p\colon L \to N$ and the inclusion $L\to DV$ is a
  homotopy equivalence we may pick a homotopy from $P$ to $P'\colon
  DV\to L \xrightarrow{p} N$, where the first map is the projection to
  the base of the vector bundle. So the fibrant replacement of $A$ is
  ex-homotopy equivalent to the fibrant replacement of $A' =
  DV:_N^{P'} SV$ - and hence have the same Serre spectral sequence.

  One may check that taking a fibrant replacement $P_NL\to N$ of $L\to
  N$ (this time not as an ex-space) and pulling back $V$ to $V_P \to
  P_NL$ and constructing
  \begin{align*}
    PA' = DV_P :_N SV_P
  \end{align*}
  we get that $PA'$ is a fibrant replacement of the ex-space $A'$. 

  Using the fact that the Thom-isomorphism is natural with respect to
  restriction we then see that since $V_P$ is oriented the Serre
  spectral sequence for the fibrant replacement of $A'$ must be the
  usual one for the fibrant replacement of $p\colon L\to N$ but
  shifted up in the fiber degree by the dimension of $V$.
\end{bevis}

The rest of the section thus concerns the construction of
$i_r$. Even though we are trying to lift a known map (we will define it
properly below) as a map of pairs involving a parametrized Conley
index - this section is mostly an exercise in \emph{un}parametrized
Conley indices of the action associated to the Hamiltonians. Notice in
particular that often we will argue using the pseudo-gradients $X^s_r$
and not the parametrized pseudo-gradients $Y^s_r$. In the proof above
this is compensated by the fact that we were allowed to homotope
the projection $P$ and get the same result - so we allow some
non-fiber-wise flows in the construction.

Let $\FL$ be defined as in Section~\ref{theex}. In particular
$a<-\norm{F}$ and $s_r\in [0,\infty]$ is an increasing sequence. We
will need $r$ to satisfy
\begin{align} \label{eq:11}
  s_r>-a \qquad \textrm{and as usual} \qquad s_r>5.  
\end{align}
So fix $r$ for the rest of this section such that this is true. The
level $\FL_r$ is defined using the Hamiltonian $H^{s_r}$ and the 
sub-division $\alpha^{s_r} \in \Delta^{r-1}$, and so we need to
consider these more carefully.

To construct the map $i_r$ we use that the family $H^s, s\in
[\epsilon,s_r]$ defines a homotopy from $H^{s_r}$ to an almost
constant Hamiltonian $H^\epsilon$. We pick this $\epsilon>0$ small
enough such that
\begin{itemize}
\item[$\epsilon$1)] {$H^\epsilon$ is sufficiently $C^2$-close to $H_\infty$ for
    $S_1^\epsilon$ to be defined.}
\item[$\epsilon$2)] {$H^\epsilon$ has only constant periodic orbits - hence its
    critical set is $\{-\epsilon,0\}$.}
\item[$\epsilon$3)] {The interval $]-s, -s+\epsilon[$ is regular for $S_r^s$ for all
    $s\in [\epsilon,s_r]$.}
\item[$\epsilon$4)] {$-s_r + \epsilon/2 < a$.}
\end{itemize}
The first is simply H3 in Appendix~\ref{cha:pd-ham} and the fact that
$f_0=0$. The next two follows easily from the ``explicit''
construction of $f_s$ and $H_\infty$ in Appendix~\ref{cha:pd-ham}. The
last follows easily from Equation~\eqref{eq:11}.

\begin{figure}[th]
  \centering
  \includegraphics{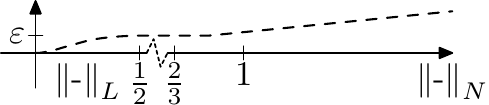}
  \caption{Hamiltonian $H^\epsilon$}
  \label{fig:Homham}
\end{figure}

If we then define the varying lower bound
\begin{align*}
  \ua(s) = -s +\epsilon/2
\end{align*}
then $\ua(s)$ (v for varying) is regular for $S_r^s$ and we thus get
by homotopy invariance that
\begin{align*}
  I_{-\epsilon/2}^b(S_r^\epsilon,X_r^\epsilon) =
  I_{\ua(\epsilon)}^b(S_r^\epsilon,X_r^\epsilon)
  \simeq
  I_{\ua(s_r)}^b(S_r^{s_r},X_r^{s_r}) =
  I_{-s_r+\epsilon/2}^b(S_r^{s_r},X_r^{s_r}). 
\end{align*}
Indeed, $b$ is a fixed upper bound on all possible critical
values for any $s$. We may then compose with the quotient map
\begin{align*}
  I_{\ua(s_r)}^b(S_r^{s_r},X_r^{s_r}) \to I_a^b(S_r^{s_r},X_r^{s_r}) =
  \FL_r /N.
\end{align*}
Indeed, by $\epsilon$4) above we have $\ua(s_r)<a$. As we will see
below the index $I_{-\epsilon/2}^b(S_r^\epsilon,X_r^\epsilon)$ is
homotopy equivalent to a Thom-suspension of $L$ and it is well-known
that the quotient  above recovers the inclusion of constant curves
into the loop space when passing to homology. So this map is the map
we need to lift to a map of pairs. The plan for the rest of the
section is:
\begin{itemize}
\item {First we prove a statement for
    $I_{\ua(\epsilon)}^b(S^\epsilon_1, X^\epsilon_1)$, which is very
    similar to what we want in the end for $\FL_r/N$.}
\item{Then we use the $TN$-suspensions in Section~\ref{cha:fib-susp} to get
    the same type of result for $I_{\ua(\epsilon)}^b(S^\epsilon_r,X^\epsilon_r)$.}
\item{Then we use the above homotopy invariance (with extra structure)
    to get the same for $I_{\ua(s_r)}^b(S^{s_r}_r,X^{s_r}_r)$.}
\item{Finally we take the quotient to get to
    $I_a^b(S^{s_r}_r,X^{s_r}_r) = \FL_r/N$.}
\end{itemize}
The idea is that we along the way have maps \emph{before} taking
quotients so that we in the end can lift this to a pair - and to get
this is a little subtle.

The first part is thus to explicitly describe a map which
induces a homotopy equivalence on quotients for $S_1^\epsilon$.

\begin{Lemma}
  \label{lem:1a}
  The embedding $j \colon DT^*L \to T^*N$ induces a homotopy
  equivalence
  \begin{align*}
    DT^*L/ST^*L \to I_{\ua(\epsilon)}^b(S^\epsilon_1,X^\epsilon_1).
  \end{align*}
\end{Lemma}

\begin{bevis}
   By construction the $S_1^\epsilon$ has critical set precisely $L$
   disjoint union with $(D_{2/3}T^*N-j(DT^*L))$ (both connected), and
   the critical values are $0$ on $L$ and $-\epsilon$ on the other
   component. Since $L$ is isolated from infinity by the other set
   this implies that $L$ consists entirely of local maxima. In fact we
   will show that in a neighborhood of $L$ $S_1^\epsilon$ is
   Morse-Bott, and then the lemma follows.

   To see that it is Morse-Bott we consider points
   $(p_0,q_0)\in \Lambda_1 T^*N=T^*N$ close to $L$ using the metric on
   $j(DT^*L)$ coming from $L$. So let $d=d(q_0,p_0)$ denote the
   distance to $L$ in $\norm{-}_L$. We may express $S_1^\epsilon$
   as a sum of two terms:
   \begin{align*}
     S_1^\epsilon(q_0,p_0) = -H^\epsilon(q_0,p_0) + \int_\gamma \lambda,
   \end{align*}
   where $\gamma \colon S^1 \to T^*N$ is the small contractible curve
   defined by the concatenation of the pieces:
   \begin{itemize}
   \item {The time-1 flow line starting at $(q_0,p_0)$.}
   \item {The horizontal geodesic in $T^*N$ (using its metric) back to
       a point in the fiber over $q_0$.}
     \item{The vertical geodesic back to the point $(q_0,p_0)$.}
   \end{itemize}
   Indeed, this is described in Section~\ref{finred}. It is not
   difficult to bound the symplectic area of this by the length of the
   first part (which is $2cd=2c(\epsilon)d$ by f2 in
   Appendix~\ref{cha:pd-ham}) squared times a constant $k$. So we in
   fact see that in a small neighborhood 
   \begin{align*}
     S_1^\epsilon(q_0,p_0) \leq -cd^2 + k(2cd)^2.
   \end{align*}
   Indeed $H^\epsilon$ is the distance squared times $c$ by $f2$ in
   Appendix~\ref{cha:pd-ham}. For small $c$ this proves that the
   Hessian of $S_1^\epsilon$ is negative definite in the normal bundle
   to $L$, and, indeed, since $c=c(\epsilon)>0$ and $c(0)=0$ we get
   this for $\epsilon$ small enough.
\end{bevis}

The next step is to use the iterations of the suspensions from
Section~\ref{cha:fib-susp}. To do this we define $V'$ as the vector
bundle over $DT^*L$ given by
\begin{align*}
  V' = (p^*TN)^{\oplus (r-1)},
\end{align*}
where $p$ is the composition $DT^*L \xrightarrow{j} T^*N
\xrightarrow{\pi} N$. Then define
\begin{align*}
  g \colon V' \to T^*\Lambda_r N
\end{align*}
by
\begin{align*}
  g(q_0,p_0,v_1,\dots,v_{r-1}) =
  (q_0,p_0,q_0,v_1,q_0,v_2,\dots,q_0,v_{r-1})
\end{align*}
This is the composition of the maps $f_1$, $f_2$, \dots $f_{r-1}$ from
Equation~\eqref{suspdef}. In the definition of the suspension maps
$\sigma_k$ we used the negative pseudo-gradient flow to move the
images into position to define a map of quotients - i.e. a map to the
index. This we do not need to do at each step and can wait until the
very end, and get an ex-homotopic (we in fact only care about
homotopic at this point) ex-map. Indeed, to see this take the
iterated suspension maps $\sigma_{1,r}$ and compose this with the
flow of $-Y_r^\epsilon$ for a very long time. This will imply that we
can by using a homotopy remove the flows inserted along the
way. Indeed we homotope by using the flow for a shorter and shorter
time. This we can do while all along defining a map of quotients. So
we pick an $R>0$ such that the image $g(S_RV')$ after composing with
the flow lies in $(S_r^\epsilon)^{-1}(]-\infty,a])$. We may identify
the disc $D_RV'$ (being a disc bundle over the disc bundle $DT^*L$) by
\begin{align*}
  l_R \colon D_RV' \cong DV,
\end{align*}
where $V \to L$ is the vector bundle defined by
\begin{align*}
  V=(p^*TN)^{\oplus (r-1)} \oplus TL.
\end{align*}
We have in fact almost proved the following lemma.

\begin{Lemma}
  \label{lem:1b}
  There is a map $DV \to T^*\Lambda_r N$, which induces a homotopy
  equivalence
  \begin{align*}
    DV / SV \to I_{\ua(\epsilon)}^b(S_r^\epsilon,X_r^\epsilon).
  \end{align*}
\end{Lemma}

\begin{bevis}
  As argued above the map $g$ after composing with a negative
  pseudo-gradient flow and taking the quotients is homotopic to: a
  Thom suspension of $j$ composed with the iterated suspension map
  $\sigma_{1,r}$, both of which are homotopy equivalences by
  Lemma~\ref{lem:1a} and Corollary~\ref{cor:ths}. This map is
  defined before taking the quotients and thus precomposing $g$ with
  $l_R$ above (and still composing with the flow) we get a map
  \begin{align*}
    DV \to T^*\Lambda_r N,
  \end{align*}
  and we notice that using $h_R$ we get $DV/SV \cong D_RV' / ((S_RV')
  \cup (D_RV')_{(ST^*L)})$, which when inspected is seen to be the
  quotient on which $g$ provides the homotopy equivalence.
\end{bevis}

The next and the most subtle step is to use the homotopy
$(H^s,\alpha^s)$ and homotope the map from $DV$ to induce a similar
homotopy equivalence to the index defined at $s=s_r$. Here we need the
parallel transport from Lemma~\ref{lem:1}, and since this assumes the
lower bound $\ua(s)$ to be constant it is convenient to shift the
Hamiltonians and approximations. This will also make the proof of
Lemma~\ref{lem:1d} easier. I.e. we define
\begin{align*}
  \oH^s   &= H^s + c - s \\
  \oS_r^s &= S_r^s -c +s \\
  \oa & = \ua(s) -c + s = - \epsilon/2 - c \\
  \ob(s) & = b -c + s.
\end{align*}
Here $c$ is chosen such that the asymptotic tangent goes through
$0$ for one $s$ - hence all $s$ because we used the same $H_\infty$
for all $s$. Notice that $\oS_r^s$ is indeed the
finite dimensional approximation associated to $(\oH^s,\alpha^s)$ and
$X_r^s$ is still a pseudo-gradient. We may ignore the fact that
$\ob(s)$ depends on $s$ because we may simply pick the common maximum
$\ob=\ob(s_r)$ for $s\in [\epsilon,s_r]$. Indeed, all we need is an
upper bound on all critical values.

The following lemma may be considered a lift of the parallel transport
construction in Lemma~\ref{lem:1} in Section~\ref{homoindex} to the
pair $(DV,SV)$.

\begin{Lemma}
  \label{lem:2}
  There is a homotopy $h_s \colon DV \to T^*\Lambda_r N$, $s\in
  [\epsilon,s_r]$ such that each $h_s$ induces a homotopy equivalence
  on the quotients
  \begin{align*}
    DV/SV \to I_{\oa}^{\ob}(\oS_r^s,X_r^s) \simeq
    I_{\ua(s)}^b(S_r^s,X_r^s).
  \end{align*}
\end{Lemma}

\begin{bevis}
  As mentioned earlier since $\oa$ and $\ob$ are regular for all $s\in
  [\epsilon,s_r]$ we get by homotopy invariance of Conley indices that
  all these indices are indeed homotopy equivalent. However, to get
  the map defined on $DV$ \emph{before we take the quotients} we need
  an elaboration on the parallel transport maps in
  Lemma~\ref{lem:1}. Indeed, for $\oa$ regular the parallel transport
  using the projection
  \begin{align*}
    T^*\Lambda_r N \times [\epsilon,s_r] \to [\epsilon,s_r]
  \end{align*}
  realizes the homotopy equivalences of the fibers of the parametrized
  index over $[\epsilon,s_r]$.

  We thus elaborate on the proof of Lemma~\ref{lem:1} by constructing
  the vector field $Z'$ appearing there with certain additional
  properties. In that proof $Y$ were the parametrized pseudo-gradient
  over the base, which in our case is simply the interval
  $[\epsilon,s_r]$, and so consistent with this we use the notation
  $Y=X_r^\sm$ viewed as a parametrized pseudo-gradient for $\oS_r^\sm$
  over $[\epsilon,s_r]$.

  Since $\oa$ is a regular value for all $t$ and Lemma~\ref{lem:comp}
  tells us that outside a compact set $Y(\oS_r^\sm)$ is bounded from
  below by a non-zero constant - we see that $Y(\oS_r^\sm)$ is in fact
  bounded from below by a positive constant $c'>0$ on the set
  \begin{align*}
    W_a = (\oS_r^\sm)^{-1}(\oa) \subset T^*\Lambda_r N \times
    [\epsilon,s_r].
  \end{align*}
  We then define
  \begin{align*}
    Z' = \pd{s} - k'(1+P)Y,
  \end{align*}
  for some large $k'\in \R$. In the proof of Lemma~\ref{lem:1}
  $\pd{s}$ were called $Z$ because we did not have a canonical lift 
  of $\pd{s}$ to the total space. Also the interval was $I$ not
  $[\epsilon,s_r]$ and we did not allow the factor in front of $Y$ to
  be unbounded, which we do here because recall that $P=\max_j
  \norm{p_j}_N$. We may pick this $k'>>0$ such that
  \begin{align}\label{eq:4bo}
    Z'_{\arz,s}(\oS_r^\sm) <  0 \qquad \textrm{for all} \qquad
    (\arz,s) \in W_a,
  \end{align}
  in addition to the needs in the proof of Lemma~\ref{lem:1}. Indeed,
  Lemma~\ref{lem:1d} below bounds the factor $\pd{s} \oS_r^\sm$ in
  such a way that it can be dominated by the factor
  $-k'(1+P)Y(\oS_r^\sm) \leq -k'(1+P)c'$.

  As in Lemma~\ref{lem:1} this means that when taking a point
  $(\arz,s)$ in the $C_\tau$ we have chosen to define the index (see
  Section~\ref{cha:paracon}) and flowing with $Z'$ then it leaves
  $C_\tau$ through the bottom $\uC_\tau$ (the exit set). In addition
  to this we have (in this special case) made sure that it stays
  out. Indeed, the bound in Equation~\eqref{eq:4bo} tells us that for
  \emph{any} point where $\oS_r^s(\arz)=a$ the flow will decrease the
  value.

  Moreover the argument in Lemma~\ref{lem:infflow} that the flow of
  $-Y$ exists for all time can be extended to the fact that the flow
  of $Z'$ exists precisely until we reach the end of the homotopy
  where $s=s_r$. This extension of Lemma~\ref{lem:infflow} requires
  the fact that any solution to $P' \leq c_1 P +c_2$ does not go to
  infinity in finite time.

  We may now take $h_\epsilon$ to be the map from Lemma~\ref{lem:1b},
  and we may then define $h_s$ as the map into $T^*\Lambda_r N \times
  \{\epsilon\} \subset T^*\Lambda_r N \times [\epsilon,s_r]$ given by
  $h_\epsilon$ composed with the flow of $Z'$ for time
  $s-\epsilon$. This will then induces the homotopy equivalences when
  going to the quotients as wanted for all $s\in
  [\epsilon,s_r]$. Indeed, it does so for $s=\epsilon$.
\end{bevis}

The next part of the construction is realizing that $h_{s_r}$ from the
lemma above actually induces a map to the \emph{parametrized} index.
Indeed, we simply do not take the full quotient but the quotient over
$N$ and get a map of pairs
\begin{align*}
  h_1 \colon (DV,SV) \to (I_{\oa}^{\ob}(\oS_r^{s_r},Y_r^{s_r})_N , N).
\end{align*}
Here the sphere is sent to $(\oS_r^1)^{-1}(]\infty,\oa])$ which is
collapsed to $N$ (instead of the point) when defining the parametrized
index. We may then compose $h_{s_r}$ with the quotient (over $N$)
\begin{align*}
  I_{\oa}^{\ob}(\oS_r^{s_r},Y_r^{s_r})_N =
  I_{\ua(s_r)}^{b}(S_r^{s_r},Y_r^{s_r})_N \to
  I_a^b(S_r^{s_r},Y_r^{s_r})_N 
\end{align*}
to define
\begin{align*}
  i_r \colon (DV,SV) \to (I_a^b(S_r^{s_r},Y_r^{s_r})_N,N) = (\FL_r,N),
\end{align*}
which then realizes the inclusion of the constant curves on the level
of homology.

We used the following technical lemma above.

\begin{Lemma}
  \label{lem:1d}
  There is a $k>0$ such that: for any $\arz \in T^*\Lambda_r N$ and
  any $s\in [\epsilon,s_r]$ we have
  \begin{align*}
    \absv{\pd{s} (\oS_r^\sm(\arz)) (s)} \leq k(P+1),
  \end{align*}
  where $P(\arz)=\max_j \norm{p_j}$.
\end{Lemma}

Note that the homotopy is a fixed homotopy for $s\in [\epsilon,s_r]$,
and the $k$ in the lemma highly depends on that. So if the domain of
$\oS_r^s$ were compact this would be trivial (and even $P$ would be
a bounded function, which of course it is not).

\begin{bevis}
  First notice that $\oS_r^s$ splits as a \emph{finite} sum
  \begin{align*}
    \oS_r^s(\arz) = \sum_{j\in \Z_r} f_j(q_j,p_j,q_{j+1},s),
  \end{align*}
  where
  \begin{align*}
    f_j(q,p,q',s) = \int_{\gamma^s_j} \lambda_N - \oH^s dt +
    p^-(\exp^{-1}_{q^-}(q')),
  \end{align*}
  and where $\gamma_j^s$ is the smooth homotopy of curves given by
  Equation~\eqref{eq:16} and the smooth homotopies $\alpha^s$ and
  $\oH^s$, and $(q^-,p^-)=\gamma^s(\alpha^s_j)$. We suppress that
  $(q^-,p^-)$ depends on $s$, but keep it in mind. We define this
  $f_j$ when
  \begin{align*}
    (q,p) \in T^*N \quad \textrm{and} \quad \dist(q',q)\leq \epsin
    \quad \textrm{and} \quad t\in I.
  \end{align*}
  This is slightly more than we need because $\oS_r^s$ is a sum of
  these with $\dist(q',q) < \epsin$, but we need the extra
  boundary for a compactness argument. Since $P\geq \norm{p_j}$ we may
  assume with out loss of generality when bounding the derivative
  $\pd{s} f_j$ that $P=\norm{p}$.  
  
  The rate of change $\pd{s} f_j$ of $f_j$ is bounded by compactness
  on the set given by $\norm{p}\leq 1$. Outside this set the
  integration term vanishes because of the adjustment we did to $H^s$
  to define $\oH^s$. Indeed, the action of \emph{any} flow curve
  outside $DT^*N$ is 0 for $\oH^s$. This implies that we need to bound
  the $s$ derivative of the function
  \begin{align*}
    g_j = f_{j\mid \{\norm{p}\geq 1\}}(q,p,q',s) =  p^-(\exp^{-1}_{q^-}(q'))
  \end{align*}
  by a constant times $P+1$. We, in fact, can bound it proportional to
  $p$. Indeed, this is easy since it scales in $p$. That is, for
  $\tau\geq 1$ we have
  \begin{align*}
    g_j(q,\tau p,q',s) = \tau p^-(\exp^{-1}_{q^-}(q'))
  \end{align*}
  because the Hamiltonian flow outside of $DT^*N$ is scale
  equivariant. So any bound on the set $\norm{p}=1$ scales with
  $\tau$ even if $p^-$ depends on $s$ - because it does so in the same
  scale equivariant way.
\end{bevis}


\chapter{The Fiber-wise Product}\label{cha:fib-w-prod}

In this section we will construct a fiber-wise product on (some
instances of) the $TN$-spectrum $\FL$ from
Section~\ref{theex}. That is, we will for some (not all) instances of
$\FL$ construct a product map of spectra.
\begin{align} \label{eq:8sder}
  \mu \colon \FL \wedge_N \FL \to \tFL
\end{align}
where $\tFL$ is another instances of $\FL$. I.e. $\tFL$ is defined
using another set of input data (see below). There are a couple of
very good reasons for using different instances of $\FL$ at this
point.
\begin{enumerate}
\item{This is how the maps naturally occur.}
\item{Had we insisted on using the same $\FL$ as the target we would
    have to ``stabilize'' the map using ideas from
    Appendix~\ref{cha:hominv}. Doing it this way we can actually get a
    morphism of $2TN$-spectra which has $k_r=r$ in the definition in
    Section~\ref{cha:fib-fl}.}
\item{All we really need for the main argument is that it realizes the
    intersection product on the constant loops on $L$, and this we
    only need homologically using the identifications in
    Corollary~\ref{cor:2}.}
\end{enumerate}
Smash products of spectra are rather delicate, and we will only define
a very naive version of these, and this will suffice for our
purpose.

\begin{Remark}
  \label{rem:1}
   We will not explicitly need the Chas-Sullivan product in this
   construction. However, we note that the construction here
   is based on it. In fact for $L=N$ (the zero section) we recover the
   spectral sequence from \cite{MR2039760} which converges
   to the Chas-Sullivan homology ring. This is not that surprising
   since below we are concatenating the loops in each fibers.
\end{Remark}

We define the smash product $\FL \wedge_N \FL$ as the $2TN$-spectrum
defined by
\begin{align*}
  (\FL \wedge_N \FL)_r = \FL_r \wedge \FL_r,
\end{align*}
with structure maps
\begin{align*}
  \sigma^\wedge_r \colon (\FL_r \wedge_N \FL_r )^{2TN} \cong
  \FL^{TN}_r \wedge_N \FL^{TN} \xrightarrow{\sigma_r \wedge_N
    \sigma_r} \FL_{r+1} \wedge_N \FL_{r+1}.
\end{align*}
We can easily consider any other instance $\tFL$ as a $2TN$ spectrum
by simply forgetting every other ex-space in the sequence
$\tFL_r$, and making the new structure maps the old
$\sigma_{2r,2r+2}$. This does not change its homology, and in fact
defines an equivalence of categories, but we wont need the latter. We
will thus define $\mu$ in the category $\Sp_N^{2TN}$.

We start by defining the fiber-wise product on the level of the
underlying smooth manifolds as maps
\begin{align*}
  G_r \colon (T^*\Lambda_r N) \times_N (T^*\Lambda_r N)
  \to T^*\Lambda_{2r} N
\end{align*}
given by simply concatenating. I.e. given two points
$\arzet=(q^1_0,p^1_0,\dots q^1_{r-1},p^1_{r-1})$ and
$\arzto=(q^2_0,p^2_0,\dots q^2_{r-1},p^2_{r-1})$ with $q^1_0=q_0^2$ (the
definition of $\times_N$ provides precisely this equality) we can
define
\begin{align} \label{eq:Grdef}
  G_r(\arzet,\arzto)=
  (q^1_0,p^1_0,\dots,q^1_{r-1},p^1_{r-1},q^2_0,p^2_0,\dots,q^2_{r-1},p^2_{r-1}).
\end{align}
This is well-defined because
$\dist(q^1_{r-1},q^2_0)=\dist(q^1_{r-1},q^1_0)<\epsin$ and similar for
$q^1_0$ and $q^2_{r-1}$.

Now assume we are given a Hamiltonian $H$ and a sub-division
$\alpha\in \Delta^{r-1}$ such that $S_r$ is defined for the pair
$(H,\alpha)$. We may then define $\tS_{2r}$ as the finite dimensional
approximation using the Hamiltonian $2H$ and the subdivision
$2_r\alpha\in \Delta^{2r-1}$ given by
\begin{align}\label{eq:23}
  (2_r\alpha)_j = (\alpha_{(j\phantom{i}\textrm{mod}\phantom{i} r)}) /2.
\end{align}
The twiddle is put on $\tS_{2r}$ because we wish to emphasize that we
are not using the same Hamiltonian and sub-division. In fact twiddles
will signify that we are dealing with structure related to the target
$\tFL$ (which will be constructed during the course of this
section). Using that $q^1_0=q^2_0$ in the definition of $G_r$ it is
easy to check that these functions satisfy
\begin{align} \label{eq:7}
  \tS_{2r}(G_r(\arzet,\arzto)) = S_r(\arzet)+S_r(\arzto).
\end{align}
The factor 2 in the Hamiltonian is compensated by the $\tfrac{1}{2}$
in the sub-division, and is explained by the fact that the action is
additive under strict concatenation. Here strict concatenation takes
two curves parametrized by an interval of length 1 and spits out a
curve parametrized by an interval of length 2. However, since this is
not the concatenation we use (we identify the result as parametrized
by an an interval of length 1 using the unique affine
reparametrization) we effectively need to scale up the Hamiltonian and
scale down the length of the pieces in the sub-division to get the
same result.

The point of these maps and Equation~\eqref{eq:7} is that $G_r$ now
induces a map of Parametrized Conley indices. Indeed, for any pair
$(C_\tau,\uC_\tau)$ in $T^*\Lambda_r N$ as in
Section~\ref{cha:paracon}, Equation~\eqref{eq:7} shows that
\begin{align*}
  G_r(C_\tau,C_\tau) &\subset (\tS_{2r})^{-1}([2a,2b]) \qquad
  \textrm{and} \\
  G_r(\uC_\tau,C_\tau) \cup
  G_r(C_\tau,\uC_\tau) &\subset (\tS_{2r})^{-1}([2a,a+b]),
\end{align*}
and thus as usual by composing with the flow of $-\tY_{2r}$ (to get the
image inside some $\tC_t$ defining the target parametrized index) we
get an induced ex-map
\begin{align*}
  \mu_r \colon I_a^b(S_r,Y_r)_N \wedge_N I_a^b(S_r,Y_r)_N
  \to I_{a+b}^{2b}(\tS_{2r},\tY_{2r})_N.
\end{align*}
Note that this, indeed, is an ex-map over $N$ because $G_r$ preserves
the projections to $N$ and so does the flow. We would like to say
that this defines a map of generating functions spectra. However,
strictly speaking it does not. Indeed, the doubling of the
sub-division in Equation~\eqref{eq:23} is not compatible with the top
face inclusions used to define the suspensions in the generating
function spectra. I.e. the diagram
\begin{align} \label{eq:8sdfsdfa}
  \xymatrix{
    \Delta^{r-1} \ar[r]^{2_r(\sm)} \ar[d]^{d_r} & \Delta^{2r-1}
    \ar[d]^{d_{2r+1} \circ d_{2r}} \\
    \Delta^r \ar[r]^{2_{r+1}(\sm)} & \Delta^{2r+1}
  }
\end{align}
does not commute. The difference is precisely a reordering of the
numbers in the sequence (moving one of the introduced zeros
around).

To be as precise as possible we now assume that we have
\begin{itemize}
\item{a smooth family of Hamiltonians $H^s \in \HS$, $s\in
    \Jinf$, defined as in Section~\ref{theex} but with $H_\infty \in
    \HS_\infty$ such that $2H_\infty \in \HS_\infty$,}
\item{a smooth family $\alpha^s\in \Delta^\infty$,$s\in \Jinf$ of
    sub-divisions,}
\item{constants $a < -2\norm{F}-1$ and $b = \norm{F} + 1$ ($\norm{F}$
    as in Section~\ref{theex}), and}
\item{a strictly increasing sequence $s_r \in \Jinf$, $r\geq r_0$}
\end{itemize}
defining $\FL$ (the source) as in Section~\ref{theex}. Here we have
put slightly more restrictions on $H_\infty$ and the bounds $a$ and
$b$ than we did in that section, and we will need this to define the
product.

We will need the target spectra $\tFL$ to be defined using a smooth
family of sub-divisions $\tA^s \in \Delta^\infty$, $s\in \Jinf$ such
that $\tA^{s_r}=2_r(\alpha^{s_r})$. Since 
\begin{align*}
  l(2_r(\alpha^s)) = \tfrac{1}{2}l(\alpha^s)  
\end{align*}
there are no problems in choosing such a family such that
\begin{align*}
  l(\tA^s) \leq \tfrac{1}{2}(\alpha^s).
\end{align*}
Basically this family could be constructed as a twisted $\alpha^s$,
where ``twisting'' takes place on each open interval $]s_r,s_{r+1}[$
and compensates for the non-commutativity of the diagram in
Equation~\eqref{eq:8sdfsdfa}. We thus use the following data to define
$\tFL$:
\begin{itemize}
\item{The smooth family of Hamiltonians $2H^s$.}
\item{A strictly increasing sequence $\ts_r$ with $\ts_{2r}=s_r$.}
\item{A smooth family $\tA^s$ such that $\tA^{s_{2r}}=2_r\alpha^{s_r}$
  satisfying the above.}
\item{the constants $\ta=a+b < -\norm{F}$ and $\tb=2b$.}
\end{itemize}
By the above we have made sure that $\tS_{2r}^s$ is defined for $s \in
[0,\ts_{2r}]$. It is no problem extending the sequence $\ts_{2r}$ to
include odd numberings $\ts_{2r+1}$. However, we are throwing the odd
levels away anyway so in light of Lemma~\ref{lem:exhomjump} it really
does not matter.

\begin{Definition}
  \label{def:7}
  The ex-maps
  \begin{align*}
    \mu_r \colon \FL_r \wedge_N \FL_r \to \tFL_{2r}
  \end{align*}
  is defined when $s_r>5$ by the map induced by $G_r$ above.
\end{Definition}

\begin{Lemma}
  The ex-maps $\mu_r$ fit together to define a morphism
  $\mu=(\mu_r,h_r,r)$ in $\Sp_N^{2TN}$ as in Equation~\eqref{eq:8sder}.
\end{Lemma}

See Section~\ref{cha:fib-fl} for definition of morphism in
$\Sp_N^{2TN}$.

\begin{bevis}
  Since the proof is highly technical and not very deep (except the
  usual issues about reordering suspension for spectra), we only
  provide a sketch of the argument. Also we do not \emph{really} need 
  this because as an alternative one can prove the main theorem using
  a compactness argument and a single $\mu_r$, provided $r$ is chosen
  large enough. However, this way of thinking makes notation much more
  compact.

  Firstly it is convenient to apply the functor from
  Appendix~\ref{cha:funcrep} to trivialize the bundles, and to do so
  in a way that makes the trivial bundles even dimensional (this
  avoids some sign issues).

  The non-commutativity in Equation~\eqref{eq:8sdfsdfa} is taken
  care of by the fact that the family $\tA^s$ interpolates between the 
  different $2_r(\alpha^{s_r})$ and
  $2_{r+1}(\alpha^{s_{r+1}})$. However, inspecting the
  Thom-suspensions from Section~\ref{cha:fib-susp} we see that when
  comparing the two sides we are in fact putting in the vector from
  the first copy of $TN$ over different $q_j$'s, and in fact shifting
  the identification of the $p_j$'s coming from the second factor by
  one. Indeed, this corresponds to the usual problem with
  non-associativity of smash-products of spectra and the rearranging
  needed to define a map from the product of two spectra to a third
  spectrum.
\end{bevis}


\chapter{The Global Product}\label{cha:glob-prod}

The products $\mu_r$ in Section~\ref{cha:fib-w-prod} was defined
fiber-wise (as ex-maps) over $N$. In this section we define associated
products on the ``total spaces'' by combining the fiber-wise products
with the intersection product on $N$, and we prove that it will be
compatible with the intersection product on $L$. We will use the
description of the intersection product as the map induced by a
Pontryagin-Thom collapse map $N\times N \to N^{TN}$. The goal is thus
to extend $\mu_r$ from Section~\ref{cha:fib-w-prod} to induce a map
\begin{align*}
  (\mu_r)_! \colon \FL_r/N \wedge \FL_r/N \to (\tFL_{2r}^{TN})/N,
\end{align*}
and then prove that for large $r$ this extends the intersection
product on constant loops in $L$ (compare with
Corollary~\ref{cor:2}). Note that we will use the fact that $\FL$
is a $TN$-spectrum below and use a suspension from
Section~\ref{cha:fib-susp} to get the $TN$ factor on the right
above. However, in general it is not important that such a fiber-wise
product is defined on a $TN$ spectrum to induce a global product as
above. We use this not so general construction because it is the 
most convenient way of proving that the product extends the
intersection product on $L$.

The additivity rule in Equation~\eqref{eq:7} of $G_r$ does not a
priori extend to a neighborhood of
\begin{align*}
  T^*\Lambda_r N \times_N T^*\Lambda_r N \subset T^*\Lambda_r
  N \times T^*\Lambda_r N.
\end{align*}
However, if we insert an extra point in the target $(q,p)\in T^*N$ we
can actually accomplish this. So extend $G_r$ to points
\begin{align*}
  \arzet &= (q_0^1,p_0^1,\dots,q_{r-1}^1,p_{r-1}^1) \qquad\textrm{and}\\
  \arzto &= (q_0^2,p_0^2,\dots, q_{r-1}^2,p_{r-1}^2)
\end{align*}
such that all three distances $\dist(q_0^1,q_{r-1}^2)$,
$\dist(q_0^2,q_{r-1}^1)$ and $\dist(q_0^1,a_0^2)$ are less than
$\epsin$. Then
\begin{align*}
  G_r(\arzet,\arzto) =
  (q_0^1,p_0^1,\dots,q_{r-1}^1,p_{r-1}^1,q_0^2,p_0^2,\dots,q_{r-1}^2,p_{r-1}^2)
\end{align*}
lies in $T^*\Lambda_{2r} N$. In the following we are implicitly using
Hamiltonians $H$ and $2H$ and sub-divisions $\alpha$ and
$\tA=2_r(\alpha)$ as assumed in Equation~\eqref{eq:7}. We wish to
construct a smooth function
\begin{align*}
  p=F(\arzet,\arzto) \in T^*_{q_0^2} N
\end{align*}
such that if we define
\begin{align}\label{eq:8sd}
  G_r^+(\arzet,\arzto)= 
  (q_0^1,p_0^1,\dots,q_{r-1}^1,p_{r-1}^1,q_0^2,p_0^2,\dots,q_{r-1}^2,p_{r-1}^2,q_0^2,p),
\end{align}
which by assumption on $\arzet,\arzto$ above lies in
$T^*\Lambda_{2r+1} N$, we get the sub-additivity rule
\begin{align} \label{eq:8lkl}
  S_r(\arzet) + S_r(\arzto) \leq \tS_{2r+1}(G_r^+(\arzet,\arzto)),
\end{align}
with equality only on the diagonal where
$q_0^2=q_0^1$. We will construct it such that the equality on the
diagonal is non-degenerate in the sense that moving all terms to one
side of the equation the diagonal is a non-degenerate critical
manifold in the sense of Morse-Bott. Here we used the usual top face
inclusion $\tA \in \Delta^{2r} \subset \Delta^{2r+1}$ to define $\tS_{2r+1}$
as we did when defining the suspension maps in
Section~\ref{cha:fib-susp}. Indeed, as mentioned in the introduction
we will use this suspension to get the extra $TN$ factor.

Using Definition~\ref{def:Ardef} we calculate (keeping in mind the
reindexing from Equation~\eqref{eq:8sd}) the difference for arbitrary
$p \in T^*_{q_0^2}N$ to be
\begin{align}\label{eq:8sdgadfha}
  \tS_{2r+1}(&\arz) - S_r(\arzet) + S_r(\arzto) = \\ & =
  (p_{r-1}^1)^-(\exp_{(q_{r-1}^1)^-}^{-1}(q_0^2)-\exp_{(q_{r-1}^1)^-}^{-1}(q_0^1))
  + p \exp_{q^2_0}^{-1}(q_0^1).
\end{align}
Indeed, since $\tA_{2r}=0$ we do not need the minuses on $q_0^2$ and
$p$ in the last term. If $q_0^2=q_0^1$ (the diagonal on which we
worked in Section~\ref{cha:fib-w-prod}) both terms vanishes and we can
put $p$ equal to anything. However, to solve when $q_0^1\neq q_0^2$ we
need that the gradient of this thing on the diagonal with respect to
both $q_0^2$ and $q_0^1$ vanishes. This is easily accomplished by
noticing that there is a unique $p$ where this is true. Indeed, the
gradient on the diagonal is anti-symmetric in $q_0^2$ and $q_0^1$, and
since the term $p\exp_{q_0^2}^{-1}(q_0^1)$ has gradient with respect
to $(q_0^1,q_0^2)$ equal to $(p,-p)$ there is such a unique $p$. Then
it is simply a matter of defining the general $p \in T^*_{q_0^2}N$ by
choosing some extension of this unique $p$ defined on the diagonal and
adding a positive smooth function times $-\exp_{q_0^2}^{-1}(q_0^1)$ to
make the quantity in Equation~\eqref{eq:8sdgadfha} non-degenerately
negative for all $\arzet$ and $\arzto$. In fact by picking this
function we can get the value of $\tS_{2r+1}$ as small as needed on
any compact set disjoint from the diagonal, which we will use in the
proof below.

For any $\beta>0$ and any two maps $p_i \colon A_i \to N$ for
$i=1,2$ we define 
\begin{align*}
  A_1 \times^{\beta}_N A_2 = \{(x,y)\in A_1\times A_2 \mid
  \dist(p_1(x),p_2(y)) \leq \beta\}.
\end{align*}
Let $C_\tau \subset T^*\Lambda_r N$ be as in Section~\ref{cha:paracon}
defining the parametrized Conley index $\FL_r$.

\begin{Lemma}
  \label{lem:prodext}
  There is a $\beta>0$ such that the map
  \begin{align*}
    G_r^+ \colon C_\tau \times_N^\beta C_\tau \to
    T^*\Lambda_{2r+1}N 
  \end{align*}
  is well-defined and being careful about the choices above it will
  induce a map
  \begin{align*}
    (\mu_r)_! \colon \FL_r /N \wedge \FL_r/N \to \tFL_{2r}^{TN}/N,
  \end{align*}
  which realizes the Pontryagin-Thom collapse map to the new factor
  $TN$.
\end{Lemma}

Note that we can not write $\tFL_{2r}^{TN}/N$ as $\tFL_{2r+1}$ since
the latter also implies a change of Hamiltonian and sub-division by
increasing $s$ from $\ts_{2r}$ to $\ts_{2r+1}$.

\begin{bevis}
  Since $C_\tau$ is compact and satisfy the inequalities assumed on
  $\arzet$ and $\arzto$ above $G_r^+$ must be well-defined on a small
  neighborhood of $C_\tau \times_N C_\tau$. This small neighborhood
  will contain $C_\tau \times_N^\beta C_\tau$ for small $\beta$. By
  using the sub-additivity (and possibly composing with the negative
  gradient flow) we get a map
  \begin{align*}
    C_\tau \times_\beta C_\tau \to \tFL_{2r+1}^{TN}/N
  \end{align*}
  induced by $G_r^+$, which sends $\uC_\tau \times_N^\beta C_\tau \cup
  C_\tau \times_N^\beta \uC_\tau$ to the base point. Inspecting
  $C_\tau \times_N^\beta C_\tau$ we see that it is neither open nor
  close. It will be a closed set if we add the compact set of
  $(\arzet,\arzto) \in C_\tau\times C_\tau$ which has
  $\dist(q_0^1,q_0^2) = \beta$. This compact set is isolated from
  the diagonal and by being careful with the choices above we can
  make sure that $G_r^+$ sends this compact set to points with value
  less than $\ta$. I.e. the lower bound used to define
  $\tFL_{2r}$. This means that the induced map extends to all of
  $C_\tau \times C_\tau$ by sending everything outside $C_\tau
  \times_N^\beta C_\tau$ to the base point. This is the usual idea of
  Pontryagin-Thom collapse maps - except we now have god control over
  what $G_r^+$ does to the action values.

  For the last statement notice that when defining the suspension in
  Section~\ref{cha:fib-susp} the vectors of $TN$ corresponds to the
  $p$ defined above and with a large factor in front of
  $\exp_{q_0^2}^{-1}(q_0^1)$ we see that this $p$ really measures the
  distance of the two points $q_0^1$ and $q_0^2$, which is how one
  constructs the Pontryagin-Thom collapse map.
\end{bevis}

The rest of this section is devoted to proving that for large $r$ this
product extends the intersection product on constant loops in
$L$. This is done by combining this construction with the construction
in Section~\ref{cha:inc-const}. We will not need the subtle two last
lemmas in that section. Indeed, the statement we need here is a purely
homological statement - without any fiber-wise concerns.

\begin{Lemma}
  \label{lem:4de}
  With $\F$ as in Corollary~\ref{cor:2} the diagram
  \begin{align*}
    \xymatrix{
      H_*(L;\F) \otimes H_*(L;\F)
      \ar[r]^{\cap} \ar[d]^{(i_r)_* \otimes (i_r)_*} &
      H_*(L;\F) \ar[d]^{(\ti_{2r})_*} \\
      H_*(\FL_r,N;\F) \otimes H_*(\FL_r,N;\F) \ar[r]^-{(\mu_r)_!} &
      H_*(\tFL^{TN}_{2r+1},N;\F)
    }
  \end{align*}
  commutes up to a possible graded sign. Here $\cap$ is the
  intersection product and the vertical maps are the inclusions of
  constant loops described in Section~\ref{cha:inc-const}.
\end{Lemma}

\begin{bevis}
  We assume that $r$ is as large as needed in the construction in
  Section~\ref{cha:inc-const}, and the notation from that section is
  used in the following. We need to identify the product on the
  inclusion constructed there, and we start by taking a closer look at
  the homotopy $(S_r^s,X_r^s),s\in[\epsilon,s_r]$.

  In Section~\ref{cha:inc-const} we considered a varying lower bound
  $a\colon [\epsilon,s_r] \to \R$, but we got rid of this
  $s$-dependence by translating the functions based on the parameter
  $s$ (which we shall generalize in Lemma~\ref{lem:6}). We will also
  need to vary the upper bound $b$ in this proof. Indeed, we will need
  the varying bounds:
  \begin{align*}
    \ua(s) & = -s +\epsilon/2 \\
    \ub(s) & = s/3
  \end{align*}
  for $s\in [\epsilon,s_r]$. This is the same $\ua(s)$ as in
  Section~\ref{cha:inc-const}, which is regular for $S_r^s$. The
  construction of the product is such that we need to define the
  bounds associated to the target $\tFL$ as:
  \begin{align*}
    \uta(s) &= \ua(s)+\ub(s) = - 2s/3 +\epsilon/2 \\
    \utb(s) &= 2\ub(s) = 2s/3.
  \end{align*}
  Varying bounds is not really problematic when dealing with
  parametrized Conley indices. However, we postpone the technical
  details of this to the following lemma, and simply suppress the $s$
  from the notation.
  
  We may construct the products $(\mu_r)_!^s$ above for each $s\in
  [\epsilon,s_r]$. In fact using a compactness argument we may define
  the family of these as an ex-map:
  \begin{align*}
    (\mu_r)_!^\sm \colon I_{\ua}^{\ub}(S_r^\sm,X_r^\sm)_{[\epsilon,s_r]}
    \wedge_{[\epsilon,s_r]}
    I_{\ua}^{\ub}(S_r^\sm,X_r^\sm)_{[\epsilon,s_r]} \to 
    I_{\uta}^{\utb}(\tS_{2r+1}^{\sm},\tX_{2r+1}^{\sm})_{[\epsilon,s_r]}.
  \end{align*}
  Lemma~\ref{lem:6} tells us that there are parallel transports
  \begin{align} \label{eq:4parpro}
    p_s^{s'} &\colon I_a^b(S_r^s,X_r^s) \to I_a^b(S_r^{s'},X_r^{s'}),
    \qquad \textrm{and} \\
    \tp_s^{s'} &\colon I_{\ta}^{\tb}(\tS_{2r+1}^s,\tX_{2r+1}^s) \to
    I_{\ta}^{\tb}(\tS_{2r+1}^{s'},\tX_{2r+1}^{s'}), \notag
  \end{align}
  both for $s\leq s' \in [\epsilon,s_r]$.

  This provides a homotopy
  \begin{align*}
    \tp_t^{s_r} \circ (\mu_r)_!^t \circ (p_\epsilon^t
    \wedge p_\epsilon^t) \colon
    I_a^b(S_r^\epsilon,X_r^\epsilon) \wedge
    I_a^b(S_r^\epsilon,X_r^\epsilon) \to 
    I_{\ta}^{\tb}(\tS_{2r}^{s_r},\tX_{2r}^{s_r})    
  \end{align*}
  for $t\in [\epsilon,s_r]$, which shows that the diagram
  \begin{align*}
    \xymatrix{
      I_\ua^\ub(S_r^\epsilon,X_r^\epsilon) \wedge
      I_\ua^\ub(S_r^\epsilon,X_r^\epsilon) \ar[rr]^-{(\mu_r)_!^\epsilon}
      \ar[d]^{p_\epsilon^{s_r} \wedge_N p_\epsilon^{s_r}} &&
      I_{\uta}^{\utb}(\tS_{2r+1}^\epsilon,\tX_{2r+1}^\epsilon)
      \ar[d]^{\tp_\epsilon^{s_r}} \\
      I_\ua^\ub(S_r^{s_r},X_r^{s_r}) \wedge
      I_\ua^\ub(S_r^{s_r},X_r^{s_r}) \ar[rr]^-{(\mu_r)_!^{s_r}} &&
      I_{\uta}^{\utb}(\tS_{2r+1}^{s_r},\tX_{2r+1}^{s_r})          
    }
  \end{align*}
  commutes. The $\ua(s_r),\ub(s_r),\uta(s_r),\utb(s_r)$ we have used
  as bounds in the lower horizontal map ($s_r$ were suppressed) are as
  in Section~\ref{cha:inc-const} all smaller than the actual chosen
  $a,b,\ta,\tb$ when we defined the product $(\mu_r)_!$ in
  Lemma~\ref{lem:prodext}. Moreover, the construction of
  $(\mu_r)_!$ is obviously compatible with quotients to larger such
  bounds if of course these satisfy $\ta=a+b$ and $\tb=2b$ before and
  after. So we get another commutative diagram
  \begin{align*}
    \xymatrix{
      I_\ua^\ub(S_r^\epsilon,X_r^\epsilon) \wedge
      I_\ua^\ub(S_r^\epsilon,X_r^\epsilon) \ar[rr]^-{(\mu_r)_!^{s_r}}
      \ar[d] &&
      I_{\uta}^{\utb}(\tS_{2r+1}^\epsilon,\tX_{2r+1}^\epsilon)
      \ar[d] \\
      I_a^b(S_r^{s_r},X_r^{s_r}) \wedge
      I_a^b(S_r^{s_r},X_r^{s_r}) \ar[rr]^-{(\mu_r)_!} &&
      I_{\ta}^{\tb}(\tS_{2r+1}^{s_r},\tX_{2r+1}^{s_r})          
    },
  \end{align*}
  which fits directly below the previous one. Putting them on top of
  each other we see that the inclusion of constant curves corresponds
  to going from top to bottom, and so we only need to see that
  $(\mu_r)_!^\epsilon$ realizes the intersection product. So we need
  to consider
  \begin{align*}
    G_r^+ \colon U \subset T^*\Lambda_rN \times T^*\Lambda_r N \to
    T^*\Lambda_{2r+1} N
  \end{align*}
  a little more careful in the case of
  $(H^\epsilon,\alpha^\epsilon)$.
  
  Since $\epsilon$ is small we have from Lemma~\ref{lem:1a} and the
  ensuing construction that both Conley indices used in the definition
  $(\mu_r^\epsilon)_!$ are Thom-spaces over $L$. In fact the
  non-degeneracy proved for small $s$ in Lemma~\ref{lem:1a} close to
  $L$ can easily be extended to $S_r^\epsilon$. Indeed, consider
  $\R^{2n}$ where the action on closed loops for $H=0$ is well
  understood and the constant loops (a copy of $\R^{2n}$) is a
  non-degenerate critical manifold - if changing $H$ by a
  sufficiently small perturbation which has a non-degenerate
  critical manifold then the action will have this same non-degenerate
  critical manifold. So let $E_-$ be a maximal sub-space of
  $T(T^*\Lambda_r N)_{\mid L}$ on which the Hessian of $S_r^\epsilon$
  is negative definite. Also let $\nu_L$ be the normal bundle of the
  diagonal $L \subset L\times L$. We have the commutative diagram
  \begin{align*}
    \xymatrix{
      U \ar[r]^-{G_r^+} & T^*\Lambda_{2r+1} N \\
      L \times L \ar[u] & L \ar[l]_\Delta \ar[u].
    }
  \end{align*}
  Now $E_- \times E_-$ is a sub-bundle in the tangent bundle of $U$
  along $L\times L$. Restriction it to the diagonal we get $E_- \oplus
  E_- \subset TU_{\mid L}$. The bundle $\nu_L$ may also be considered
  such a sub-bundle by $\nu_L \subset T(L\times L)_{\mid L} \subset
  TU_{\mid L}$, and these two sub-bundles intersects trivially. Indeed,
  the Hessian of $S_r^\epsilon(\arzet) + S_r^\epsilon(\arzto)$ is zero
  on $L\times L$ so $\nu_L$ is part of the kernel and cannot be part
  of any negative bundle.

  Now $G_r^+$ is an embedding and by the sub-additivity and its
  non-degeneracy on the diagonal we see that $E_- \oplus E_- \oplus
  \nu_L$ is send to a negative sub-space, which by dimension counting
  is maximal. This means that the map realizes the Printmaking-Thom
  collapse map from $L\times L$ to $L^{\nu_L}$, and the extra negative
  bundle of $E_-\oplus E_-$ is by assumption oriented - so on the
  level of homology the lemma follows.
\end{bevis}

\begin{Lemma}
  \label{lem:6}
  The parallel transport maps in Equation~\eqref{eq:4parpro} exists.
\end{Lemma}

\begin{bevis}
  The varying bounds are easy to get rid of by simply replacing
  $S_r^s$ with the function
  \begin{align*}
    \frac{S_r^s-\ua(s)}{\ub(s)-\ua(s)}.
  \end{align*}
  This does not change the fact that $X^s_r$ is a pseudo-gradient, and
  now the bounds $\ua(s)$ and $\ub(s)$ are replaced by the constants
  $0$ an $1$. We can do the same for the target approximations
  $\tS_{2r+1}^s$ using $\utb(s)$ and $\uta(s)$.

  The trouble with the right hand side is that $\uta$ probably has to
  cross several critical values because it is far from
  $\ts_{2r}=s_r$. However, this is were the parallel transport in
  Lemma~\ref{lem:1} comes in handy again.

  Indeed, for there to be a parallel transport forward in $s$ on the
  right hand side we need that
  \begin{align*}
    (\pd{s})^\arz \pare*{\frac{\tS_{2r+1}^\sm-\uta}{\utb-\uta}} <0
  \end{align*}
  for any $\arz \in T^*\Lambda_r N$ regular for $\tS_{2r+1}^s$ and such
  that $\tS_{2r+1}^s(\arz)=\uta$. This translates into the need for
  \begin{align} \label{eq:dsda}
    (\pd{s})^\arz \tS_{2r+1}^{\sm} < \pd{s} \uta
  \end{align}
  reflecting that $\uta$ is not constant. By Lemma~\ref{lem:closact}
  the left hand side is $-1$ and by the choice of $\uta$ the right
  hand side is $-2/3$. The intuition is that even though the lower
  bound $\uta$ goes down with speed $-2/3$, the critical values of
  $\tS_r^\sm$ moves down faster - i.e. with speed $-1$. So, indeed,
  there is a parallel transport
  \begin{align*}
    \tp_s^{s'} \colon I_{\ta}^{\tb}(\tS_{2r+1}^s,\tX_{2r+1}^s) \to
    I_{\ta}^{\tb}(\tS_{2r+1}^{s'},\tX_{2r+1}^{s'}).
  \end{align*}
  The source is much easier to handle since $\ua$ is always
  regular. In fact it is a limit case of the above idea where the
  bound $\ua$ moves with same speed $-1$ as $s$, but then we still
  have a parallel transport because it is always regular (it moves
  parallel to the critical values without being one).
\end{bevis}


\chapter{The Product on the Serre Spectral
  Sequence}\label{cha:serprod} 

In this section we describe how the products define a product on the
Serre spectral sequence. First we consider an external smash product
and the associated map on spectral sequences, which will be an
appropriate derivation on all pages inducing the map on subsequent
pages. Then we consider the Pontryagin Thom collapse map, and finally
we use the map $\mu$ defined in Section~\ref{cha:fib-w-prod}.

As in the proof of Proposition~\ref{Spectral} this will work with any
coefficient ring $\F$, but for notational purposes we suppress this
and consider only $\F=\Z$. We have borrowed some ideas from
\cite{MR2039760}.

First we define the exterior smash-product $\wedge^e$ as a functor
from $\Sp_W^\beta \times \Sp_{W'}^{\beta'}$ to $\Sp_{W\times
  W'}^{\beta\times \beta'}$ by
\begin{align*}
  (\As \wedge^e \Bs)_r = (\As_r \times \Bs_r) :_{W\times W'} (\As_r
  \times W' \cup W \times \Bs_r).
\end{align*}
So that the fiber at $(x,x') \in W\times W'$ is $\As_{r\mid x}
\wedge \Bs_{r\mid x'}$. The structure maps are the obvious ones
defined using the identification
\begin{align*}
  (\As \wedge^e \Bs)_r^{\beta\times \beta'} \cong
  (\As_r^\beta \times \Bs_r^{\beta'}) :_{W\times W'}
  (\As^\beta_r \times W' \cup W \times \Bs^{\beta'}_r).
\end{align*}
The Eilenberg-Zilber operators defined in Section~\ref{serre} induce
maps on the chains relative to the sections and we thus have maps
\begin{align*}
  C_*(\As_r , W) \otimes C_*(\Bs_r,W') \to C_*((\As \wedge^e
  \Bs)_r,W\times W').
\end{align*}
When both $\beta$ and $\beta'$ are trivialized of respective dimension
$l$ and $l'$ (which we may assume) and one of them are even
dimensional the suspensions $\Sigma_*$ defined on the chain level in
Equation~\eqref{eq:1023} is compatible with these maps. I.e. the
diagram
\begin{align*}
  \xymatrix{
    C_*(\As_r , W) \otimes C_*(\Bs_r,W')
    \ar[r] \ar[d]^{\Sigma_* \otimes \Sigma'_*} &
    C_*((\As \wedge^e \Bs)_r,W\times W') \ar[d]^{\Sigma''_*} \\
    C_{*+l}(\As_{r+1} , W) \otimes C_{*+l'}(\Bs_{r+1},W') \ar[r] &
    C_{*+l+l'}((\As \wedge^e \Bs)_{r+1},W\times W')
  }
\end{align*}
commutes. Here $\Sigma_*''$ is the suspension on chains induced by the
structure maps of $\As \wedge^e \Bs$. Indeed, this follows from the
commutativity mentioned in Remark~\ref{rem:6} and the way the
suspension maps on chains were defined. Furthermore, this is graded
from the bi-graded tensor complex to the graded complex. In addition
it adds the filtration degrees. This implies that it induces a map
from the tensored spectral sequence which on each page is a derivation
inducing the product on the next page (see e.g. \cite{MR1793722}).

The next step is the spectral sequence version of the Pontryagin-Thom
collapse map. For this we need a relative version of the spectral
sequence. So assume $B'\subset B\subset M$ are closed Neighborhood
retracts of any compact manifold $M$ and $\As$ is a 3S fibrant
$l$-spectrum over $M$. Then we have a relative Serre spectral sequence
given by the relative chain complex
\begin{align*}
  C_*(\As_{\mid B},\As_{\mid B'}) = \colim_{r\to\infty}
  C_{*+lr}(\As_{r\mid B} , \As_{r\mid B'} \cup B)
\end{align*}
which on page two can be identified with
$H_*(B,B';H_*(A_{\mid \bullet}))$. Indeed, this is a
relative version of Proposition~\ref{Spectral}.

For any closed proper sub-manifold $M'\subset M$ with oriented (in
general with respect to $\F$) normal bundle $\nu$ of dimension $k$ we
would like to create a spectral sequence version of the
``Pontryagin-Thom collapse'' map
\begin{align} \label{ptcwish}
  H_*(\As) \to H_{*-k} (\As_{\mid M'}),
\end{align}
where $\As$ is any 3S fibrant $l$-spectrum over $M$. We will start by
replacing the spectral sequence with a spectral sequence which is
isomorphic from page 2 and onwards.

We may identify $D\nu$ with a closed tubular neighborhood of $M$. Let
$U\supset D\nu$ be a neighborhood whose closure deformation retracts
onto $D\nu$, and let $D\nu^c$ denote the open complement of $D\nu$. We
may use Proposition 2.21 in Hatcher~\cite{MR1867354} together with the
identification of page 2 in our proof of Proposition~\ref{Spectral},
and conclude that the inclusion
\begin{align} \label{altspec}
  C_*^{U,(D\nu)^c}(\As)
  :=\colim_{r\to \infty} C^{\As_{r\mid
      U},\As_{r\mid(D\nu)^c}}_{*+lr}(\As_r,N) \to C_*(\As),
\end{align}
induces an isomorphism on page 2 and thus also on higher pages (the
sub-complex is given the induced filtration). So we may replace
the spectral sequence for $\As$ with this.

Now pick a deformation retraction of the closure of $U$ onto $D\nu$ we
can as usual use Lemma~\ref{lem:4} to make sure that it is
path-smooth. We may also make sure that everything that start outside
$D\nu$ never enters the interior of $D\nu$, and thus everything
outside $D\nu$ must be deformed onto the sphere bundle $S\nu \subset
D\nu$. Using this with Definition~\ref{3S} and restricting to $1\in I$
we get a map
\begin{align*}
  p \colon \As_{\mid \overline{U}} \to \As_{\mid D\nu},
\end{align*}
where all the fibers outside $D\nu$ are mapped to fibers over points
in $S\nu$. We may use this on the chain complex in~\eqref{altspec}
by sending all simplices $\alpha$ living over $U$ to $p_*(\alpha)$
and the rest to $0$. This is a chain map into $C_*(\As_{\mid
  D\nu},\As_{\mid S\nu})$ because the new boundary we create by
``ignoring'' some simplices in $D\nu^c$ is all sent to $C_*(Y_{\mid
  S\nu})$. It is also filtration preserving because the projections
of simplices is governed by the deformation retraction.

To continue this map into $C_{*-k}(\As_{\mid M})$ we simply do as in
\cite{MR2039760} and pick a Thom class $\tau$ in $C^*(D\nu,S\nu)$
which vanishes on simplices contained in $S\nu$ and on
degenerate simplices. Denoting $\tau$ pulled back to
$\As_{r\mid D\nu}$ by $\tau'$ we get 
\begin{align*}
  \tau' \cap (-) \colon C_{*+lr}(\As_{r\mid D\nu},\As_{r\mid S\nu}) \to
  C_{*-k+lr}(\As_{r\mid D\nu}),
\end{align*}
which is a chain map for the usual reasons and because $\tau'$ is zero
on simplices in $\As_{r\mid S\nu}$. This lowers filtration by $k$ for
the same reasons as in \cite{MR2039760} although we above made the
choice to cap away the first $c$ indices. This choice makes it commute
with our suspensions (which adds to the other end of the indices
defining the simplices). So this is compatible with the limit as $r$
tends to infinity, and thus induces a map of spectral sequences which
lowers the degree by $(k,0)$. To end up in $C_*(\As_{\mid M'})$ we
again simply use a path-smooth deformation retraction of $D\nu$ onto
$M'$ with Definition~\ref{3S} and restrict to $1\in I$ to get a
morphism
\begin{align*}
  \As_{\mid D\nu} \to \As_{\mid M},
\end{align*}
which preserves filtration.

Using the above steps on $\FL \wedge^e \FL$ and the diagonal $N\subset
N\times N$ we get a map from the tensor product of the spectral
sequence defined using $\FL$ with itself to the spectral sequence
defined by $\FL \wedge_N \FL$. Indeed, $\FL \wedge_N \FL$ is the
restriction of $\FL \wedge^e \FL$ to the diagonal. Composing this with
the product $\mu$ from Section~\ref{cha:fib-w-prod} we get a product
on the spectral sequence associated to $\FL$. It lowers degree by
$(d,0)$ due to the Thom isomorphism used at one of the steps - so if
unital the unit would have degree $(d,0)$ (the fundamental class for
$L$ is in fact a unit, but we do not prove this at any point). Since
Section~\ref{cha:glob-prod} and the construction here both realizes
the Pontryagin-Thom construction combined with the fiber-wise product
we have the following lemma.

\begin{Lemma}
  \label{lem:14asd}
  The product above induces on non-filtered homology the same map as
  the product in Section~\ref{cha:glob-prod} does up the the shift
  by $d$ coming from the Thom-isomorphism.
\end{Lemma}


\chapter{Proof of the Main Theorem} \label{cha:proof-theor-refh}

As usual $N$ is a closed $d$ dimensional manifold with a Riemannian
structure and $L \subset T^*N$ is an exact Lagrangian.

\begin{Proposition}
  Let $\F$ be either $\Q$ or $\F_p$ for some $p$ such that the
  requirements in Corollary~\ref{cor:2} is satisfied. Assume also that
  $\pi_1(L) \to \pi_1(N)$ is surjective. Then $L\to N$ is a homology
  equivalence with coefficients in $\F$.
\end{Proposition}

\begin{bevis}
  Under the assumptions in Corollary~\ref{cor:2} the induced map
  \begin{align*}
    (i_r)_* \colon H_*(L) \to H_*(\FL) \cong H_\bullet(\Lambda
    L,\F) 
  \end{align*}
  from Section~\ref{cha:inc-const} is an injection on the level of
  homology. Indeed, evaluation at base-point defines a splitting. In
  Lemma~\ref{lem:14as} we proved that this gives a map of Serre
  spectral sequences - one coming from fibrantly replacing $L\to N$
  and the one from $\FL$ from Proposition~\ref{Spectral}. More
  precisely this defines two (at this point possibly different)
  filtrations on $H_p(L)$, one coming from the Serre filtration over
  $N$ and one coming from applying $(i_r)_*(\alpha)$ and restricting
  the filtration associated to $\FL$. Denote these; the
  Serre-filtration $F^S_nH_*(L)$ and the $\FL$-filtration
  $F^{\FL}_nH_*(L)$. The fact that Lemma~\ref{lem:14as} says that we
  have a map of the spectral sequences proves that
  \begin{align*}
    F^S_n H_*(L) \subset F^{\FL}_nH_*(L).
  \end{align*}
  In particular
  \begin{align} \label{eq:4}
    H_n(L)\subset F^{\FL}_nH_*(L),  
  \end{align}
  which was not a priori clear because that spectral sequence
  associated to $\FL$ can be non-trivial in the 4. quadrant.

  Claim: $F^S_{n-1}H_n(L)=0$. Proof of claim:
  assume for contradiction that $0\neq \alpha \in H_n(L)$ is in
  Serre-filtration $n-1$. This implies by the above that it is also in
  $\FL$-filtration $n-1$. By assumption $\F$ is a field for which the
  intersection product on $H_*(L)$, which is Poincare dual (up to
  a sign) to the cup product, is a perfect parring. So there exists a
  $\beta \in H_{d-n}(L)$ such that $\beta\cdot\alpha = 1_L\in
  H_0(L)$. By Equation~\eqref{eq:4} $\beta$ is in filtration
  $d-n$. The construction of the product in Section~\ref{cha:serprod},
  which by Lemma~\ref{lem:4de} and Lemma~\ref{lem:14asd} realizes the
  intersection product, tell us that the intersection product
  preserves the $\FL$-filtration. This implies that $1_L$ is in
  filtration $-1$ in the $\FL$-filtration. This is a contradiction,
  since the entire spectral sequence is supported in the 1. and 4.
  quadrant.

  The above claim implies that the Serre spectral sequence associated
  to the map $L \to N$ has support entirely on the 1. axis in the
  abutment.

  The fact that $\pi_1(L) \to \pi_1(N)$ is surjective implies that the
  homotopy fiber of $L\to N$ is connected, and combining this with the
  fact that the $E_\infty$ term of the Serre spectral sequence for
  $L\to N$ is concentrated on the first axis we get that $H_*(L)$
  injects into $H_*(N)$. Since $L$ and $N$ are both $d$-dimensional
  manifolds this implies that the degree is invertible and thus the
  maps is also surjective on homology.
\end{bevis}

This now implies the main theorem rather easily.

\begin{MainTheorem}
  If $N$ is oriented and the induced map $p_*$ on fundamental
  groups is surjective then the induced map $p_*$ on homology is an
  isomorphism.
\end{MainTheorem}

\begin{mainbevis}
  Since the assumptions in Corollary~\ref{cor:2} is satisfied for
  $\F=\F_2$ the above proposition proves that $p_*$ is a $\F_2$-homology
  equivalence. It is a fact that any $\F_2$ cohomology equivalence of
  manifolds preserves the Stiefel-Whitney classes of the tangent
  bundles (this also mentioned in \cite{FSS} and used in
  \cite{Abou1} in the exact same context). This implies that $L$ is
  also oriented. It also implies that $p\colon L \to N$ is relative
  spin. This means that the assumptions in Corollary~\ref{cor:2} is
  satisfied for any $\F$.

  The above proposition thus works for all $\F=F_p$ and $\F=\Q$, and
  it is well-known that a map which is a homology equivalence for all
  these coefficients is a homology equivalence with $\Z$ coefficients.
\end{mainbevis}


\appendix
\chapter{Construction of Hamiltonians} \label{cha:pd-ham}

\section*{Asymptotic behavior and $\HS_\infty$.}

In the constructions throughout the paper we use functions $H_\infty
\colon T^*N \to \R$ with certain properties. We now construct the
space $\HS_\infty\subset \HS$ in which these should lie. We need this
space to satisfy several assumptions, and the easiest way to get this
is to explicitly construct it and then prove all of these.

So fix any $h\colon \R_{\geq 0} \to \R$, which satisfies
\begin{itemize}
\item {$h$ is convex,}
\item {$h(t)$ is zero for $t\in [0,2/3]$, and}
\item {$h(t) = \mu t+ c$ for $t\in [1,\infty[$ for some constants
    $\mu,-c \in \R_+$.}
\end{itemize}
For some $\beta>0$ small enough we define $H_\infty$ to be in
$\HS_\infty$ if
\begin{align*}
  H_\infty(q,p) = \epsilon h(\norm{p})
\end{align*}
for any $0<\epsilon<\beta$.

For $\beta>0$ small enough we have for any $H_\infty\in \HS_\infty$
that
\begin{itemize}
\item[H1)] {$H_\infty$ is 0 on $D_{2/3}T^*N$.}
\item[H2)] {$t H_\infty \in \HS_\infty$ for all $t \in ]0,1]$.}
\item[H3)] {Any time-1 Hamiltonian flow line for $H_\infty$ starting
    and ending in the same fiber $T^*_qN$ is constant.}
\item[H4)] {The finite dimensional approximation $S_1$ from
    Section~\ref{finred} is defined using the unique sub-division in
    $\Delta^0$ and any Hamiltonian in $\HS$ sufficiently $C^2$-close
    to $H_\infty$.}
\end{itemize}
Obviously H1) and H2) holds by construction. H3) follows because the
Hamiltonian flow of $H_\infty$ is a reparametrization of the geodesic
flow with constant speed $\norm{\nabla H_\infty}$. So if $\norm{\nabla
  H_\infty}\leq \delta_0$ (smaller than the injective radius) then it
can not return to the same starting point - unless it had speed 0, in
which case it is constant. H4) is a simple matter of making sure that
\begin{align*}
  C_1^H + C_2^H < \delta,
\end{align*}
with $\delta$ from Definition~\ref{def:Ardef}. 

\section*{Constructing $f_s$ from $f$.}

In this section we construct the smooth family of smooth increasing
functions $f_s\colon \R_{\geq 0} \to \R$, $s\in [0,\infty[$ used in
Section~\ref{theex} and later. Throughout the paper we need this to
satisfy several properties.

We assume that we are given a smooth function $f\colon \R_+ \to \R$
such that
\begin{itemize}
\item { $f(t) \to -\infty$ when $t\to 0$.}
\item { $f(t) = 0$ when $t\in [1,\infty]$.}
\item { $f''(t) < 0$ when $t\in ]0,1[$ - i.e. $f$ is strictly concave
    on $]0,1[$.}
\end{itemize}
For each $s>0$ we define $t_s$ to be the unique point $t_s\in ]0,1[$
such that the tangent of $f$ at $t_s$ intersects the $2^\textrm{nd}$
axis at $-s$ (see figure~\ref{fig:fs}).

The first properties we want (see figure~\ref{fig:fs}) for all $f_s$
is
\begin{itemize}
\item[f1)]{$f_s\colon \R_{\geq 0} \to \R$ is a smooth family for
    $s \in[0,\infty[$ of smooth \emph{increasing} functions.}
\item[f2)]{For $s>0$ and $t$ close to $0$ we have $f_s(t)=ct^2$ for
    some $c=c(s)>0$.}
\item[f3)]{For $t\geq 1$ we have $f_s(t) = s$.}
\item[f4)]{For $s>0$ the restriction of $f_s(t)$ to $]0,1[$ is
    strictly increasing.}
\end{itemize}
Of course these forces $f_0=0$. In the construction of $\FL$ we need
the following for $s\geq 5$ (here 5 is a rather arbitrary choice).
\begin{itemize}
\item[f5)]{The tangent to $f_s$ at any $t\in[0,t_s]$ intersects the
    $2^\textrm{nd}$ axis in $]-1,0]$.}
\item[f6)]{For $s\geq 5$ we have $f_s(t)=f(t)+s$ when $t\geq t_s$.}
\end{itemize}
At other points in the paper we will need for all $s\in
[0,\infty[$ that:
\begin{itemize}
\item[f7)]{The tangent to $f_s$ at any $t\in [0,1]$ intersects the 
    $2^\textrm{nd}$ axis above $-s/4$.}
\item[f8)]{Lemma~\ref{lem:closact} is true for all $s\in [0,\infty[$.}
\end{itemize}
The last one uses the definition of the Hamiltonians $H^s$ in
Section~\ref{theex}, and depends on the fixed choice of Riemannian
structure on $L$ from that section, which we will thus assume is
given.

\begin{figure}[ht]
  \centering
  \includegraphics{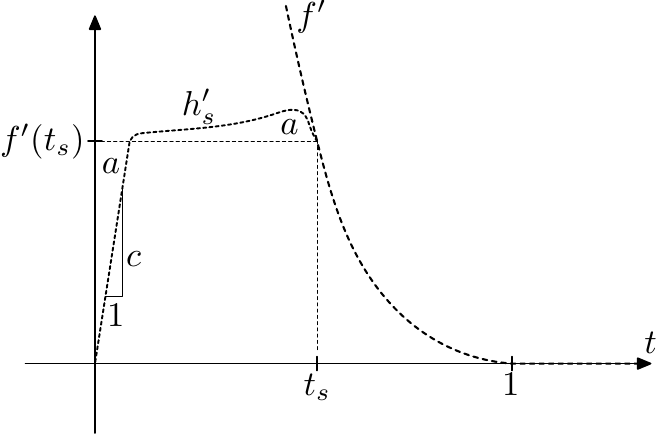}
  \caption{Constructing $h'_s$}
  \label{fig:bumpfpr}
\end{figure}

As Figure~\ref{fig:bumpfpr} suggests it is not difficult using bump
functions and cut-off functions together with $f'$ to construct a
smooth family $h_s ,s > 0$ satisfying f1-f7 - and even f6 for all
$s>0$. Indeed, pick $h_s$ as the unique anti-derivative to the pictured
$h'_s$ (which is equal to $f'$ for $t\geq t_s$) satisfying f2. We get
f4 trivially, and by making sure that the two areas marked with $a$ is
the same this will automatically satisfy f6 and thus also
f3. Properties f5 and f7 follow if we make $h'_s$ have a unique
critical point (a maximum) close to $t_s$ - as suggested in the
figure. Indeed this controls the inflection point of $h_s$. We can
make this maximum value as close to $f'(t_s)$ as we would like by
increasing $c$ and thus decreasing the areas $a$. We may thus assume
that
\begin{itemize}
\item {$h'_s$ tends to 0 when $s$ tends to 0.} 
\end{itemize}
The proof of Lemma~\ref{lem:closact} then works as written for
this family for all $s>0$. The problem is having f6 satisfied close to
0, and retain smoothness ($C^2$ would be enough, but this seems no
easier). Indeed, in the paper we used smoothness at $s=0$ to get
bounds on the first and second derivatives of $f_s$ for $s$ close to
$0$ - so this is crucial.

We thus use that $h'_s$ tends to zero when $s$ does. Indeed, pick some
$s_0>0$ small enough such that $h'_s(t)$ is less than the shortest
non-zero geodesic length of $L$ for all $t$ and all
$s\in[0,s_0]$. Then the family
\begin{align*}
  \tfrac{s}{s_0} h_{s_0} \colon \R \to \R  
\end{align*}
satisfies all the requirements for $s<5$ (f6 is now considered void
because $s<5$) because $h_{s_0}$ does so. All properties f2-f7 are for
fixed $s$ preserved by convex combinations. So we may choose any
smooth function $\psi \colon \R_{\geq 0} \to [0,1]$ with $\psi(s)=0$
for $s\leq s_0/2$ and $\psi(s)=1$ for $s\geq s_0$ and define
\begin{align*}
  f_s(t) = (1-\psi(s))\tfrac{s}{s_0} h_{s_0}(t) + \psi(s) h_s(t).
\end{align*}
This then satisfies all the requirements. Indeed, it is smooth (f1)
and f2-f7 is handled by the convex combination argument. Even f8 is
true simply because $f_s$ for $s\in [0,s_0]$ has $f'_s$ bounded by the
length of any non-zero geodesic on $L$, and so only the constant
geodesics contributes, and by f4 these only happens for $\norm{p}_L=0$
or $\norm{p_L}\geq 1/2$, which were already considered in the proof
of Lemma~\ref{actions}.


\chapter{A Functor from $\beta$-spectra to Parametrized
  Spectra}\label{cha:funcrep}

In Section~\ref{cha:gen-spec} we defined $\beta$-spectra and their
homology. In this Appendix we describe a functor $F$ from
$\beta$-spectra over $W$ to $l$-spectra,  were $l\in \N$ denotes the
trivial bundle $W\times \R^l$. This functor preserves the homology in
the cases that we care about, and we used it to simplify the
construction of the spectral sequence in Section~\ref{serre}.

Let $i \colon \beta \to W \times \R^l$ be an embedding of the bundle
$\beta$ over $W$. Then we get a canonical isomorphism
$\beta^\perp \oplus \beta \cong l=W\times \R^l$, where $\beta^\perp$ is
the orthogonal complement bundle of $\beta$ in $W\times \R^l$. We then
define the functor $F$ on an object $\As\in \Sp_W^\beta$ by
\begin{align*}
  (F\As)_r = \As_r^{r\beta^\perp},
\end{align*}
The new structure maps are defined by
\begin{align*}
  (F\sigma_r) = \sigma_r^{\beta^\perp} \colon 
  \As_r^{r\beta^\perp \oplus l} \cong
  \As_r^{r\beta^\perp \oplus \beta^\perp \oplus \beta} \cong 
  \As_r^{\beta\oplus (r+1)\beta^\perp} \to
  \As_{r+1}^{(r+1)\beta^\perp}.
\end{align*}
This is easily extended to morphisms since the fiber-wise Thom
suspension is functorial on ex-spaces as described in
Section~\ref{cha:gen-spec}.

Global homology is easily seen to be preserved in all cases of
interest. Indeed, we made sure that the sections are cofibrations and
that $\beta$ (and hence $\beta^\perp$) is oriented. We could also
argue that for a 3S fibrant $\beta$-spectrum also the homologies of
the fibers are preserved, i.e. the non-genericity of the fibers goes
away. However, we will not actually need this since we never use what
exactly the homology of the fibers are in any of the arguments.

\begin{Lemma}
  \label{lem:9}
  The functor $F$ preserves 3S fibrancy.
\end{Lemma}

\begin{bevis}
  This is an easy consequence of the fact that any stable lift for
  $\As$ can by using parallel transport in (copies of) $\beta^\perp$
  be turned into a stable lift for $F(\As)$. 
\end{bevis}


\chapter{Homotopy Invariance of $\FL$}\label{cha:hominv}

In this appendix we discus why the construction of $\FL$ is unique up
to unique isomorphism in the homotopy category of parametrized spectra
over $N$. A similar discussion can be done for the product, but we
omit this.

Following notation from Section~\ref{theex} we assume for the moment
that $f$ and the capping off family $f_s$ is fixed. The extra choices
of $(a,b,(s_r,q_r,\sigma_r)_{r\in \N})$ in defining an instance of $\FL$
is contractible - that is \emph{if} we forget the assumption that $a$
should be a regular value for each $S_r^{s_r}$. We only assumed this
so that the section would be cofibrations making subsequent
calculations easier - so we will simply not assume this here. The
argument in Section~\ref{cha:fib-fl} that $\FL$ is 3S fibrant can be
extended to any compact smooth family of such a set of choices. That
is if $K$ is a compact subset of a smooth manifold which smoothly
parametrizes choices as above then we get a 3S fibrant $TN$-spectrum
over $N\times K$. So in fact the object in the homotopy category does
not depend on these choices up to unique isomorphism.

This can be extended to families of $(f,f_s)$ if we ask that the
family of $f$'s agree close to $0$, and that the same is true for
$f_s$ for large $s>0$ (on the same given neighborhood of $0$). Indeed,
with this we can see that the differences in the functions is pushed
away from what happens a critical values $a$ if $s$ is large enough,
and so the fibrancy still holds at least for large $s$.

For two different choices $(f^i,f^i_s),i=1,2$ which does not
agree like above one can create interpolations and construct a
zig-zagging argument. Indeed, we can construct a pair $(F,F_s)$ such
that away from a small neighborhood of $0$ $F=f^1$ but on an even
smaller neighborhood we have $F=f^2$, we may thus construct the
sequence $F_s$ such that for small $s$ it looks like $f_s^1$ and for
large $s$ it looks like $f^2_s$ close to $0$. Sequences of these which
agrees with $f^1$ on smaller and smaller sets can be used to construct a
map from the instance defined by $(f^1,f_s^1)$ to the one defined by
$(f^2,f_s^2)$. To proved that this is an equivalence we create
interpolations from $(F,F_s)$ to $(f^1,f_s^1)$, which on smaller an
smaller sets agree with $(F,F_s)$, but which looks like $(f^2,f_s^2)$
on a set $(\epsilon,\epsilon')$ (here one needs to choose $\epsilon$
dependent on $\epsilon'$ for the last part to work) and for some range
of $s$ such that we can factor the identity on the instance defined by
$(f^1,f_s^1)$ through the map defined by the first sequence.


\chapter{Coherent Orientations and Finite Dimensional
  Approximations}\label{cha:coor}

To understand the issue in Remark~\ref{rem:coef} better we will in
this appendix describe why coherent orientations can be subtle from
the point of view of finite dimensional approximations. We will also
relate this to stable homotopy types and discuss why in some sense
Viterbo functoriality is not natural for the stable homotopy types -
unless extra structure is considered as in \cite{hejeh} where we
realize it as a map of spectra. 

Consider the following abstract situation: assume that $i\colon M \to
M'$ is a finite dimensional manifold $M$ embedded into an infinite
dimensional manifold $M'$. Also assume that we have a Morse-Smale
function $f\colon M' \to \R$ for which some sort of infinite
dimensional Morse homology with $\Z$ coefficients can be
defined - even though Morse indices happens to be
infinite. Furthermore, assume that all critical points of $f$ and all
flow lines between such (defined by a PDE and not a flow) are
contained in the image of $i$ as an isolated invariant set $S \subset
M$ (see \cite{MR511133} for definition of isolated invariant set) for
$f\circ i$. To define this we use some pseudo-gradient for $f\circ
i$.

In such a situation we cannot conclude that the homology of the
Conley index of $S$ equals the Morse homology defined for $f$. Indeed,
to define Morse homology of $f$ with $\Z$ coefficients we need some sort
choice of coherent orientations that depends on more than just the
isolated invariant set. To see this take any finite
dimensional vector bundle $E\to M$, whose restriction to $S$ is not
orientable. Then because the Morse indices of
$f$ are infinite there will often (if not always) be ``room'' to do
the following; we may extend $i$ to an embedding  $i'\colon E \to M'$
in such a way that $f \circ i'$ is strictly concave on the fibers with
unique maximum on the zero-section. With a little precision one can
prove that the Conley index of $f \circ i'$ is the Thom space of $E$
on the Conley index of $f\circ i$ (relative to the base point). This
means that the homologies are not the same because $E$ is not
orientable on $S$.

The above discussion seems to indicate that it is difficult to capture
the ``correct'' orientation on the homology theories when doing finite
dimensional approximation. However, in cotangent bundles there is a
canonical choice of Lagrangian foliations which from the stable
homotopy theory point of view settles the above problem in more
generality than just $\Z$ coefficients. Indeed, one may heuristically
think of this choice of Lagrangian at each point in $T^*N$ as
canonically identifying a subset of the finite dimensional
approximations as above, which are only related by an $E$ (also as
above) if this $E$ is stably parallizable. So in fact there is a
canonical way of making a spectrum. However, it turns out that the
orientation coming from this choice is not the one usual used in Floer
theory, and, indeed, this is why the symplectic cohomology does not
agree with the homology of this spectrum. Even more, when one
considers Viterbo functoriality from $T^*N$ to $DT^*L$ the Lagrangian
foliations induced by the cotangent bundle structures do not agree, and
the different choices induce the bundle $\eta$ in
Equation~\eqref{eq:8sdafuc}.


\chapter{Nearby Lagrangians are homotopy equivalent (by
  M. Abouzaid)}\label{cha:nlahe}

As in the introduction, let $N$ be a closed smooth manifold, $L$ a closed
exact Lagrangian in $\TN$, $\tN$ the universal cover of $N$, and $\tL$ the
inverse image of $L$ in $\TtN$.  In this short appendix, we prove the
following result
\begin{Proposition}\label{prop:main}
If  $\tL$ has vanishing Maslov class, then so does  $L$.
\end{Proposition}
As a corollary, we conclude
\begin{Thm} \label{thm:main}
If $L$ is a closed exact Lagrangian in $\TN$, then $L$ is homotopy
equivalent to $N$.
\end{Thm}
\begin{bevis}
  By Corollary~\ref{cor:1gfr}, every closed exact Lagrangian satisfies
  the hypothesis of Proposition \ref{prop:main}, so $L$ has vanishing
  Maslov index. The result therefore follows from \cite{Abou1} which
  proves that every closed exact Lagrangian in $\TN$ with trivial
  Maslov class is homotopy equivalent to the base.
\end{bevis}
The starting point for the proof of Proposition \ref{prop:main} is the
main result  of \cite{abou2} that the wrapped Fukaya category of
$\TN$ is split-generated by a fibre $\Tn$.  In particular, the Floer
cohomology of $L$ with $\TN$ does not vanish.  To be more precise,
\cite[Theorem 1.1]{abou2} shows that the Fukaya category consisting
of Lagrangians of vanishing Maslov index is generated by a fibre.  As
an intermediate step, one prove a split-generation result, and that
proof applies, without modification, to the more general case in which
the Maslov index may not necessarily vanish.  To see this, we use the
fact that the proof of Theorem 1.1. of \cite{abou2} itself relies
on Theorem 1.1  of \cite{Abou0}, which asserts that the Fukaya
category of a symplectic manifold is split-generated by a Lagrangian
$L$ whenever a certain geometric criterion involving the Hochschild
homology of the Floer complex $CF^{*}(L,L)$ is satisfied.  The fact
that the cotangent fibre in a cotangent bundle satisfies this property
is verified in Proposition 1.6 of \cite{abou2}. Once this criterion
is satisfied, the annulus degeneration argument described in
\cite[Section 6]{Abou0} shows that the self-Floer cohomology of any
Lagrangian in $\TN$ vanishes if and only if its Floer cohomology with
a fibre vanishes.  This result in no way uses integral gradings in
Floer cohomology.  If $L$ is an exact Lagrangian, its self-Floer
cohomology is isomorphic to its ordinary cohomology, so we conclude:
\begin{Proposition}\label{prop:non-zero-HF}
If $L$ is an exact Lagrangian in $\TN$, then the Floer cohomology
$HF(L,\TN)$ does not vanish. 
\end{Proposition}

In order to derive Proposition \ref{prop:main}, we need to discuss the
construction of gradings in Lagrangian Floer cohomology.  Starting with a
Riemannian metric on any smooth manifold $N$ (not necessarily closed), we
obtain a density which assigns to a basis of tangent vectors at a point the
square of the volume of the corresponding parallelepiped. Whenever $N$ is
oriented, this density is the square of the volume form determined by the
metric and the orientation.  By complexifying this density, we obtain a
\emph{quadratic complex volume form} $\eta$ on $\TN$; in local volume
preserving coordinates $(q_1, \ldots, q_n)$ with dual coordinates $(p_1,
\ldots, p_n)$ we have the expression
\begin{align}  \label{eq:quadratic_volume_equation}
  \eta = \left( \left( dp_1 + \sqrt{-1}dq_1 \right) \wedge \cdots \wedge 
    \left(dp_n + \sqrt{-1} dq_n \right) \right)^{\otimes 2}
\end{align}
where the symbol $\otimes 2$ means that we assign to an $n$-tuple of
tangent vectors on $\TN$ the square of the complex number obtained by
applying the complex-valued $n$-form written in coordinates inside the
parentheses.

Given any  Lagrangian $L \subset \TN$, we may evaluate $\eta$ on a basis
of tangent vectors to obtain a complex number.  The Lagrangian condition
implies that this number does not vanish, so we obtain a \emph{complex
phase map}
\begin{align}
  \label{eq:phase_map}
  \theta_{L} \co L & \to \bR/ 2\pi \bZ\\
  e^{ \sqrt{-1} \theta_{L} (x)}  & =  \frac{\eta(v_1 \wedge \cdots \wedge
    v_n)}{| \eta(v_1 \wedge \cdots \wedge v_n) |},
\end{align}
where $(v_1, \ldots, v_n)$ is an arbitrary basis of tangent vectors for
$L$ at $x$.
\begin{Definition}
  The \emph{Maslov class} of $L$ is the integral first cohomology class
  represented by $\theta_{L}$.  If $L$ has vanishing Maslov class then a
  \emph{grading} is a choice $\tilde{\theta}_{L}$ of an $\bR$-valued lift of
  $\theta_{L}$:
  \begin{align}
    \xymatrix{L \ar[r] \ar[dr] &  \bR \ar[d] \\ & S^1 = \bR /2 \pi \bZ.} 
  \end{align}
\end{Definition}

Vanishing of the Maslov class is the condition needed  in order to be able
to assign an integer (up to global shift) to each intersection point
between a pair of Lagrangians, and the choice of a grading determines such
an integer and hence equips Lagrangian Floer complexes and cohomology 
groups with  $\bZ$-gradings.  We shall only need the following case: 
Assume that we are given a Lagrangian $L$ which intersect $\Tn$
transversely as a point $x$.  Since $L$ is Lagrangian, one may choose local
coordinates $(q_1, \ldots, q_n)$ on $N$ so that the tangent space of $L$ is
spanned by products of $n$ lines, each lying in a plane spanned by $\partial_{p_i}$ and $
\partial_{q_i}$, with phase $\alpha_{i}(x) \in [0,\pi)$.

Note that $\theta_{\Tn}$ is independent of the point on $\Tn$, so that we
can choose $\tilde{\theta}_{\Tn}$ to be identically $0$.  In this case, the
following formula for the Maslov index is given in Section 3.2 of
\cite{MR1957663}:
\begin{align} \label{eq:thomas-yau-formula}
 \mu(x) = \frac{1}{\pi} \left(  \sum_{i=1}^{n} \alpha_{i}(x) -
\theta_{L}(x) \right).
\end{align}

With this formula at hand, we can finish the proof of the result announced
at the beginning of this note. 

\begin{bevisE1}
  Assume that $\tL$ has vanishing Maslov class.  In \cite{Abou1}, a
  wrapped Fukaya category was assigned to the universal cover of $\TN$, which
  is also the cotangent bundle $\TtN$ of the universal cover of $N$.  By
  comparing holomorphic discs in $\TN$ and $\TtN$, we find that the Floer
  cohomology of $L$ with a fibre $\Tn$ is isomorphic to that of $\tL$ with
  any cotangent fibre $\Ttn$, and in particular has finite rank. Upon
  choosing a grading on $\tL$, we conclude from Proposition \ref{prop:non-zero-HF} that
  \begin{equation}
    \parbox{30em}{$HF^{*}(\tL, \Ttn)$  is non-zero, and is supported in finitely many
      cohomological degrees.}
  \end{equation}
  Since all cotangent fibres are Hamiltonian isotopic via an isotopy which
  preserves the vanishing of the phase, and Floer cohomology together
  with its grading is invariant under such isotopies, the group
  $HF^{*}(\tL, \Ttn)$ is moreover independent of $q'$.

  We shall use Equation \eqref{eq:thomas-yau-formula} to derive a
  contradiction to the bounded support of $ HF^{*}(\tL, \Ttn) $  if $L$ does
  not also have vanishing Maslov class.

  Choose a quadratic volume form on $\TtN$ which is pulled back from $\TN$. 
  In particular,  a grading $  \tilde{\theta}_{\tL}  $ on $\tL$ fits in a
  commutative diagram
  \begin{equation}
    \xymatrix{  \tL \ar[r]^{\tilde{\theta}_{\tL}} \ar[d] & \bR \ar[d]  \\  
      L \ar[r]^{\theta_{L}} & S^1.}
  \end{equation}
  Given any point $\tx \in \tL$, the commutativity of the above diagram
  implies that for any loop $\gamma \in \pi_{1}(N,q)$, with associated action
  on $\tN$ and $\TtN$ by deck transformations, we have
  \begin{equation}
    \theta_{L}(\gamma \cdot \tx) =  \theta_{L}(\tx)  +
    [\theta_{L}]([\gamma])
  \end{equation}
  where $[\gamma]  $ is the class represented by $\gamma$ in $H_{1}(Q) \cong
  H_{1}(L)$, and $ [\theta_{L}] $  is the Maslov class of $L$.

  Let us now assume that $L$ intersects $\Tn$ transversely, which may be
  achieved after a small perturbation.   Applying Equation
  \eqref{eq:thomas-yau-formula}, we conclude that
  \begin{equation}
    \mu(\gamma \cdot \tx) =   \mu(\tx) -  [\theta_{L}]([\gamma]).
  \end{equation}
  Since $\gamma \cdot (\tL \cap \Ttn) = \tL \cap T^{*}_{\gamma \cdot \tq}
  \tN $, with a corresponding identification of holomorphic curves for
  appropriate complex structures, we conclude that
  \begin{equation}
    HF^{*}( T^{*}_{\gamma \cdot \tq} \tN   , \tL) = HF^{* - 
      [\theta_{L}]([\gamma])}( \Ttn , \tL),
  \end{equation}
  which, because of the finite dimensionality and invariance of this
  group, implies that it vanishes, contradicting the output of
  Proposition \ref{prop:non-zero-HF}. 
\end{bevisE1}


\bibliographystyle{plain}
\bibliography{../Mybib}

\def\cprime{$'$}
\begin{thebibliography}{10}

\bibitem{MR2190223}
Alberto Abbondandolo and Matthias Schwarz.
\newblock On the {F}loer homology of cotangent bundles.
\newblock {\em Comm. Pure Appl. Math.}, 59(2):254--316, 2006.

\bibitem{abou2}
Mohammed Abouzaid.
\newblock {A cotangent fibre generates the Fukaya category}.
\newblock {\em arXiv:1003.4449v5}, 2010.

\bibitem{Abou0}
Mohammed Abouzaid.
\newblock A geometric criterion for generating the {F}ukaya category.
\newblock {\em Publ. Math. Inst. Hautes \'Etudes Sci.}, (112):191--240, 2010.

\bibitem{Abou1}
Mohammed Abouzaid.
\newblock {Maslov 0 nearby Lagrangians are homotopy equivalent}.
\newblock {\em arXiv:1005.0358v1}, 2010.

\bibitem{MR765426}
Marc Chaperon.
\newblock Une id\'ee du type ``g\'eod\'esiques bris\'ees'' pour les syst\`emes
  hamiltoniens.
\newblock {\em C. R. Acad. Sci. Paris S\'er. I Math.}, 298(13):293--296, 1984.

\bibitem{MR2039760}
Ralph~L. Cohen, John D.~S. Jones, and Jun Yan.
\newblock The loop homology algebra of spheres and projective spaces.
\newblock In {\em Categorical decomposition techniques in algebraic topology
  ({I}sle of {S}kye, 2001)}, volume 215 of {\em Progr. Math.}, pages 77--92.
  Birkh\"auser, Basel, 2004.

\bibitem{MR511133}
Charles Conley.
\newblock {\em Isolated invariant sets and the {M}orse index}, volume~38 of
  {\em CBMS Regional Conference Series in Mathematics}.
\newblock American Mathematical Society, Providence, R.I., 1978.

\bibitem{MR0192475}
D.~B.~A. Epstein.
\newblock The degree of a map.
\newblock {\em Proc. London Math. Soc. (3)}, 16:369--383, 1966.

\bibitem{FSS}
Kenji Fukaya, Paul Seidel, and Ivan Smith.
\newblock Exact {L}agrangian submanifolds in simply-connected cotangent
  bundles.
\newblock {\em arXiv:math/0701783v1}, 2007.

\bibitem{MR1867354}
Allen Hatcher.
\newblock {\em Algebraic topology}.
\newblock Cambridge University Press, Cambridge, 2002.

\bibitem{HatchSpec}
Allen Hatcher.
\newblock {\em Spectral Sequences in Algebraic Topology}.
\newblock \url{www.math.cornell.edu/~hatcher/SSAT/SSATpage.html}, 2004.

\bibitem{hejeh}
Thomas Kragh.
\newblock The viterbo transfer as a map of spectra.
\newblock {\em ArXiv:0712.2533v1}, 2009.

\bibitem{MR1090163}
Fran{\c{c}}ois Lalonde and Jean-Claude Sikorav.
\newblock Sous-vari\'et\'es lagrangiennes et lagrangiennes exactes des fibr\'es
  cotangents.
\newblock {\em Comment. Math. Helv.}, 66(1):18--33, 1991.

\bibitem{MR2271789}
J.~P. May and J.~Sigurdsson.
\newblock {\em Parametrized homotopy theory}, volume 132 of {\em Mathematical
  Surveys and Monographs}.
\newblock American Mathematical Society, Providence, RI, 2006.

\bibitem{MR1793722}
John McCleary.
\newblock {\em A user's guide to spectral sequences}, volume~58 of {\em
  Cambridge Studies in Advanced Mathematics}.
\newblock Cambridge University Press, Cambridge, second edition, 2001.

\bibitem{MR0163331}
J.~Milnor.
\newblock {\em Morse theory}.
\newblock Based on lecture notes by M. Spivak and R. Wells. Annals of
  Mathematics Studies, No. 51. Princeton University Press, Princeton, N.J.,
  1963.

\bibitem{Nadler}
David Nadler.
\newblock {Microlocal Branes are Constructible Sheaves}.
\newblock {\em arXiv:math/0612399v4}, 2007.

\bibitem{MR797044}
Dietmar Salamon.
\newblock Connected simple systems and the {C}onley index of isolated invariant
  sets.
\newblock {\em Trans. Amer. Math. Soc.}, 291(1):1--41, 1985.

\bibitem{Salweb}
Dietmar~A. Salamon and Joa Weber.
\newblock Floer homology and the heat flow.
\newblock {\em GEOM. FUNCT. ANAL}, 16:1050--1138, 2005.

\bibitem{MR1957663}
R.~P. Thomas and S.-T. Yau.
\newblock Special {L}agrangians, stable bundles and mean curvature flow.
\newblock {\em Comm. Anal. Geom.}, 10(5):1075--1113, 2002.

\bibitem{MR1403954}
Claude Viterbo.
\newblock Generating functions, symplectic geometry, and applications.
\newblock In {\em Proceedings of the {I}nternational {C}ongress of
  {M}athematicians, {V}ol.\ 1, 2 ({Z}\"urich, 1994)}, pages 537--547, Basel,
  1995. Birkh\"auser.

\bibitem{MR1617648}
Claude Viterbo.
\newblock Exact {L}agrange submanifolds, periodic orbits and the cohomology of
  free loop spaces.
\newblock {\em J. Differential Geom.}, 47(3):420--468, 1997.

\end{thebibliography}

\end{document}